\documentclass[11 pt, reqno]{article}

\usepackage{amsmath,amsfonts,amssymb,amsthm}
\usepackage[colorlinks=true,hyperindex=true]{hyperref}
\usepackage{cancel,bbm}
\usepackage{mathabx} 
\usepackage{comment}
\usepackage{enumitem}
\usepackage{cite}

\usepackage{chngcntr}
\counterwithin*{equation}{section}

\usepackage[a4paper, left=2cm, right=2cm, top=2 cm, bottom= 2 cm]{geometry}
\usepackage{color}
\usepackage{fancyhdr}
\usepackage{latexsym}

\newtheorem{ccounter}{ccounter}[section]
\newtheorem{thm}[ccounter]{Theorem}
\newtheorem{lem}[ccounter]{Lemma}
\newtheorem{cor}[ccounter]{Corollary}
\newtheorem{defn}[ccounter]{Definition}
\newtheorem{prop}[ccounter]{Proposition}
\newtheorem{ass}[ccounter]{Assumption}
\newtheorem{ex}[ccounter]{Example}

\def\bet{\begin{thm}}
\def\eet{\end{thm}}
\def\bel{\begin{lem}}
\def\eel{\end{lem}}
\def\bas{\begin{ass}}
\def\eas{\end{ass}}
\def\bec{\begin{cor}}
\def\eec{\end{cor}}
\def\bed{\begin{defn}}
\def\eed{\end{defn}}
\def\bep{\begin{prop}}
\def\eep{\end{prop}}
\def\beq{\begin{equation}}
\def\eeq{\end{equation}}
\def\proof{\noindent {\bf Proof.}\ \ }
\def\bea{\begin{equation*}}
\def\eea{\end{equation*}}

\def\bex{\begin{ex}}
\def\eex{\end{ex}}
\def\remark{\noindent{\bf Remark. }}

\def\rr{\mathbb{R}}
\def\zz{\mathbb{Z}}
\def\cc{\mathbb{C}}
\def\1{\boldsymbol{1}}
\def\Im{\mathrm{Im}}
\def\Re{\mathrm{Re}}
\def\e{\mathrm{e}}
\def\i{\mathrm{i}}
\def\del{\partial}
\def\d{\mathrm{d}}
\def\eps{\varepsilon}
\renewcommand\leq\varleq
\renewcommand\geq\vargeq
\def\ee{\mathrm{E}}

\def\F{\mathcal{F}}
\def\O{\mathcal{O}}

\def\ee{\mathbb{E}}
\def\om{\omega}
\def\pp{\mathbb{P}}

\def\etast{\eta_*}
\def\phist{\phi_*}
\def\supp{\mathrm{supp}}

\def\mfct{m_{\mathrm{fc}, t}}
\def\rhofct{\rho_{\mathrm{fc}, t}}

\def\hatm{\hat{m}}
\def\hatrho{\hat{\rho}}
\def\hatmt{\hat{m}_{\mathrm{fc}, t}}
\def\hatE{\hat{E}}
\def\hatrhot{\hat{\rho}_{\mathrm{fc}, t}}
\def\hatD{\hat{\mathcal{D}}}
\def\muo{\mu_1}
\def\aom{a_{-, 1}}
\def\mut{\mu_2}
\def\atm{a_{-, 2}}
\def\Fo{F_1}
\def\Ft{F_2}
\def\xiom{\xi_{-, 1}}
\def\xitm{\xi_{-, 2}}
\def\Eom{E_{-, 1}}
\def\Etm{E_{-, 2}}
\def\bo{b_1}
\def\bt{b_2}
\def\ao{a_1}
\def\at{a_2}
\def\xio{\xi_1}
\def\xit{\xi_2}
\def\Eo{E_1}
\def\Et{E_2}
\def\I{\mathcal{I}}
\def\J{\mathcal{J}}

\def\tilz{\tilde{z}}

\def\hatz{\hat{z}}

\def\sumAi{\sum^{\A, (i) }}
\def\sumAci{\sum^{\A^c, (i)}}
\def\B{\mathcal{B}}
\def\V{\mathcal{V}}
\def\gamhat{\hat{\gamma}}
\def\hatgam{\hat{\gamma}}
\def\zetaz{\zeta^{(0)}}

\def\A{\mathcal{A}}

\def\D{\mathcal{D}}

\def\J{\mathcal{J}}
\def\L{\mathcal{L}}
\def\UL{\mathcal{U}^{\mathcal{L}}}

\begin{document}
\title{Deformed GOE}


\begin{table}
\centering

\begin{tabular}{c}
\multicolumn{1}{c}{\Large{\bf Edge statistics of Dyson Brownian motion}}\\
\\
\\
\end{tabular}
\begin{tabular}{c c c}
Benjamin Landon & & Horng-Tzer Yau\\
\\
 \multicolumn{3}{c}{ \small{Department of Mathematics} } \\
 \multicolumn{3}{c}{ \small{Harvard University} } \\
\small{landon@math.harvard.edu} &  & \small{htyau@math.harvard.edu}  \\
\\
\end{tabular}
\\
\begin{tabular}{c}
\multicolumn{1}{c}{\today}\\
\\
\end{tabular}

\begin{tabular}{p{15 cm}}
\small{{\bf Abstract:} We consider the edge statistics of Dyson Brownian motion with deterministic initial data.  Our main result states that if the initial data has a spectral edge with rough square root behavior down to a scale $\eta_* \geq N^{-2/3}$ and no outliers, then after times $t \gg \sqrt{ \eta_*}$, the statistics at the spectral edge agree with the GOE/GUE.  In particular we obtain the optimal time to equilibrium at the edge $t = N^{\eps} / N^{1/3}$ for sufficiently regular initial data.  
Our methods rely on eigenvalue rigidity results similar to those appearing in \cite{kevin1,kevin2}, the coupling idea of \cite{homogenization} and the energy estimate of \cite{bourgade2014edge}.   }
\end{tabular}
\end{table}
{\let\thefootnote\relax\footnote{The work of H.-T. Y. is partially supported by NSF Grant DMS-1307444,  DMS-1606305 and a Simons Investigator award.}}
\section{Introduction}
One of the guiding principles of random matrix theory is that of universality.  This states that the limiting behavior of the eigenvalues of large random matrices is universal.  Said differently, universality is the observation that often the eigenvalue fluctuations of two seemingly unrelated random matrix ensembles converge, as the size of the matrices $N \to \infty$, to the same limiting object.

Perhaps the most fundamental class of random matrix ensembles are Wigner matrices.  These are $N \times N$ symmetric matrices $W$ with independent (up to the constraint $W=W^*$) centered entries of identical variances.  We consider two symmetry classes of Wigner matrices; real symmetric, in the case that the entries are real; or complex Hermitian in the case that the entries are complex.  The prototypical examples of real symmetric and complex Hermitian Wigner matrices are the Gaussian Orthogonal and Unitary ensembles (GOE/GUE).  These are constructed by taking the entries to be standard real or complex Gaussians, respectively.  In the case of the GOE/GUE, the limiting eigenvalue behavior can be computed explicitly.  The universality conjecture for Wigner matrices can be rephrased as the fact that these formulas are true, in the limit $N \to \infty$, for all Wigner matrices regardless of the details of the distribution of the matrix elements, and depend only on the symmetry class (real symmetric or complex Hermitian) of the ensemble.

There has been significant progress on the understanding of bulk universality for random matrix ensembles. Bulk universality refers to the behavior of eigenvalues contained in the interior of the support of the eigenvalue density.  Bulk universality for Wigner matrices of all symmetry classes was proven in the works \cite{renyione, erdos2012spectral,erdos2010bulk,erdos2011universality,Gap,erdos2012rigidity}.  Parallel results were established in certain cases in \cite{tao2010random,tao2011random}, with the key result being a ``four moment comparison theorem."  

One of the major contributions of the works  \cite{renyione, erdos2012spectral,erdos2010bulk,erdos2011universality,Gap,erdos2012rigidity}  was to establish a general, robust framework within which to establish bulk universality for random matrix ensembles.  This {\it three-step strategy} is as follows.
\begin{enumerate}
\item Prove a local law for the random matrix ensemble under consideration.  This local law is used to establish high probability {\it rigidity} estimate for the eigenvalue locations.
\item Given a random matrix ensemble $H$, establish bulk universality for a Gaussian divisible ensemble of the form $H + \sqrt{t} G$, where $G$ is a GOE/GUE matrix, and $t = o(1)$ is interpreted as time.
\item A density argument comparing the eigenvalue statistics of a general ensemble to a Gaussian divisible ensemble for which universality was established in the previous step.
\end{enumerate}
The first step is model-dependent, and gives strong a-priori estimates crucial to the next two steps.  The third step is a perturbation argument, typically based on comparing the expectation of the Green's functions, or a use of the Ito lemma.  

For the present article, the second step is most relevent.  It is based on an analysis of the local ergodicity of {\it Dyson Brownian motion}.  In the seminal work \cite{dyson}, Dyson introduced a matrix-valued stochastic process and calculated the resulting eigenvalue dynamics.  At each fixed time $t$, the eigenvalues of this matrix-valued Brownian motion are distributed as the eigenvalues of a Gaussian divisible ensemble.  Crucially,  he found that the eigenvalues satisfy a closed system of stochastic differential equations, for which the GOE/GUE is the equilibrium measure.

While global eigenvalue statistics reach equilibrium under DBM only for long times $t \gg 1$, Dyson conjectured that local statistics reach equilibrium after a much shorter time.  The works \cite{localrelaxation,erdos2011universality,erdos2012bulk} established that if Dyson Brownian motion is started from a Wigner matrix, then the bulk GOE statistics are reached in a very short time $t \sim N^{-1}$.  This is one of the crucial elements in proving bulk universality for Wigner matrices.  Moreover, it provides a crucial conceptual understanding of the origin of universality - that of the local ergodicity of Dyson Brownian motion.

More recently, there has been much success in applying the above three-step strategy to the bulk statistics of random matrix ensembles which go beyond the Wigner class.  As in the case of Wigner matrices, one of the important ingredients has been a study of the behavior of Dyson Brownian motion for more general initial data \cite{ES, landonyau,fixed}.  The main contribution of these works is to establish bulk universality for the gap statistics \cite{ES,landonyau} and correlation functions \cite{fixed} of Dyson Brownian motion with very general deterministic or random initial data, going far beyond the class of Wigner matrices.  In terms of the three-step strategy outlined above and in the context of obtaining bulk universality, this takes care of the second step for many random matrix ensembles.

While the third step is quite robust and applies in most cases, obtaining a local law is typically model dependent and challenging.  We mention here some of the recent works on local laws and bulk universality for random matrix ensembles going beyond the Wigner class.  The adjacency matrices of $d$-regular graphs were considered in \cite{roland1,roland2}, and other sparse random graph ensembles in \cite{hl,huang,erdos2012spectral,renyione,ziliang}.   A very general class of Wigner-type matrices were studied in a series of papers \cite{ajanki1,ajanki2,ajanki3}, and matrices with correlated entries were studied in \cite{ziliang2,ajankicor, ajankigaus}.  An additive model of random matrices was studed in \cite{bao1,bao2,bao3,chel}.  The proof of the local ergodicity of Dyson Brownian motion \cite{ES,landonyau,fixed} for general initial data has been an important element in obtaining the eigenvalue universality in many of these works.

So far, the above discussion has been limited to the behaviour of the bulk eigenvalues, or those confined to the interior of the spectrum.  A natural question is also to investigate the behavior of the extremal eigenvalues of random matrices.  In the case of the Gaussian ensembles, Tracy and Widom \cite{tw1,tw2} calculated 
\beq
\lim_{N \to \infty} \pp \left[ N^{2/3} ( \lambda_N -2 ) \leq s \right] = F_\beta(s)
\eeq
where $\lambda_N$ is the largest eigenvalue of a GOE/GUE matrix and $F_\beta$ can be characterized in terms of Painlev\'e equations and $\beta=1, 2$ respectively for the real symmetric and complex Hermitian ensembles.  

Edge universality was first obtained by Soshnikov for a wide class of Wigner ensembles via the moment method \cite{sosh}.  This method required symmetry of the distribution of the individual matrix elements.  This requirement was partially removed in \cite{pechesosh}.  A different approach to edge universality based on direct comparison to Gaussian ensembles was developed in \cite{tao2010random,erdos2012rigidity}.  The work \cite{erdos2012rigidity} proves edge universality for general Wigner matrices under only a high finite moment condition.  
 This latter condition  was partially removed in \cite{erdos2012spectral}. In the work \cite{leeyin} it was shown that finiteness of the fourth moment is an essentially optimal condition for edge universality of Wigner matrices.  An almost-optimal necessary condition was obtained earlier in \cite{auf}.  

In contrast with the existing work on bulk universality of Wigner matrices and the three step strategy outlined above, these works make no use of Dyson Brownian motion.  As shown in \cite{erdos2012rigidity}, edge universality can be proven by directly comparing the eigenvalue statistics of an arbitrary Wigner matrix to that of the corresponding Gaussian ensemble.  To put it another way, the Green function comparison theorem can be used even if the Gaussian component is large, i.e., $t$ order $1$.  This is in contrast to the bulk where it is required that $t \ll 1$ in order to match Green's functions.

The work \cite{bourgade2014edge} implemented the above three-step strategy to prove edge universality of {\it generalized} Wigner matrices, which are ensembles in which the matrix element variances are allowed to vary.  Such ensembles cannot be directly matched to the GOE/GUE and so it was required to establish first edge universality for a Gaussian divisible ensemble.  The approach there uses a certain random walk representation of the correlation functions as well as the uniqueness of Gibbs measures of local log gasses.  Moreover, the work \cite{bourgade2014edge} establishes the edge universality of $\beta$-ensembles for general potentials and $\beta \geq 1$.

The three-step strategy offers an attractive route to analyze the edge universality of general random matrix models.  For general ensembles the moment method may fail, especially if the spectral edge is not an extremal eigenvalue - for example, models whose limiting spectral density has multiple intervals of support, or adjacency matrices of random graphs which may have outlying eigenvalues seperated from the density of states. 

In the present work we analyze the local ergodicity of the edge statistics of Dyson Brownian motion for general initial data $V$.  Our main result is that if the initial data has an edge with square root behavior down to a scale $\eta_* \geq N^{-2/3}$, then the distribution of the first eigenvalue at that spectral edge is given by the Tracy-Widom law in the limit $N \to \infty$ for $t \gg \sqrt{ \eta_*}$.  The assumption on ``square root behavior'' is quantified in terms of the imaginary part of the Stieltjes transform of $V$.

In particular, this result solves the DBM component of applications of the three-step strategy to edge universality of general random matrix ensembles.  In a joint work with J. Huang \cite{sparseedge}, we apply the main result of the present paper to analyze the edge statistics of sparse random matrix ensembles.

Another approach to edge universality has been developed by Lee-Schnelli in the works \cite{schnelliedge1,schnelliedge2,schnelliedge3}.  In this approach, the edge statistics of ensembles are calculated by direct comparison of the Green's function to the GOE via a continuous interpolation.  To rephrase this in our language, the change in expectation of Green's function elements under the DBM flow are carefully calculated for long times $t \gg 1$.  
 In the work \cite{erdos2012rigidity} it was noted that naive power counting arguments can be improved using higher order cancellations in the trace of the Green's function.  The works \cite{schnelliedge1,schnelliedge2,schnelliedge3}  can be interpreted as systematically extending this cancellation to arbitrarily high orders.

 
 The methodology and scope of \cite{schnelliedge1,schnelliedge2,schnelliedge3} are unrelated to the present work.  Our work concerns models of the form $A + \gamma B$ for GOE $B$ and arbitrary $A$, whereas in  \cite{schnelliedge1,schnelliedge2,schnelliedge3}, $A$ is a random matrix with independent entries (up to $A=A^*$), such as a Wigner or sparse random matrix.  We allow $\gamma = N^{-c}$ and find the optimal size which can be viewed as Dyson's conjecture of local ergodicity for edge statistics.




 
  
  
  A short-time edge universality result for general initial data is useful for random matrix theory.  Our result is applied in a joint work with J. Huang on the edge statistics of sparse Erd\H{o}s-R\'enyi matrices \cite{sparseedge}.  In \cite{sparseedge}, Green's function methods are used to uncover a Gaussian shift to the position of the extremal eigenvalues in the regime that the edge probability satisfies $ N^{-7/9} \ll p \ll N^{-2/3}$.  After subtracting this shift, the edge statistics are compared to a Gaussian divisible ensemble; the present work implies that the latter has Tracy-Widom fluctuations.  This implies that the eigenvalue gaps near the edge have the same distribution as the GOE, and that at $p=N^{-2/3}$, the extremal eigenvalues converge to a sum of Gaussian and Tracy-Widom distribution.


Our anaylsis of DBM is based around the coupling idea of \cite{homogenization} and matching idea of \cite{landonyau}.  At a time $t_0 \gg \sqrt{ \eta_*}$ we find that the free convolution, which gives the macroscopic eigenvalue density of DBM, has a square root behavior at the edge.  We re-scale and shift the DBM so that this edge matches the semicircle density of the GOE.  For times $ t \geq t_0$, the Dyson Brownian motion is then coupled (as in \cite{homogenization}) to DBM with initial data the matching GOE ensemble.  Under this coupling, the difference between these two stochastic eigenvalue locations obeys a simple discrete parabolic equation. 

In order to analyze the discrete parabolic equation, it is convenient to localize our analysis by introducing a ``short-range approximation'' to Dyson Brownian motion. Coupling the short-range approximations instead of the full DBMs leads to a parabolic equation whose heat kernel has rapid off-diagonal decay.  Such finite speed estimates first appeared in \cite{que}.  These rapid decay estimates offer several advantages.  Due to the fact that the free convolution law varies on the scale $t^2$ near the edge, the DBM may only be matched to a GOE matrix locally, and not globally.  The finite speed estimates then allow for a cut-off of these non-matching elements.  Moreover, as our assumptions on $V$ are only local, it is possible that the DBM evolution away from the edge is irregular - the finite speed estimates ensure that this does not affect the behavior at the edge.

Our use of the short-range approximation and the finite speed estimates of \cite{que} ensure that we do not rely on level repulsion estimates.  Level repulsion estimates were proven for DBM in the bulk in \cite{landonyau}.  Near the edge, these kinds of results for Dyson Brownian motion are unknown and do not appear to be a direct generalization of the method in the bulk.

Our analysis of the resulting parabolic equation is based around the energy estimate of \cite{bourgade2014edge}.  The input of this into our work is that the $\ell^\infty$ norm of the solution to the parabolic equation is much smaller than $N^{-2/3}$ for times $t \gg N^{-1/3}$.  Hence, we find that the edge statistics of DBM match those of the coupled GOE ensemble down to scale $N^{-2/3}$, yielding the universality.

An additional input of our work is rigidity for Dyson Brownian motion.  Rigidity for long times $t \asymp 1$ were established in \cite{kevin1} via matrix methods, and in the bulk this method was extended to short times in \cite{landonyau}.  In \cite{burger} rigidity in the bulk was established using purely dynamical methods.  In this work we need to establish rigidity at the edge; we chose to do so using matrix methods which is a straightforward extension of \cite{kevin1,kevin2}.  The purely dynamical approach does not appear to be easily extended to the edge. 

The remainder of the work is organized as follows.  In \ref{sec:main} we define precisely our model and state our main results.  Our main contribution is in Section \ref{sec:dbm}, which is our analysis of DBM via coupling as outlined above.  Section \ref{sec:finitespeed} contains an auxilliary calculation needed in Section \ref{sec:dbm}.  The local deformed law is proven in Sections \ref{sec:ll1} and \ref{sec:ll2}.  Finally we analyze the regularity of the deformed semicircle law in Section \ref{sec:fcanal}.

\subsection{Asymptotic notation}
The fundamental large parameter in our paper is $N$, the size of our matrices, and all asymptotic notation is wrt $N$. 

We use $C$ to denote generic $N$-independent constants, the value of which may change from line to line in proofs.

For two possibly $N$-dependent nonnegative parameters $X$, $Y$ we use the notation
\beq
X \asymp Y
\eeq
to denote the fact that there is a constant $C>0$ so that
\beq
\frac{1}{C} X \leq Y \leq C X.
\eeq
For two possibly $N$-dependent parameters $X$ and $Y$, with $X$ possibly complex and $Y$ nonnegative, the notation
\beq
X \leq \O (Y) , \qquad X = \O (Y)
\eeq
means that $|X| \leq C Y$ for some constant $C>0$.  The notation $X = Y + \O (Z)$ means $X-Y = \O (Z)$.

\section{Definition of model and main results} \label{sec:main}

In this paper we will consider models of the form
\beq
H_t := V + \sqrt{t} G
\eeq
where $V$ is a deterministic diagonal matrix, and $G$ is a GOE matrix.   We define the Stieltjes transform of $V$ by
\beq
m_V (z) = \frac{1}{N} \sum_i \frac{1}{V_i - z}.
\eeq
WLOG we assume that $V_i$ are in increasing order.
We consider the following class of $V$.
\bed \label{def:main}
Let $\etast$ be a parameter satisfying
\beq
\etast := N^{-\phist}
\eeq
for some $0 < \phist \leq 2/3$. 
We say that $V$ is $\etast$-regular if 
\begin{enumerate}
\item \label{item:ass1} There is a constant $C_V >0$ such that
\beq \label{eqn:mvass1}
\frac{1}{C_V}  \frac{ \eta}{ \sqrt{ |E| + \eta } } \leq \Im [ m_V (E + \i \eta ) ] \leq C_V \frac{ \eta}{ \sqrt{ |E| + \eta } }, \qquad -1 \leq E \leq 0, \quad \etast \leq \eta \leq 10,
\eeq
and
\beq \label{eqn:mvass2}
\frac{1}{C_V}  \sqrt{ |E| + \eta } \leq  \Im [ m_V (E + \i \eta ) ] \leq C_V \sqrt{ |E| + \eta }, \qquad 0 \leq E \leq 1, \quad \etast^{1/2} \sqrt{ |E| } + \etast \leq \eta \leq 10.
\eeq
\item \label{item:ass2} There are no $V_i$ in the interval $[-1, - \etast]$.
\item We have $||V|| \leq N^{C_V}$ for some $C_V >0$.
\end{enumerate}
\eed
There are many possible reformulations of the above assumptions.  We summarize a few of these observations in the remarks below.  

\remark \begin{enumerate}
\item The motivation for the upper and lower bounds of assumption \ref{item:ass1} is as follows.  If $m(z)$ is the Stieltjes transform of a measure with density $\rho (x)$ such that $\rho (x) \asymp \1_{ \{ x \geq 0 \} } \sqrt{x}$ then the estimates of assumption \ref{item:ass1} holds for $\Im [ m (z)]$ all $z = E + \i \eta$ near $0$.
\item The first two assumptions \ref{item:ass1}, \ref{item:ass2} are equivalent (up to constants) to the estimates of assumption \ref{item:ass1} holding, as well as the estimate \eqref{eqn:mvass1} holding  in the larger  regime $10 \geq \eta \geq 0$ and $-1 \leq E \leq \etast$. 
\item The second assumption together with only the estimate \eqref{eqn:mvass2} imply that \eqref{eqn:mvass1} holds on a slightly smaller domain $-1 + c \leq E \leq 0$, any $c>0$.  For clarity, we have chosen to list both estimates of assumption \ref{item:ass1} regardless of this redundancy.
\item The choice of the interval $[-1, 1]$ and the $10$ in $10 \geq \eta$ in the first two assumptions plays no role, it is only important that the estimates hold in an order $1$ interval near $0$.  We just use the above for notational simplicity.  
\item Our main result below, Theorem \ref{thm:main}, in fact holds under weaker assumptions regarding the width of the interval on which we assume estimates for $V$ - that is, the width of the interval may go to $0$ with $N$.  For simplicity of proofs we have not explored the optimal assumption.
\item The above set-up is for an extremal eigenvalue at a left edge.  Of course, one can consider also a right edge, etc.
\end{enumerate}
We now state our main result.  We denote the eigenvalues of $H_t$ by $\lambda_i$.
\bet \label{thm:main}
Let $V$ be $\etast$-regular.  Let $t$ satisfy $N^{-\eps} \geq t \geq N^{\eps} \etast$.  Recall that $V_i$ are indexed in increasing order.  Let $i_0$ be the index of the first $V_i \geq -1/2$.  Fix $k$ a nonnegative integer, and let $F : \rr^{k+1} \to \rr$ be a test function such that
\beq
|| F||_\infty \leq C, \qquad || \nabla F ||_\infty \leq C.
\eeq
 There are deterministic scaling factors $\gamma_0$ and $E_-$ depending only on $V$ such that,
\begin{align}
\bigg| \ee[ F (  \gamma_0 N^{2/3}& (\lambda_{i_0} - E_- ), \cdots  , N^{2/3} \gamma_0 ( \lambda_{i_0+k} - E_- ) ) ]\notag\\
 &- \ee^{(GOE)}[ F( N^{2/3} ( \mu_1 +2 ), \cdots ,  N^{2/3} ( \mu_{1+k } + 2 ) ) ] \bigg| \leq N^{-c}
\end{align}
for some $c>0$.  The latter expectation is with respect to the eigenvalues $\mu_i$ of a GOE.  The scaling factor $\gamma_0 \asymp 1$ and the magnitude of $E_-$ is bounded above.
\eet

In the next subsection we define the scaling factors $\gamma_0$ and $E_-$.  They are defined in terms of the {\it free convolution} of $V$ and the semicircle distribution.

\subsection{Free convolution} \label{sec:fcint}
The eigenvalue density of $H_t$ is described by the free convolution of $V$ with the semicircle distribution at time $t$, which we denote by $\rhofct$.  The free convolution is well studied and we refer the reader to, e.g., \cite{biane} for its properties. It is defined by its Stieltjes transform which is the unique solution to the following fixed point equation, that has  the property $|m(z) | \sim |1/z|$,
\beq
\mfct = m_V ( z + t \mfct (z) ) = \frac{1}{N} \sum_i \frac{1}{V_i - z -t \mfct (z) }.
\eeq
For $t>0$, the density $\rhofct$ is analytic and may be recovered by the Stieltjes inversion formula,
\beq
\rhofct (E) = \lim_{\eta \to 0 } \frac{\Im [ \mfct (E + \i \eta ) ]}{ \pi }.
\eeq

The following lemma will be proven in Section \ref{sec:fcanal}.  It follows from Lemma \ref{lem:rhofctexpand}.
\bel
Let $V$ be $\etast$-regular and $N^{-\eps} \geq t \geq N^{\eps} \etast$.  The support of the restriction of $\rhofct$ to $[-3/4, 3/4]$ consists of a single interval of the form $[E_-(t), 3/4]$ where 
\beq
|E_- (t) | \leq C t \log(N)
\eeq
for a constant $C$.  For $|E-E_- | \leq c t^2$ we have the expansion
\beq \label{eqn:rhofcsq}
\rhofct (E) = \1_{ \{ E \geq E_- \}} \gamma_0^{-1/2} \sqrt{E-E_-}\left(1 + t^{-2} \O ( |E-E_- | ) \right)
\eeq
where the scaling factor $\gamma_0$ is defined by
\beq
\gamma_0 := \left( \frac{t^3}{N} \sum_i \frac{1}{ (V_i - E_- - t \mfct (E_- ) )^3} \right)^{-1/3} \asymp 1.
\eeq
The scaling factors $\gamma_0$ and $E_-(t)$ are those given in the statement of Theorem \ref{thm:main}.
\eel

\subsubsection{Free convolution properties and associated notation} \label{sec:fcnot}
In this short subsection we summarize some of the key properties of $\mfct$.   We also introduce some notation which will be used when dealing with the free convolution.  Define
\beq
\kappa := |E - E_- |, \qquad \xi (z) := z + t \mfct (z),
\eeq
and
\beq
g_i := \frac{1}{ V_i - z - t \mfct }, \qquad R_k := \frac{1}{N} \sum_i g_i^k
\eeq
The following lemma is proved in Section \ref{sec:fcanal}.  It follows from Lemma \ref{lem:immfct1}.
\bel
Let $V$ and $t$ be as above.  We have
\beq
\Im [ \mfct ] \asymp \frac{ \eta}{ \sqrt{ \kappa + \eta } }, \qquad -3/4 \leq E \leq E_-,
\eeq
and
\beq
\Im [ \mfct ] \asymp \sqrt{ \kappa + \eta}, \qquad E_- \leq E \leq 3/4.
\eeq
\eel

We will also use the notion of {\it overwhelming probability}.
\bed
We say that an event $\F$ or possibly a family of events $\F_u$ with $u$ in some index set $\I$, hold with overwhelming probability if
\beq
\inf_{u \in \I} \pp [ \F_u ] \geq 1 - N^{-D}
\eeq
for any $D>0$, for large enough $N$.  We say that an event $\F_1$ holds with overwhelming probability on an event $\F_2$ if
\beq
\pp [ \F_1^c \cap \F_2 ] \leq N^{-D}
\eeq
for any $D>0$ for large enough $N$.  We also have a similar notion as above for families of events $\F_{i, u}, u \in \I$.
\eed

\section{DBM calculations} \label{sec:dbm}

In this section we fix two time scales
\beq
t_0 = \frac{N^{\om_0}}{N^{1/3}}, \qquad t_1 = \frac{N^{\om_1}}{N^{1/3}},
\eeq
with
\beq
0 < \om_1 < \om_0/100.
\eeq
The purpose of the introduction of these scales is to first ``regularize'' the global eigenvalue density that $H_t$ follows.  For times $t_0 \leq t \leq t_0 + t_1$ we use a coupling idea of \cite{homogenization}, and couple the DBM process started from $H_{t_0}$ to a GOE ensemble.  For times $t_0 \leq t \leq t_0 + t_1$ we will have to track the evolution of the edge of the ensemble, but the scaling factor (i.e., the size of the eigenvalue gaps) will remain approximately constant due to the fact that $t_1 \ll t_0$.

Recall the definition of $i_0$ as in the statement of Theorem \ref{thm:main}.  For $ t \geq 0$ we define the process $\lambda_i$ by
\beq
\d \lambda_i = \frac{ \d B_{i-i_0+1}}{ \sqrt{N}} + \frac{1}{N} \sum_j \frac{1}{ \lambda_i - \lambda_j } \d t,
\eeq
with initial data 
\beq
\lambda_i (0) = \lambda_i ( \gamma_0 H_{t_0} ),
\eeq
where $\gamma_0$ is as defined in Section \ref{sec:fcint}.
Note that for each time $t$, the process $\left\{ \lambda_i (t) \right\}_i$ is distributed like the eigenvalues of the matrix $\gamma_0 V + \sqrt{\gamma_0^2 t_0 + t} G$ where $G$ is a GOE matrix.  Above $\left\{ B_i \right\}_{-N \leq i \leq  N}$ are standard independent Brownian motions.  In terms of the free convolution law, the edge of the $\lambda_i (t)$ is given by
\beq
E_\lambda (t) := \gamma_0 E_- (t_0 + t^2/ \gamma_0^2 ).
\eeq
Note that $\gamma_0$ is fixed - it was chosen from time $t_0$, but the edge continues to evolve in time.

The purpose of the re-scaling by $\gamma_0$ is so that constant scaling the square root in \eqref{eqn:rhofcsq} is changed to $1$ at time $t_0$, matching the semicircle.


We now define $\mu_i$ as the unique strong solution to the SDE,
\beq
\d \mu_i = \frac{ \d B_{i  }} { \sqrt{N}} + \frac{1}{N} \sum_j \frac{1}{ \mu_i - \mu_j } \d t,
\eeq
with initial data $\mu_i (0)$ being distributed like the eigenvalues of a GOE matrix independent of $H_{t_0}$.  Then for each time $t$, the $\left\{ \mu_i (t) \right\}_i$ are distributed like the eigenvalues of $\sqrt{1+t} G$ for $G$ a GOE matrix, and their edge is given by
\beq
E_\mu (t) = -2 \sqrt{1+t}.
\eeq

The main result of this section is the following.
\bet \label{thm:maindbm}
Let $t_1$ be as above.  With overwhelming probability, we have
\beq
| (\lambda_{i_0 + i-1} (t_1) - E_\lambda (t_1 ) ) - ( \mu_i - E_\mu (t_1 ) ) | \leq \frac{1}{N^{2/3+c} }
\eeq
for a $c>0$ and for any finite $1 \leq i \leq K$.
\eet

Note that the Brownian motions for $\mu_1$ and $\lambda_{i_0}$ are identical.  At this point, we want to take the difference of the $\mu_i$ and $\lambda_j$, but we need to pad each system with dummy particles so that the difference is defined for all indices $i$.

More precisely, we let the system $\left\{ x_i \right\}_{-N \leq i \leq N}$ of $2N+1$ particles be defined by
\beq
\d x_i = \frac{ \d B_i }{ \sqrt{N}} + \frac{1}{N} \sum_j \frac{1}{ x_i - x_j } \d t,
\eeq
with initial data 
\beq
x_i (0) = \begin{cases} -3 N^{C_V} + i N, & -N  \leq i \leq 1- i_0 \\
\lambda_{i +i_0 - 1} (0), &  2 - i_0 \leq i \leq N +1 - i_0 \\
3 N^{C_V} + i N, & N + 2 - i_0 \leq i \leq N 
\end{cases},
\eeq
and $y_i$ be defined by
\beq
\d y_i = \frac{ \d B_i }{ \sqrt{N}} + \frac{1}{N} \sum_j \frac{1}{ y_i - y_j } \d t,
\eeq
with initial data
\beq
y_i (0) = \begin{cases} -3 N^{C_V} + i N, & -N \leq i \leq 0 \\
\mu_{i} (0), & 1 \leq i \leq  N \end{cases}.
\eeq
The following follows from Appendix C of \cite{fixed}, and we omit the proof.
\bel
With overwhelming probability, the following estimates hold.  We have,
\beq
\sup_{0 \leq t \leq 1 }  \sup_{ 2 - i_0 \leq i \leq  N + 1 - i_0 } |x_i (t) - \lambda_{i+i_0 - 1 } (t) | \leq \frac{1}{N^{100}},
\eeq
and
\beq
\sup_{0 \leq t \leq 1 } \sup_{ 1 \leq i \leq  N} |y_i (t) - \mu_{i} (t) | \leq \frac{1}{N^{100}},
\eeq
and
\begin{align}
\sup_{0 \leq t \leq 1 } x_{1 - i_0 }(t) &\leq - 2 N^{C_V}  \notag\\
\inf_{0 \leq t \leq 1} x_{N+2 - i_0 } (t) & \geq 2 N^{C_V} \notag\\
\sup_{0 \leq t \leq 1 } y_{0} (t) & \leq - 2 N^{C_V}.
\end{align} 
\eel

\subsection{Interpolation} 

While we would like to take the difference $x_i - y_i$ directly, the equation that this function satisfies is quite singular.  Instead, we introduce an interpolation that results in a better equation.  A similar interpolation appeared in \cite{fixed}.

We define the following interpolating processes for $0 \leq \alpha \leq 1$.
\beq
\d z_i (t,  \alpha ) = \frac{ \d B_i}{ \sqrt{N}} + \frac{1}{N} \sum_i \frac{1}{ z_i (t, \alpha )- z_j (t, \alpha ) } \d t 
\eeq
with initial data
\beq
z_i (0, \alpha ):= \alpha x_i (0) + (1 - \alpha ) y_i (0).
\eeq
We define
\beq
m_N (t, \alpha ) = \frac{1}{N} \sum_i \frac{1}{ z_i (t, \alpha ) - z }.
\eeq
Define $\rho_{x, t}$ to be the free covolution of $\gamma_0 V$ with the semicircle at time $t_0 +t/\gamma_0^2$.  The edge is given by $E_x (t) = E_\lambda (t)$ as above.  Denote $\rho_{y, t}$ to be the semicircle at time $\sqrt{1+t}$ with edge $E_y (t) = E_\mu (t)$ as above.

By our choice of $\gamma_0$ (see Lemma \ref{lem:rhofctexpand}) we have
\beq \label{eqn:rhox}
\rho_{x, 0} (E + E_x (0) ) = \rho_{y, 0} (E+E_y (t) ) \left(1 + \O \left( \frac{|E|}{t_0^2} \right) \right), 0 \leq  E\leq c t_0^2
\eeq
for some $c>0$.  
Let $\gamma_{x, i} (t)$ and $\gamma_{y, i} (t)$ be the quantiles - more precisely, we define them by
\beq
\frac{i}{N} = \int_{-1/2}^{\gamma_{x, i} (t) } \rho_{x, t} (E ) \d E = \int_{-10}^{\gamma_{y, i} (t) } \rho_{y, t} (E) \d E.
\eeq
 By the local law estimates of Section \ref{sec:ll1} (as well as a stochastic continuity argument as in Appendix B of \cite{fixed} to pass from fixed times $t$ to all times) we know that there is a $k_* \asymp N$ so that
\beq \label{eqn:xirig1}
 \sup_{0 \leq t \leq 10 t_1} | z_i (t, 1 ) - \gamma_{x, i} (t) | \leq \frac{N^{\eps}}{i^{1/3} N^{2/3}}, \quad 1 \leq i \leq k_*
\eeq
with overwhelming probability, for any $\eps >0$, and a similar estimate for $z_i (t, 0) = y_i (t)$.  An easy calculation using \eqref{eqn:rhox} gives that
\beq \label{eqn:gamxy}
|  ( \gamma_{x, i} (0) - E_x (0) ) - ( \gamma_{y, i}(0) - E_y(0) ) | \leq C \frac{ i^{4/3}}{N^{2 \om_0} N^{2/3} },
\eeq 
for $1 \leq i  \leq N^{6 \om_0/5}/C$.

\subsubsection{Construction of a density for interpolating ensembles}

We will need measures with well-behaved square root densities that give the eigenvalue density of the interpolating ensembles.  For this, we construct the following measures.  First, let $\rho_{x, 0}$ and $\rho_{y, 0}$ as above.

Define the eigenvalue counting functions near $0$ by
\beq
n_{x} (E) = \int_{-1/2}^E \d \rho_{x, 0} (E') , \qquad n_{y} (E) = \int_{-10}^E \d \rho_{y, 0} (E') ,
\eeq
and then the eigenvalue counting functions by
\beq
n_x ( \varphi_x (s) ) = s, \qquad n_y ( \varphi_y (s) ) = s.
\eeq
Define now
\beq
\varphi (s, \alpha ) := \alpha \varphi_x (s) + (1- \alpha ) \varphi_y (s)
\eeq
on the domain
\beq
\varphi (s, \alpha ) : [0, k_*/N ] \to [\alpha E_{x} (0) + (1-\alpha ) E_y (0), \alpha \varphi_x (k_*/N) + (1-\alpha ) \varphi_y (k_*/N) ] := F_\alpha.
\eeq
Define now the inverse function $n ( E, \alpha ) : F_\alpha \to [0, k_*/N]$ by
\beq
n ( \varphi (s, \alpha ), \alpha ) = s,
\eeq
and finally the density $\rho (E, \alpha )$ on the interval $F_\alpha$ by
\beq
\rho (E, \alpha ) := \frac{\d }{ \d E} n (E, \alpha).
\eeq
From the inverse function theorem we see that
\beq
\rho (E, \alpha ) = \left( \alpha ( \rho_x ( \varphi_x ( n (E, \alpha ) ) ) )^{-1} + (1- \alpha )( \rho_y ( \varphi_y ( n (E, \alpha ) ) ) )^{-1} \right)^{-1}
\eeq
from which we can deduce that
\beq
\rho (E +E_- ( \alpha ) , \alpha ) = \rho_y (E + E_{y} (0) ) \left( 1 + \O \left( |E|/t_0^2 \right) \right), \qquad 0 \leq E \leq c t_0^2 .
\eeq
Using these, we need to construct measures $\mu (E, \alpha )$ to which $m_N (0, \alpha )$ are close.  We construct these as follows.
\beq
 \d \mu (E, \alpha ) = \rho (E, \alpha ) \1_{ \left\{ E \in F_\alpha \right\} } \d E + \frac{1}{N} \sum_{ i \leq 0}   \delta_{ z_i (0, \alpha ) } (E)  + \frac{1}{N} \sum_{ i > k_* } \delta_{ z_i (0, \alpha ) } (E).
\eeq
The motivation for this definition is as follows.  We would like to take a deterministic density that matches $z_i (0, \alpha)$ with which to take a free convolution with.  However, we have no real control on the density away from $0$, so we can only find a deterministic density near $0$ - this is the role of $\rho (E, \alpha)$ - this density matches $z_i (0, \alpha)$ near $0$.  For the remaining particles which are an order $1$ distance away from our point of interest we just take $\delta$ functions.  While this part of the measure is random, we can control the effect that it has on deterministic quantities that we need.

We let now $\rho_t (E, \alpha )$ be the free convolution of $ \mu ( \alpha )$ with the semicircle at time $t$, and denote the Stieltes transform by $m_t (z, \alpha)$.  The properties of these measures are studied in Section \ref{sec:matching}.  In particular, they have a square root density which we denote by  $\rho_t (E, \alpha )$ with an edge which we denote by $E_{-} (t, \alpha)$.  

With $\mu$ as constructed above, it is not hard to see that the difference $m_N (z, 0, \alpha ) - m_0 (z, \alpha )$ obeys the estimates outlined at the start of Section \ref{sec:ll2}.  Let $\gamma_i (t, \alpha )$ be the classical eigenvalue locations wrt $\rho_t (E, \alpha )$.   To be more precise, they are defined by
\beq
\frac{i}{N} = \int_{E_- (t, \alpha ) }^{ \gamma_i (t, \alpha ) } \rho_t (E, \alpha ) \d E.
\eeq
 As a consequence of Section \ref{sec:ll2} we have the following.
\bel \label{lem:riginterpolating}
There is an $i_* \asymp N$ so that the following estimates hold.  We have,
\beq \label{eqn:rigzi}
\sup_{ 0 \leq \alpha \leq 1 } \sup_{ 0 \leq t \leq 10 t_1 } |z_i (t, \alpha ) - \gamma_i (t, \alpha ) | \leq \frac{ N^{\eps}}{N^{2/3} i^{1/3}}
\eeq
for $1 \leq i \leq i_*$, for any $\eps >0$ and with overwhelming probability.
\eel
While the measures $\rho_t (E, \alpha)$ are random, they obey estimates wrt deterministic quantities with overwhelming probability.  In particular, we see from Section \ref{sec:matching} (more precisely, Lemmas \ref{lem:match1} and \ref{lem:match2}) that
\bel \label{lem:matching}
For all $0 \leq E \leq c N^{ - 2 \eps  } t_0^2$ we have
\beq \label{eqn:densitymatch}
\rho_t (E + E_{-} (t, \alpha ) ) = \rho_{y, t} (E + E_{y, t} ) \left( 1 + \O \left( N^{\eps} t / t_0 \right) \right),
\eeq
and 
\beq \label{eqn:mtmatch}
\left| \Re[ m_t (E + E_{-} (t, \alpha ) ) - m_t (E_{-} (t, \alpha ) ) ] - \Re [ m_{y, t} (E + E_{y, t} ) - m_{t, y} (E_{y, t} ) ] \right| \leq C \frac{ |E| N^{\eps}}{ t_0},
\eeq
and for $c N^{-2 \eps} t_1 t_0 \leq E \leq 0$,
\beq
\left| \Re[ m_t (E + E_{-} (t, \alpha ) ) - m_t (E_{-} (t, \alpha ) ) ] - \Re [ m_{y, t} (E + E_{y, t} ) - m_{t, y} (E_{y, t} ) ] \right| \leq C\frac{|E|^{1/2} (t_1)^{1/2} N^{\eps}}{t_0^{1/2} }.
\eeq
as well as 
\beq \label{eqn:hatgamcompare}
\left| \gamhat_i (\alpha, t) - \gamhat_i (0, t) \right| \leq N^{\eps} \frac{i^{2/3} t}{ N^{2/3} t_0 }.
\eeq
\eel
Additionally, we have the following esimates, as proven in Appendix \ref{a:inter}.
\bel \label{lem:edges}
We have, 
\beq
| E_{-} (t, 1 ) - E_{x, t} | \leq N^{\eps} \left( t^3 + \frac{t}{N^{1/2}} \right)
\eeq
and
\beq
|E_{-} (t, 0) - E_{y, t} | \leq N^{\eps} \left( t^3 + \frac{t}{N^{1/2}} \right).
\eeq
\eel

\subsection{Short-range approximation}

It will be convenient to introduce the shifted $z_i (t, \alpha)$,
\beq
\tilz_i (t, \alpha ) := z_i (t, \alpha ) - E_{-} (t, \alpha),
\eeq
which obey the SDE
\beq
\d \tilz_i (t, \alpha ) = \frac{ \d B_i}{ \sqrt{N}} + \frac{1}{N} \sum_j \frac{1}{ \tilz_i (t, \alpha ) - \tilz_j (t, \alpha ) } \d t + \Re [ m_t (E_- (t, \alpha ), \alpha ) ] \d t.
\eeq
We also introduce the corresponding
\beq
\hatgam_i (t, \alpha ) := \gamma_i (t, \alpha ) - E_{-} (t, \alpha ).
\eeq
We now construct a ``short-range'' set of indices $\A \subseteq [[-N, N ]] \times [[-N, N]]$.  We choose $\A$ to be symmetric, i.e., $(i, j) \in \A$ iff $(j, i) \in \A$.  The definition of $\A$ requires the choice of
\beq
\ell := N^{\om_\ell}.
\eeq
We let
\begin{align}
\A := \left\{ (i, j) : i, j >0, |i-j| \leq \ell ( 10 \ell^2 + i^{2/3} + j^{2/3} ) \right\} \bigcup \left\{ (i, j) : i, j >i_*/2 \right\} \bigcup \left\{ (i, j ) : i, j \leq 0 \right\}.
\end{align}
The following is not essential, but it will simplify notation.  It is an exercise.
\bel
For each $i$, the set $\left\{ j : (i, j) \in \A \right\}$ is an interval of natural numbers.
\eel
It will be convenient to introduce the following short-range summation notation.  We define,
\beq
\sumAi_j := \sum_{ j : (i, j) \in \A}, \qquad \sumAci_j := \sum_{j : (i, j ) \notin \A }.
\eeq
We also need the integral analogs of the above.  For each $i$, we define the interval
\beq
\I_i ( \alpha, t) := [ \hatgam_{j_-} (\alpha, t) , \hatgam_{j_+} ( \alpha, t)]
\eeq
where for each $i$,  $\left\{ j : (i, j) \in \A \right\} := [[ j_-, j_+ ]]$.  We remark that we are only going to use the classical eigenvalue locations $\gamma_i ( \alpha, t)$ for indices $1 \leq i \leq c N$ for some small $c>0$ - for such locations, the measure $\rho_t ( E, \alpha )$ is well-behaved (i.e., here it has a square root density).  In particular, the behavior of the above intervals is relatively tame.

The short-range approximation to $\tilz$ is the process $\hatz$ defined as the solution to the following SDE.  It requires the choice of an additional parameter $N^{\om_A}$.  We introduce the notation
\beq
\J ( \alpha, t) := [-1/2, \hatgam_{3 i_*/4 } ( \alpha, t ) ],
\eeq
where $i_*$ is as in Lemma \ref{lem:riginterpolating}.  
Recall,
\beq
i_* \asymp N.
\eeq
For $i \leq 0$ we let
\beq
\d \hatz_i (t, \alpha ) = \frac{ \d B_i } { \sqrt{N}} + \frac{1}{N} \sumAi_j \frac{1}{ \hatz_i (t, \alpha) - \hatz_j (t, \alpha ) } + \frac{1}{N} \sumAci_j \frac{1}{ \tilz_i (t, \alpha ) - \tilz_j (t, \alpha ) } + \Re [ m_t (E_- (t, \alpha ), \alpha ) ],
\eeq
for $1 \leq i \leq N^{\om_A}$,
\begin{align}
\d \hatz_i (t, \alpha ) &= \frac{ \d B_i }{ \sqrt{N}} + \frac{1}{N} \sumAi_j \frac{1}{ \hatz_i (t, \alpha ) - \hatz_j (t, \alpha ) } + \int_{ \I^c_i (0, t) } \frac{1}{ \hatz_i (t, \alpha ) - E } \rho_t (E + E_- (t, 0), 0 ) \d E\notag\\
& + \Re[ m_t (E_- (t, 0), 0) ] \d t ,
\end{align}
for $N^{\om_A} \leq i \leq i_*/2$,
\begin{align}
\d \hatz_i (t, \alpha ) &= \frac{ \d B_i }{ \sqrt{N}} + \frac{1}{N} \sumAi_j \frac{1}{ \hatz_i (t, \alpha ) - \hatz_j (t, \alpha ) } +  \int_{ \I^c_i (\alpha, t) \cap \J ( \alpha, t) } \frac{1}{ \hatz_i (t, \alpha ) - E } \rho_t (E + E_- (t, \alpha ), \alpha ) \d E \notag\\
&+ \sum_{j \geq 3 i_*/4, j \leq 0} \frac{1}{ \tilz_i (t, \alpha ) - \tilz_j (t, \alpha ) }+ \Re[ m_t (E_- (t, \alpha), \alpha) ] \d t 
\end{align}
and for $i_*/2 \leq i \leq N$,
\beq
\d \hatz_i (t, \alpha ) = \frac{ \d B_i } { \sqrt{N}} + \frac{1}{N} \sumAi_j \frac{1}{ \hatz_i (t, \alpha) - \hatz_j (t, \alpha ) } + \frac{1}{N} \sumAci_j \frac{1}{ \tilz_i (t, \alpha ) - \tilz_j (t, \alpha ) } + \Re [ m_t (E_- (t, \alpha ), \alpha ) ],
\eeq
with initial data
\beq
\hatz_i (0, \alpha ) = \tilz_i (0, \alpha ).
\eeq
We pause to review the hierarchy of scales that we have introduced.  We have
\beq
t_0 = \frac{N^{\om_0}}{N^{1/3}}, \quad t_1 = \frac{N^{\om_1}}{N^{1/3}} , \quad \ell = N^{\ell}, \quad N^{\om_A}
\eeq
and
\beq
0 < \om_1 \ll \om_\ell \ll \om_A \ll \om_0.
\eeq
The purpose of the scale $\ell$ is to cut-off the effect of the initial data far from the edge - that is, for large $i$, we do not have $\tilz_i (t=0, \alpha=1 ) \sim \tilz_i (t=0, \alpha= 0)$, whereas for small $i$ they do in fact match.  Secondly, we will want to differentiate in $\alpha$.  The quantities that depend on $\alpha$, i.e., $\rho_t (E, \alpha )$ are approximately $\alpha$ independent near $E=0$ (after accounting for the drift $\del_t E_- (t, \alpha)$).  This regularity scale comes $\om_0$, and so we fix a scale $\om_A \ll \om_0$ on which we replace the $\alpha$-dependent quantities by $\alpha$-independent quantities.

We now show that $\hatz_i (t, \alpha )$ is a good approximation to $\tilz_i (t, \alpha )$.  Make the following choices of parameters: choose $\om_1 < \om_\ell / 10$, $\om_\ell < \om_A /20$ and $\om_A < \om_0 / 20$.   With these choices, the error term on the RHS of \eqref{eqn:shortrange} is $\ll N^{-2/3}$.
\bel \label{lem:shortrange}
With overwhelming probability,
\beq \label{eqn:shortrange}
\sup_{ 0 \leq \alpha \leq 1 } \sup_{ i } \sup_{ 0 \leq t \leq 10 t_1 } | \tilz_i (t, \alpha ) - \hatz_i (t, \alpha ) | \leq \frac{N^{\eps}}{N^{2/3}} \left(  \frac{ N^{\om_1}}{N^{2 \om_\ell} }+ \frac{N^{\om_1}}{N^{1/6}} + \frac{ N^{ 2 \om_A/3 +\om_1} }{N^{\om_0}} + \frac{ N^{\om_A/3 + 2 \om_1}}{ N^{ \om_0 } }\right)
\eeq
for any $\eps >0$.
\eel
\remark The important error on the RHS is the term $N^{\om_1} / N^{2 \om_\ell}$.  More precisely, it is important that $\om_\ell$ appear with a power strictly greater than $1$, for later estimates.

\proof Let $v_i = \tilz_i - \hatz_i$.  We have the equation
\beq
\del_t v = \B_1 v + \V_1 v + \zeta,
\eeq
where $\B_1$ is the operator defined by
\beq
( \B_1 v )_i = \frac{1}{N} \sumAi_j \frac{ v_i - v_j }{ (\tilz_i (\alpha) - \tilz_j ( \alpha )) ( \hatz_i (\alpha )- \hatz_j (\alpha))}
\eeq
the operator $\V_1$ is diagonal: $( \V_1 v )_i = \V_1 (i) v_i$, where for $\V_1 (i) = 0$ for $i \leq 0$ or $i \geq i_*/2$ and for $1 \leq i \leq N^{\om_A}$,
\beq
\V_1 (i) = - \int_{ \I^c_i (0, t) } \frac{ \rho_t (E+E_- (t, 0), 0) }{ ( \tilz_i ( \alpha, t) - E) ( \hatz_i ( \alpha, t) - E)}
\eeq
and for $N^{\om_A} \leq i \leq i_*/2$, 
\beq
\V_1 (i) = -\int_{ \I^c_i (\alpha, t) \cap \J ( \alpha, t) } \frac{ \rho_t (E+E_- (t, \alpha), \alpha) }{ ( \tilz_i ( \alpha, t) - E) ( \hatz_i ( \alpha, t) - E)}.
\eeq
In particular 
\beq
\V_1 \leq 0,
\eeq
and so the semigroup of $\B_1 + \V_1$ is a contraction on every $\ell^p$ space.  Hence, for the difference we have,
\beq
v(t) = \int_0^t \mathcal{U}^{\B_1 + \V_1 } (s, t) \zeta (s) \d s
\eeq
and so
\beq
||v (t) ||_\infty \leq t  \sup_{ 0 \leq s \leq t} || \zeta(s) ||_\infty,
\eeq
and so we must estimate $|| \zeta||_\infty$. 

The error term $\zeta$ is given by $\zeta_i = 0$ for $i \leq 0$ and $i \geq i_*/2$, and for $1 \leq i \leq N^{\om_A}$, it is given by
\begin{align}
\zeta_i &:= \int_{ \I^c_{i} ( 0, t) } \frac{ \rho_t (E + E_{-} (0, t) , 0) }{ \tilz_i (t, \alpha ) -E} - \frac{1}{N} \sumAci_j \frac{1}{ \tilz_i (t, \alpha ) - \tilz_j (t, \alpha ) } \notag\\
&+ \Re[ m_t ( E_- (t, \alpha ), \alpha ) ] - \Re [ m_t ( E_- (t, 0 ), 0) ] ,
\end{align}
and for $N^{\om_A} \leq i \leq i_*/2$ by
\beq
\zeta_i := \int_{ \I^c_{i} ( \alpha, t) \cap \J ( \alpha, t) } \frac{ \rho_t (E + E_{-} (\alpha, t) , \alpha) }{ \tilz_i (t, \alpha ) -E } - \frac{1}{N} \sumAci_{j \leq 3 i_*/4} \frac{1}{ \tilz_i (t, \alpha ) - \tilz_j (t, \alpha ) }.
\eeq
We need to estimate $|| \zeta||_\infty$.  
This term 
 is controlled by
\begin{align}
\left|  \int_{\J ( \alpha, t) \backslash \I_i ( \alpha, t) }  \frac{ \rho_t (E + E_{-} (\alpha, t) , \alpha) }{ \tilz_i (t, \alpha ) -E} - \frac{1}{N} \sumAci_{0 < j < 3 i_*/4} \frac{1}{ \tilz_i (t, \alpha ) - \tilz_j (t, \alpha ) } \right| &\leq \frac{ N^{\eps}}{N^{5/3} } \sumAci_{ 0 < j < 3 i_*/4} \frac{1}{ ( \hatgam_i - \hatgam_j )^2 j^{1/3} } \notag\\
&\leq \frac{ CN^{\eps}} {N^{1/3}}  \sumAci_{ 0 < j < 3 i_*/4} \frac{ i^{2/3} + j^{2/3} }{ (i-j)^2 j^{1/3}}.
\end{align}
We first estimate,
\begin{align}
\sumAci_{j>i} \frac{ i^{2/3} + j^{2/3} } { (i - j )^{2} j^{1/3}} &\leq   C\sumAci_{j>i} \frac{ j^{1/3} } { (i - j )^{2}} \notag \\
&\leq  C\sumAci_{j>i} \frac{ i^{1/3}}{ (i-j)^2} + \frac{1}{ (i-j)^{5/3}} \notag \\
&\leq C \frac{i^{1/3}}{ \ell ( \ell^2 + i^{2/3} ) } + \frac{C}{ ( \ell ( \ell^2 + i^{2/3} ) )^{2/3} } .
\end{align}
We also have
\begin{align}
\sumAci_{j<i} \frac{ i^{2/3} + j^{2/3} } { (i - j )^{2} j^{1/3}} &\leq C \sumAci_{j<i/2} \frac{ i^{2/3} } {  (i + \ell ( \ell^2 + i^{2/3} ) )^2 j^{1/3}} + C \sumAci_{i>j>i/2} \frac{ i^{1/3}} { (i - j )^{2} } \notag \\
&\leq C \frac{ i^{4/3}}{ (i + \ell ( \ell^2 + i^{2/3} ) )^2 } + C \frac{ i^{1/3}}{ \ell ( \ell^2 + i^{2/3} ) } \leq C \frac{ i^{1/3}}{ \ell ( \ell^2 + i^{2/3} ) }.
\end{align}
Hence,
\beq
\frac{ CN^{\eps}} {N^{1/3}}  \sumAci_{ 0 < j < 3 i_*/4} \frac{ i^{2/3} + j^{2/3} }{ (i-j)^2 j^{1/3}} \leq \frac{ C N^{\eps}}{N^{1/3} N^{2 \om_\ell} }.
\eeq
In conclusion, for $i \geq N^{\om_A}$,
\beq
|\zeta_i | \leq N^{\eps} \left( \frac{1}{N^{1/3} N^{2 \om_\ell }} \right)
\eeq
for any $\eps >0$ with overwhelming probability.  We now need to bound $\zeta_i$ for $1 \leq i \leq N^{\om_A}$.  We rewrite it as follows.
\begin{align}
\zeta_i &= \left( \int_{ \I^c_{i} ( \alpha, t) } \frac{ \rho_t (E + E_{-} (\alpha, t) , \alpha) }{ \tilz_i (t, \alpha ) -E } - \frac{1}{N} \sumAci_j \frac{1}{ \tilz_i (t, \alpha ) - \tilz_j (t, \alpha ) } \right) \notag\\
&+\Re[ m_t (E_- (t, \alpha ), \alpha ) - m_t ( \tilz_i (t, \alpha ) + E_- (t, \alpha ), \alpha ) - m_t (E_- (t, 0 ), 0 ) + m_t ( \tilz_i (t, \alpha ) + E_- (t, 0) , 0 ) ] \notag\\
&+\left( \int_{ \I_i (  \alpha, t ) } \frac{ \rho_t (E + E_- ( \alpha, t), \alpha ) }{ \tilz_i (t, \alpha ) - E } - \int_{ \I_i ( 0, t) } \frac{ \rho_t (E  + E_- (0, t ), 0  )}{ \tilz_i (t, \alpha ) - E } \right) =: B_1 + B_2 + B_3.
\end{align}
For the term $B_1$ we rewrite it as
\begin{align}
B_1 &= \int_{\J ( \alpha, t) \backslash \I_i ( \alpha, t) }  \frac{ \rho_t (E + E_{-} (\alpha, t) , \alpha) }{ \tilz_i (t, \alpha ) -E} - \frac{1}{N} \sumAci_{0 < j < 3 i_*/4} \frac{1}{ \tilz_i (t, \alpha ) - \tilz_j (t, \alpha ) } \label{eqn:zeta1} \\
&+ \int_{\J^c ( \alpha, t) } \frac{ \rho_t (E + E_{-} (\alpha, t) , \alpha) }{ \tilz_i (t, \alpha )-E } - \frac{1}{N} \sumAci_{j < 0, j \geq 3 i_*/4 }\frac{1}{ \tilz_i (t, \alpha ) - \tilz_j (t, \alpha ) }  \label{eqn:zeta2}
\end{align}
The term \eqref{eqn:zeta1} was handled above and is bounded by $N^{\eps} / (N^{1/3+2\om_\ell})$.  Now we estimate the term \eqref{eqn:zeta2} we write.  Fix for the moment an auxilliary scale $\eta_1$.  We write \eqref{eqn:zeta2} as
\begin{align}
&\int_{\J^c ( \alpha, t) } \frac{ \rho (E + E_{-} (\alpha, t) , \alpha) }{ \tilz_i (t, \alpha )-E } - \frac{1}{N} \sumAci_{j < 0, j \geq 3 i_*/4 }\frac{1}{ \tilz_i (t, \alpha ) - \tilz_j (t, \alpha ) } \notag\\
 = &\left( \int_{\J^c ( \alpha, t) } \frac{ \rho (E + E_{-} (\alpha, t) , \alpha) }{ \tilz_i (t, \alpha )-E } -\int_{\J^c ( \alpha, t) } \frac{ \rho (E + E_{-} (\alpha, t) , \alpha) }{ \tilz_i (t, \alpha )-E + \i \eta_1 } \right) \notag\\\\
+& \left(\sumAci_{j < 0, j \geq 3 i_*/4 }\frac{1}{ \tilz_i (t, \alpha ) - \tilz_j (t, \alpha ) + \i \eta_1} -\sumAci_{j < 0, j \geq 3 i_*/4 }\frac{1}{ \tilz_i (t, \alpha ) - \tilz_j (t, \alpha ) } \right)  \notag\\
+ &\left( \frac{1}{N} \sum_{ 1 \leq j \leq 3 i_*/4 } \frac{1}{ \tilz_i ( \alpha, t) - \tilz_j ( \alpha, t) + \i \eta_1 } - \int_{ J ( \alpha, t ) } \frac{ \rho (E + E_- ( \alpha, t), \alpha ) }{ \tilz_i (t, \alpha ) - E + \i \eta_1 } \right) \notag\\
+ & \left( m_N ( \tilz_i ( \alpha, t ) + \i \eta , t, \alpha ) - m_t ( \tilz_i ( \alpha, t) + \i \eta_1, \alpha ) \right)
 =: A_1 + A_2 + A_3 + A_4.\end{align}
By the local law we have 
\beq
|A_4| \leq \frac{ N^{\eps}}{N \eta_1}
\eeq
with overwhelming probability.  The term $A_2$ is bounded by
\beq
|A_2| \leq \eta_1 C \Im [ m_N (\tilz_i (t, \alpha ) + \i , \alpha ) ] \leq \eta_1 C
\eeq
with overwhelming probability.  Similarly, $|A_1| \leq \eta_1 C$.  A similar calculation to above yields
\beq
|A_3| \leq \frac{ N^{\eps}}{N \eta_1}.
\eeq
We optimize and choose $\eta_1 = N^{-1/2}$.  Hence, we see that
\beq
|B_1| \leq N^{\eps} \left( \frac{1}{ N^{1/3} N^{2 \om_\ell}} + \frac{1}{N^{1/2}} \right).
\eeq
From Lemma \ref{lem:matching} we see
\beq
|B_2| \leq \frac{ N^\eps}{N^{1/3}} \frac{ N^{ 2 \om_A/3}}{N^{\om_0}} + \frac{N^{\eps} N^{\om_1/2}}{N^{1/3} N^{\om_0/2}}.
\eeq
Here we used the fact that the largest index $j_+ (i)$ such that $(j, i) \in \A$ can be bounded by
\beq
j_+ (i) \leq C ( \ell^3 + i ) \leq C N^{\om_A}.
\eeq
We rewrite $B_3$ as
\begin{align}
B_3 &= \left( \int_{ \I_i ( \alpha , t ) } \frac{ \rho_t (E + E_- ( \alpha, t), \alpha ) }{ \tilz_i (t, \alpha ) - E } -\int_{ \I_i ( \alpha , t ) } \frac{ \rho_t (E + E_- ( 0, t), 0 ) }{ \tilz_i (t, \alpha) - E }  \right) \notag\\
&+\left(  \int_{ \I_i ( \alpha , t ) } \frac{ \rho_t (E + E_- ( 0, t), 0 ) }{ \tilz_i (t, \alpha) - E } - \int_{ \I_i (0, t ) } \frac{ \rho_t (E + E_- ( 0, t), 0 ) }{ \tilz_i (t, \alpha ) - E }  \right) =: D_1 + D_2.
\end{align}
Starting with $D_2$ we first see that
\beq
| \I_i ( \alpha, t) \Delta \I_i (0, t) | \leq \frac{ N^{\eps} N^{\om_1} ( N^{ 2 \om_\ell} +i^{2/3} )}{ N^{2/3} N^{\om_0} }
\eeq
where $\Delta$ is symmetric difference, and we used \eqref{eqn:hatgamcompare}.  On the above symmetric difference, the integrand is bounded by
\beq
\left|  \frac{ \rho_t (E + E_- ( 0, t), 0 ) }{ \tilz_i (t, \alpha) - E } \right| \leq \frac{ N^{1/3} (\ell + i^{1/3}) )}{\ell ( \ell^2 + i^{2/3} ) } ,
\eeq
and so
\beq
|D_2| \leq \frac{C N^{\eps} }{N^{1/3} } \frac{N^{\om_A/3} N^{\om_1}}{ N^{\om_\ell} N^{\om_0}},
\eeq
with overwhelming probability.  For $D_1$, the integral is a principal value so we have to do some minor case analysis to deal with the logarithmic singularity.  First assume $i \geq N^{\delta}$ for a $\delta < \om_{\ell}/10$.  Then in particular we know that $\tilz_i (t, \alpha)$ is at least distance $N^{-2}$ away from the boundary of $\I_i ( \alpha, t)$, and also that $| \tilz_i (t, \alpha ) > N^{-2}$.  Then, using \eqref{eqn:densitymatch} we have
\begin{align}
|D_1| &\leq \frac{ C N^{\eps} N^{\om_1} ( \ell^2 + i^{2/3} )}{N^{1/3} N^{\om_0} } \int_{ \I_i ( \alpha, t ) \cap |E-\tilz_i (\alpha, t ) | >N^{-50} } \frac{1}{ |\tilz_i ( \alpha, t ) - E | } \notag\\
&+ \left|  \int_{ |\tilz_i ( \alpha , t ) - E | < N^{-50} } \frac{ \rho_t (E + E_- ( \alpha, t), \alpha ) - \rho_t (E + E_- (0, t), 0 ) }{ \tilz_i (t, \alpha ) - E } \right|.
\end{align}
For the second term, we have on the domain of integration that $| \rho_t' (E + E_- ( \alpha, t) , \alpha ) | \leq C |\tilz_i (t, \alpha ) |^{-1/2} \leq CN$, and so
\beq
 \left|  \int_{ |\tilz_i ( \alpha , t ) - E | < N^{-50} } \frac{ \rho_t (E + E_- ( \alpha, t), \alpha ) + \rho_t (E + E_- (0, t), 0 ) }{ \tilz_i (t, \alpha ) - E } \right| \leq N^{-10}.
\eeq
Therefore, in the case that $i > N^{\delta}$ we have
\beq
|D_1| \leq \frac{ C N^{\eps} N^{\om_1} ( N^{ 2 \om_A/3} )}{N^{1/3} N^{\om_0} } 
\eeq
for any $\eps >0$ with overwhelming probability.  We now consider $i \leq N^{\delta}$.  First of all, if $| \tilz_i (t, \alpha ) | \geq N^{-100}$, then the above argument goes through and we obtain the same bound for $D_1$.  So we assume that $| \tilz_i (t, \alpha ) | \leq N^{-100}$.  We break the integral up into a few pieces (note that in this case $ 0 \in \I_i ( \alpha, t) $).  Denote the integrand by $G(E)$.  We write $D_1$ as
\begin{align}
D_1 &= \int_{\I_{i} ( \alpha, t ) \cap E > N^{-50} } G(E) +  \int_{3 \tilz_i /2 < E < N^{-50} }  G(E) +  \int_{\tilz_i /2 < E < 3 \tilz_i /2 }  G(E) + \int_{0 \leq E \leq \tilz_i/2}  G(E) \notag\\
&=: I_1  +I_2 +I_3 +I_4.
\end{align}
The term $I_1$ can be handled as in the case $i > N^{\delta}$.  For $I_2$, we just use that the integrand is bounded by $C |E|^{-1/2}$ and so $|I_2 | \leq C N^{-25}$.  If $\tilz_i ( \alpha, t | \leq 0$ then $I_3 = I_4 = 0$.  So we consider the case that $\tilz_i ( \alpha, t) > 0$. For $I_3$ that $| \rho_t' (E + E_{-} (t, \alpha), \alpha )| \leq C |\tilz_i ( \alpha, t) |^{-1/2}$ on the domain of integration of $I_3$ to obtain $|I_3 | \leq C |\tilz_i (t, \alpha )|^{1/2} \leq C N^{-25}$.  For $I_4$ we bound the integrand $|G(E) | \leq C |E|^{-1/2}$ and obtain $|I_4 | \leq CN^{-25}$.

This proves that
\beq
|D_1| \leq \frac{ C N^{\eps} N^{\om_1} ( N^{ 2 \om_A/3} )}{N^{1/3} N^{\om_0} } 
\eeq
for any $i$. We get the claim. \qed

The above implies the following, due to the fact that the choices of our parameters ensure that the RHS of \eqref{eqn:shortrange} is bounded by
\beq
\frac{1}{N^{2/3} }\frac{N^{\eps} N^{\om_1}}{N^{2 \om_\ell}},
\eeq
i.e., the first error term is the largest. 
\bel \label{lem:zhatrig}
Let $i \leq N^{3 \om_\ell + \delta}$ for $\delta < \om_\ell - \om_1$.  Then,
\beq
\sup_{ 0 \leq t \leq 10 t_1} \left| \hatz_i (\alpha, t)  - \hatgam_i (\alpha, t) \right| \leq \frac{N^{\eps}}{N^{2/3} i^{1/3} }
\eeq
with overwhelming probability.
\eel

\subsection{Differentiation}
Let $u_i (t, \alpha ) := \del_\alpha \hatz_i (t, \alpha )$.  We see that $u$ satisfies the equation,
\beq
\del_t u = \L u + \zetaz,
\eeq
where the operator $\L$ is as follows, and $\zetaz$ is a forcing term as follows.  First, 
\beq
\L = \B + \V
\eeq
where
\beq
( \B u )_i = \frac{1}{N} \sumAi_j \frac{ u_j - u_i}{ ( \tilz_i ( \alpha, t) - \tilz_j ( \alpha, t) )^2}, \qquad ( \V u)_i = \V_i u_i
\eeq
where for $1 \leq i \leq N^{\om_A}$,
\beq
\V_i = - \int_{ \I_i (0, t ) } \frac{\rho_t (E + E_- (0, t), 0 ) }{ ( \tilz_i ( \alpha, t) - E )^2 } 
\eeq
for $ N^{\om_A} \leq i \leq i_*/2$,
\beq
\V_i = - \int_{ \I_i (\alpha, t ) } \frac{\rho_t (E + E_( \alpha , t) , \alpha ) }{ ( \tilz_i ( \alpha, t) - E )^2 } 
\eeq
and $\V_i$ is $0$ otherwise.  The error term $\zetaz$ is $0$ unless $i \leq 0$ or $i \geq N^{\om_A}$, and it comes from the $\del_\alpha$ derivative hitting all the other terms.  It is not too hard to check that 
\beq
| \zetaz_i | \leq N^{C},
\eeq
for some $C>0$,
with overwhelming probability.  We also make the definition
\beq
\B_{ij} = - \frac{1}{N} \frac{1}{ ( \hatz_i (t, \alpha ) - \hatz_j (t, \alpha ) )^2}.
\eeq
Additionally, we introduce the $\ell^p$ norms,
\beq
||u||_p =\left(  \sum_i |u_i|^p \right)^{1/p} , \qquad ||u||_\infty = \max_i |u_i| = \lim_{p \to \infty} ||u||_p.
\eeq

\subsection{Long range cut-off}
First, we have the following finite speed of propogation estimate.  It follows from Lemma \ref{lem:fstech}.
\bel \label{lem:fs} For all small $\delta >0$, we have the following.  Let $a \leq N^{3 \om_\ell+\delta}$ and $b \geq N^{ 3 \om_\ell + 2 \delta}$.  Then
\beq
\sup_{ 0 \leq s \leq t \leq 10 t_1 } \UL_{ab} (s, t) + \UL_{ba} (s, t) \leq N^{-D}
\eeq
for any $D>0$ with overwhelming probability.
\eel

Fix now a small $\delta_v >0$. 
Define $v_i$ to be the solution to
\beq
\del_t v = \L v, \qquad v_i (0) = u_i (0) \1_{ \left\{ 1 \leq i \leq N^{3 \om_\ell + \delta_v } \right\} }.
\eeq
By Lemma \ref{lem:fs} and the fact that $\UL_{ij} = 0$ for $i \geq 1$ and $j \leq 0$ or $i \leq 0$ and $j \geq 1$, we immediately see the following.
\bel \label{lem:uvcutoff}
We have
\beq
\sup_{ 0 \leq t \leq 10 t_1} \sup_{ 1 \leq i \leq \ell^3} |u_i (t) - v_i (t) | \leq N^{-100}.
\eeq
\eel

\subsection{Energy estimate}

We require the following energy estimate.  
\bel \label{lem:energy}
Let $\delta_1 >0$ be small.
Let $w \in \rr^N$, $w_i = 0$ for $i \geq \ell^3 N^\delta$ or $i \leq 0$.  Then for all $\eps >0$ and all $\eta>0$ there is a constant $C$ (independent of $\eps$ and $\eta)$ s.t. for $0 \leq t \leq 2 t_1$,
\beq
|| \UL(0, t) w ||_\infty \leq C (p, \eta ) \left( \frac{ N^{C \eta + \eps}}{ N^{1/3} t } \right)^{3(1-6\eta)/p} || w||_p.
\eeq

\eel

For its proof we need the following Sobolev-type inequality from \cite{bourgade2014edge}.
\bel
For all $\eta >0$ there exists a $c_\eta >0$ s.t. 
\beq
\sum_{i \neq j \in \zz_+} \frac{ (u_i - u_j )^2}{ | i^{2/3} - j^{2/3} |^{2 - \eta } } \geq c_\eta \left( \sum_{i \in \zz_+} |u_i |^p \right)^{2/p}.
\eeq
\eel

The above lemma is used in the following which is proved via the Nash method, from which Lemma \ref{lem:energy} follows.  It is very similar to that appearing in \cite{bourgade2014edge}. 
\bel \label{lem:l2l1} Let $0 \leq s \leq t \leq t_1$.  Let $\delta_1 >0$ satisfy
\beq
0 < \delta_1 < \om_\ell - \om_1.
\eeq
Let $w$ be a vector s.t. $w_i = 0$ for $i \geq \ell^3 N^{\delta_1}$ and $i \leq 0$.
Let $\eta>0$ and $\eps>0$.  There is a $C>0$ independent of $\eps$ and $\eta$ and a constant $c_\eta$ s.t. the following hold with overwhelming probability for all $0 \leq s \leq t \leq 5 t_1$.
\beq \label{eqn:l2l1}
|| \UL (s, t) w ||_2 \leq \left( \frac{1}{ ( c_\eta N^{- C \eta - \eps} N^{1/3} (t-s) )^{3/2} } \right)^{1-6 \eta } || w ||_1
\eeq
and
\beq \label{eqn:l2l1T}
||\left( \UL (s, t)\right)^T w ||_2 \leq \left( \frac{1}{ ( c_\eta N^{- C \eta - \eps} N^{1/3} (t-s) )^{3/2} } \right)^{1-6 \eta } || w ||_1
\eeq

\eel

\proof We start with \eqref{eqn:l2l1}.  This is a modification of the argument that uses the usual Nash approach.  We provide all the details for completeness.  Fix $\delta_2$ and $\delta_3$ s.t.
\beq
0 < \delta_1 < \delta_2 < \delta_3 < \om_\ell - \om_1.
\eeq
For notational simplicity we just do $s=0$.  Fix $\eta >0$. Let
\beq
p = \frac{3}{ 1 + \eta }.
 \eeq
  Assume that
  \beq
  || w(0) ||_1 =1.
  \eeq
  We can also assume that
\beq
|| w( u) ||_p \geq \frac{1}{N^{100}}, \qquad 0 \leq u \leq t
\eeq
or else 
\beq
|| w (t) ||_\infty \leq \frac{1}{N^{95}}
\eeq
by the contraction property.  We have
\begin{align}
|| w(u) ||_p^2 & \leq \sum_{i, j \in \zz_+ } \frac{ ( w_i - w_j )^2}{ | i^{2/3} - j^{2/3} |^{2 - \eta } } \notag \\
& \leq \sum_{ (i, j) \in \A }  \frac{ ( w_i - w_j )^2}{ | i^{2/3} - j^{2/3} |^{2 - \eta } } + C \sum_i \sumAci_j \frac{ w_i^2 }{ | i^{2/3} - j^{2/3} |^{2 - \eta }  } .
\end{align}
We have
\begin{align}
\sum_{ (i, j) \in \A }  \frac{ ( w_i - w_j )^2}{ | i^{2/3} - j^{2/3} |^{2 - \eta } }  \leq \sum_{ \substack{(i, j) \in \A \\ i \mbox{ or } j \leq \ell^3 N^{\delta_2 } }}  \frac{ ( w_i - w_j )^2}{ | i^{2/3} - j^{2/3} |^{2 - \eta } }  + \sum_{ \substack{ (i, j) \in \A \\ i, j \geq \ell^3 N^{\delta_2 }} }  \frac{ ( w_i - w_j )^2}{ | i^{2/3} - j^{2/3} |^{2 - \eta } } .
\end{align}
By Lemma \ref{lem:fs},
\beq
\sum_{ \substack{ (i, j) \in \A \\ i, j \geq \ell^3 N^{\delta_2 }} }  \frac{ ( w_i - w_j )^2}{ | i^{2/3} - j^{2/3} |^{2 - \eta } } \leq \frac{1}{N^{400}}. \label{eqn:eng2}
\eeq
If $(i, j) \in \A$ and $i$ or $j$ is less than $\ell^3 N^{\delta_2}$ then both $i$ and $j$ are less than $  \ell^3 N^{\delta_3}$.  By Lemma \ref{lem:zhatrig} we have
\beq
| \hatz_i - \hatz_j | \leq \frac{N^{\eps} | i^{2/3} - j^{2/3} |  }{N^{2/3} }
\eeq
for all $\eps>0$ for such $i$ and $j$.  Therefore,
\beq
\sum_{ \substack{(i, j) \in \A \\ i \mbox{ or } j \leq \ell^3 N^{\delta_2 } }}  \frac{ ( w_i - w_j )^2}{ | i^{2/3} - j^{2/3} |^{2 - \eta } }  \leq - N^{-1/3} N^{\eps+C\eta} \sum_{i, j} \B_{ij} (w_i - w_j )^2 = - N^{-1/3} N^{\eps+C \eta} \langle w , \B w \rangle. \label{eqn:eng3}
\eeq
(recall $\B_{ij}$ are negative).
Similarly,
\begin{align}
 \sum_i \sumAci_j \frac{ w_i^2 }{ | i^{2/3} - j^{2/3} |^{2 - \eta }  } &\leq  \sum_{i \leq \ell^2 N^{\delta_3} } \sumAci_j \frac{ w_i^2 }{ | i^{2/3} - j^{2/3} |^{2 - \eta }  }  +\frac{1}{N^{200}} \notag \\
 &\leq - N^{C \eta} N^{-1/3} \sum_{i \leq \ell^3 N^{\delta_3}} w_i^2 \V_i + \frac{1}{ N^{400}} \notag \\
 &\leq - N^{C \eta} N^{-1/3} \langle w , \V w \rangle + \frac{1}{N^{400}}. \label{eqn:eng4}
 \end{align}
From \eqref{eqn:eng2}, \eqref{eqn:eng3} and \eqref{eqn:eng4} we obtain
\beq
|| w ( u) ||_p^2 \leq - N^{\eps + C \eta } N^{-1/3}<w , \L w> + \frac{1}{10}||w (u)||_p^2 \
\eeq
Therefore
\begin{align}
\del_u || w(u) ||_2^2 &=  \langle w, \L w \rangle \notag \\
&\leq - c_\eta N^{-1/3} N^{-\eps-C \eta} || w(u) ||_p^2 \notag \\
& \leq -  c_\eta N^{-1/3} N^{-\eps-C \eta} ||w(s)||_2^{ \frac{8-4 \eta}{3}} ||w(t)||_1^{\frac{-2-4\eta}{3}} \notag \\
& \leq   - c_\eta N^{-1/3} N^{-\eps-C \eta} ||w(s)||_2^{ \frac{8-4 \eta}{3}} ||w(0)||_1^{\frac{-2-4\eta}{3}}
\end{align}
where in the second inequality we used Holder and in the last we used the $\ell^1$ contractivity of $\UL$.  

From this we obtain
\beq
|| w(t) ||_2 \leq \left( \frac{1}{ ( c_\eta N^{- C \eta - \eps} N^{1/3} t )^{3/2} } \right)^{1-6 \eta } || w (0) ||_1
\eeq
as desired.

The proof of \eqref{eqn:l2l1T} is identical and follows by duality.  One considers the function
\beq
w(u) = \UL(u, t) w
  \eeq
  which satisfies
  \beq
  \del_u w(u) = \L(u) w(u).
  \eeq
  Note that the only inputs in the proof of \eqref{eqn:l2l1} are time-independent lower bounds on $\L$ and Lemma \ref{lem:fs} which holds for both $\UL$ and $(\UL)^T$. \qed

\noindent{\bf Proof of Lemma \ref{lem:energy}.}  Let
\beq
0 < \delta_1 < \delta_2 < \om_\ell - \om_1.
\eeq
Let $\chi_2 (i) = \1_{ \left\{ 1 \leq i \leq \ell^3 N^{\delta_2 } \right\} }$.  Let $v$ have $||v||_1 =1$.  We have,
\beq
\langle \UL (0, t) w, v \rangle = \langle w, (\UL)^T v \rangle = \langle w, ( \UL)^T \chi_2 v \rangle + \langle w , (\UL)^T (1 - \chi_2 ) v \rangle.
\eeq
We have by Lemma \ref{lem:fs},
\beq
| \langle w , (\UL)^T (1 - \chi_2 ) v \rangle | \leq \frac{1}{N^{100}} ||w||_2 ||v||_1.
\eeq
By Cauchy-Schwartz and Lemma \ref{lem:l2l1}
\begin{align}
\left| \langle w, ( \UL)^T \chi_2 v \rangle \right| &\leq ||w||_2 || (\UL)^T \chi_2 v ||_2 \notag \\
&\leq ||w||_2 \left( \frac{1}{ ( c_\eta N^{- C \eta - \eps } N^{1/3} t)^{3/2} } \right)^{1 - 6 \eta } ||v||_1.
\end{align}
Hence,
\beq
||\UL (0, t) w ||_\infty \leq \left( \frac{1}{ (c_\eta N^{- C \eta - \eps} N^{1/3} t )^{3/2} } \right)^{1-6\eta} ||w||_2
\eeq
and so by the semigroup property,
\beq
||\UL (0, t) w ||_\infty \leq \left( \frac{1}{ (c_\eta N^{- C \eta - \eps} N^{1/3} t )^{3} } \right)^{1-6\eta} ||w||_1.
\eeq
The rest follows from interpolation. \qed

\subsection{Proof of Theorem \ref{thm:maindbm}}
For notational simplicity we just do $i=1$.  
By Lemmas \ref{lem:edges} and \ref{lem:shortrange} we have,
\beq
| \lambda_{i_0} (t_1)- E_\lambda (t_1)  - ( \mu_1 (t_1) - E_\mu (t_1 ) ) | \leq | \hatz_1 (1, t_1) - \hatz_1 (0, t_1) | + \frac{1}{N^{2/3+c} },
\eeq
with overwhelming probability.  By the definition of $u_i (t)$ we have
\beq
 \hatz_1 (1, t_1) - \hatz_1 (0, t_1) = \int_0^\alpha u_1 (t_1, \alpha ) \d \alpha.
\eeq
By Markov's inequality and Lemma \ref{lem:uvcutoff} we have with overwhelming probability,
\beq
\left| \int_0^1 u_1 (t_1, \alpha) \d \alpha \right| \leq  N^{-10} + \int_0^1 v_1 (t_1, \alpha ) \d \alpha.
\eeq
Note that by  \eqref{eqn:gamxy} we see that
\beq
||v (0) ||_4 \leq \frac{ N^{\eps_1}}{N^{2/3}}
\eeq
for any $\eps_1 >0$, with overwhelming probability.  Hence by Lemma \ref{lem:energy} and Markov again we get that
\beq
\left| \int_0^1 v_1 (t_1, \alpha ) \d \alpha \right| \leq \frac{ N^{2 \eps_1}}{ N^{c \om_1} } \frac{1}{N^{2/3}}
\eeq
with overwhelming probability.  This yields the claim. \qed

\section{Finite speed calculations} \label{sec:finitespeed}

The following method for getting bounds on $\UL$ originates from \cite{que}.  
\bel
Let $\L$ and $\UL$ be as above.  Fix a small $\delta >0$ satisfying $\delta < \om_\ell - \om_1$.  Let $\om_1 < \om_\ell$.  For any $a \geq N^{3 \om_\ell+\delta}$ and $b \leq N^{ 3 \om_\ell + \delta}/2$ and fixed $s$ we have
\beq
\sup_{s \leq t \leq 10 t_1 } \UL_{ab} (s, t) + \UL_{ba} (s, t) \leq N^{-D}
\eeq
for any $D>0$.
\eel
\proof For notational simplicity we take $s=0$. 
Let $\psi (x)$ be a function as follows.  We let 
\beq
\psi = -x, \qquad |x| \leq \frac{ N^{2 \om_\ell} N^{\delta2/3}}{N^{2/3}},
\eeq
and
\beq
\psi' (x) = 0, \qquad |x| > 2 \frac{ N^{2 \om_\ell} N^{\delta2/3}}{N^{2/3}},
\eeq
and demand that $| \psi (x) - \psi (y) | \leq |x-y|$ and $| \psi' (x) | \leq 1$, that $\psi$ be decreasing, and 
\beq
| \psi'' (x) | \leq C \frac{ N^{2/3}}{ N^{2 \om_\ell + 2 \delta/3}}.
\eeq
We consider a solution of 
\beq
\del_t f = \L f
\eeq
with initial data $f_i (0) = \delta_{p_*}$ for any $p_* \geq N^{\delta} N^{3 \om_\ell} =: p$.  We consider the function
\beq
F(t) := \sum_k f_k^2 \e^{ \nu \psi ( \hatz_k (t, \alpha ) - \gamhat_p (t, \alpha ) ) } =: \sum_k f_k^2 \phi_k^2 =: \sum_k v_k^2
\eeq
This is obeying the equation
\begin{align}
\d F &= - \sum_{(i, j) \in \A} \B_{ij} ( v_i - v_j )^2 \d t - \sum_i \V_i v_i^2 \label{eqn:fs1} \d t \\
&+ \sum_{ (i, j) \in \A } \B_{ij} v_i v_j \left( \frac{ \phi_i}{\phi_j} + \frac{ \phi_j}{\phi_i } - 2 \right)  \label{eqn:fs2} \\
&+ \nu \sum_i v_k^2 \psi' ( \hatz_i - \hatgam_p ) \d ( \hatz_i  - \gamhat_i )  \label{eqn:fs3}\\
&+ \sum_i v_i^2 \left( \frac{ \nu^2}{N} ( \psi' ( \hatz_i - \hatgam_p ) )^2 + \frac{ \nu}{N} \psi'' ( \tilz_i - \hatgam_p ) \right) \label{eqn:fs4} \d t.
\end{align}
We define a stopping time $\tau$ as follows.  We can take $\tau_1$ to be a stopping time so that for $t < \tau_1$ the rigidity estimates of Lemma \ref{lem:shortrange} and \eqref{eqn:rigzi} hold, for a small $\eps >0$ with $\eps < \delta/1000$.  With overwhelming probability, $\tau \geq 10 t_1$.   Take $\tau_2$ to be the first time that $F \geq 5$, and then $\tau = \tau_1 \wedge \tau_2 \wedge 10 t_1$.  We want to prove $\tau = 10 t_1$ with overwhelming probability.
  In the remainder of the proof we work with times $t < \tau$.  
  
  We note that $\psi' ( \hatz_i - \gamhat_i ) = 0$ unless $i \leq C N^{ 3 \om_\ell + \delta}$ for some $C>0$.  Moreover, if $i \leq C N^{ 3 \om_\ell + \delta}$ and $(i, j) \in \A$, then $j \leq C' N^{3 \om_\ell + \delta}$.  From this we see that the nonzero terms in the sum in \eqref{eqn:fs2} have both $i, j \leq C N^{3 \om_\ell + \delta}$.  For such terms we have
\beq
|\hatz_i - \hatz_j | \leq \frac{ \ell^2 N^{\delta/3}}{N^{2/3}}.
\eeq
Hence, if 
\beq \label{eqn:nucond}
\frac{ \nu \ell^2 N^{\delta/3}}{ N^{2/3}} \leq C,
\eeq
then
\begin{align}
\sum_{ (i, j) \in \A } \B_{ij} v_i v_j \left( \frac{ \phi_i}{\phi_j} + \frac{ \phi_j}{\phi_i } - 2 \right) &\leq C \frac{ \nu^2}{N} \sum_i v_i^2 \sum_{j : (j, i ) \in \A } \1_{ \left\{ \phi_j \neq \phi_i \right\} } \notag\\
&\leq \frac{ \nu^2 \ell^3 N^{2 \delta/3}}{N} F(t).
\end{align}
The term \eqref{eqn:fs4} is easily bounded by
\beq
\left| \sum_i v_i^2 \left( \frac{ \nu^2}{N} ( \psi' ( \hatz_i - \hatgam_p ) )^2 + \frac{ \nu}{N} \psi'' ( \hatz_i - \hatgam_p ) \right) \right| \leq C \left( \frac{ \nu^2}{N} + \frac{ \nu}{N^{1/3}  N^{2 \om_\ell+2\delta/3}} \right) F(t).
\eeq
Now we deal with \eqref{eqn:fs3}.  Recall that $\psi' ( \hatz_i - \hatgam_p ) \neq 0$ only if $i \leq C N^{3 \om_\ell + \delta } \ll N^{\om_A}$.  Hence, for such $i$ we have
\begin{align} \label{eqn:dzgam}
\d ( \hatz_i ( \alpha, t) - \hatgam_p (\alpha, t) ) &= \frac{ \d B_i} { \sqrt{N}} + \frac{1}{N} \sumAi_j \frac{1}{ \hatz_i (t, \alpha ) - \hatz_j ( t, \alpha ) } + \int_{\I^c ( 0, t) } \frac{1}{ \hatz_i (t, \alpha ) - E } \rho_t (E + E_- (0, t), 0 ) \d E \notag\\
& + \Re[m_t (E_- (t, 0 ), 0 ]  + \Re [ m_t ( \hatgam_p ( \alpha, t) + E_- ( \alpha, t ), \alpha ) ] - \Re[ m_t ( E_- ( \alpha, t), \alpha ) ].
\end{align}
By the definition of $\tau$ and the BDG inequality, we have with overwhelming probability,
\beq
\sup_{ 0 \leq t \leq \tau} \int_0^t \sum_i v_i^2 \psi' ( \hatz_i - \gamhat_p ) \nu \frac{ \d B_i } { \sqrt{N}} \leq C N^{\eps} \left( \frac{ \nu^2  N^{\om_1}}{N^{4/3} } \right)^{1/2}.
\eeq
Next,
\beq
\frac{\nu}{N} \sum_{(i, j) \in \A} \frac{ \psi' ( \hatz_i - \gamhat_p ) v_i^2}{ \hatz_i - \hatz_j } = \frac{\nu}{2N} \sum_{(i, j) \in \A} \frac{ \psi' ( \hatz_i - \gamhat_p ) ( v_i^2-v_j^2)}{ \hatz_i - \hatz_j } + \frac{\nu}{2N} \sum_{ (i, j) \in \A} v_i^2 \frac{ \psi' ( \hatz_i - \gamhat_p ) - \psi' ( \hatz_j - \gamhat_p ) }{ \hatz_i - \hatz_j }.
\eeq
The latter sum is bounded by
\begin{align}
 \frac{\nu}{2N} \sum_{ (i, j) \in \A} v_i^2 \frac{ \psi' ( \hatz_i - \gamhat_p ) - \psi' ( \hatz_j - \gamhat_p ) }{ \hatz_i - \hatz_j } &\leq C \frac{ \nu}{N^{1/3} ( N^{2 \om_\ell +2 \delta/3} ) } \sum_i v_i^2 \sumAi_j \1_{ \left\{ \psi' (\hatz_i -\gamhat_p ) \neq \psi' ( \hatz_j - \gamhat_p ) \right\} } \notag\\
 & \leq C \frac{ \nu N^{\om_\ell}}{N^{1/3}} F(t) .
\end{align}
For the first term we use the Schwarz inequality, obtaining
\begin{align}
\frac{\nu}{2N} \sum_{(i, j) \in \A} \frac{ \psi' ( \hatz_i - \gamhat_p ) ( v_i^2-v_j^2)}{ \hatz_i - \hatz_j } &= \frac{\nu}{2N} \sum_{(i, j) \in \A} \frac{ \psi' ( \hatz_i - \gamhat_p ) ( v_i - v_j)(v_i + v_j ))}{ \hatz_i - \hatz_j }  \notag\\
&\leq \frac{1}{100N} \sum_{(i, j) \in \A} \frac{ ( v_i- v_j)^2}{( \hatz_i - \hatz_j )^2} \notag\\
&+\frac{ C \nu^2}{2N} \sum_{(i, j) \in \A}  \psi' ( \hatz_i - \gamhat_p )^2 (v_i^2 + v_j^2)
\end{align}
The second sum is bounded by
\beq
\frac{ C \nu^2}{2N} \sum_{(i, j) \in \A}  \psi' ( \hatz_i - \gamhat_p )^2 (v_i^2 + v_j^2) \leq \frac {C \nu^2 N^{3 \om_\ell} N^{ 2 \delta/3} }{N} F(t).
\eeq
In summary, we have estimated
\beq
\frac{\nu}{N} \sum_{(i, j) \in \A} \frac{ \psi' ( \hatz_i - \gamhat_p ) v_i^2}{ \hatz_i - \hatz_j } \leq \frac{1}{ 100  } \sum_{ (i, j) \in \A } \B_{ij} (v_i -v_j)^2 + F(t) C \left( \frac{ \nu N^{\om_\ell}}{N^{1/3}} + \frac{ \nu^2 N^{ 3 \om_\ell} N^{2 \delta/3}}{ N} \right).
\eeq
Returning to \eqref{eqn:dzgam} we write
\begin{align}
 &\int_{\I^c ( 0, t) } \frac{1}{ \hatz_i (t, \alpha ) - E } \rho_t (E + E_- (0, t), 0 ) \d E \notag\\
 + & \Re[m_t (E_- (t, 0 ), 0 ]  + \Re [ m_t ( \hatgam_p ( \alpha, t) + E_- ( \alpha, t ), \alpha ) ] - \Re[ m_t ( E_- ( \alpha, t), \alpha ) ] \notag\\
 = & \left( \int_{\I^c ( 0, t) } \frac{1}{ \hatz_i (t, \alpha ) - E } \rho_t (E + E_- (0, t), 0 ) \d E  + \Re[ m_t ( \hatgam_p  ( \alpha, t) + E_- ( 0, t), 0 ) ]  \right)\notag\\
 +&\bigg[ \Re[m_t (E_- (t, 0 ), 0 ]  + \Re [ m_t ( \hatgam_p ( \alpha, t) + E_- ( \alpha, t ), \alpha ) ] \notag\\
 -& \Re[ m_t ( E_- ( \alpha, t), \alpha ) ]  - \Re[m_t ( \hatgam_p ( \alpha, t) + E_- (0, t), 0)]  \bigg] =:A_1 + A_2.
\end{align}
By \eqref{eqn:mtmatch} we have
\beq
|A_2| \leq N^{\eps} \frac{ N^{2 \om_\ell + 2 \delta/3} }{N^{1/3} N^{\om_0} }.
\eeq
Using $\rho_t ( \hatgam_{2p} ( \alpha, t)  + E_- (0, t ), 0 ) \leq C N^{\om_\ell} N^{\delta/3} / N^{1/3}$ we see that
\beq
|A_1| \leq \left( \int_{ E \geq \hatgam_{2p} ( \alpha, t) } \left( \frac{1}{ \hatz_i (t, \alpha ) - E } - \frac{1}{ \hatgam_p (\alpha, t) - E } \right) \rho_t (E + E_- (0, t), 0 ) \right) + C \frac{ N^{\om_\ell + \delta/3}}{N^{1/3} }.
\eeq
It is easy to see that the first term is bounded by
\begin{align}
 &\left( \int_{ E \geq \hatgam_{2p} ( \alpha, t) } \left( \frac{1}{ \hatz_i (t, \alpha ) - E } - \frac{1}{ \hatgam_p (\alpha, t) - E } \right)\rho_t (E + E_- (0, t), 0 ) \right) \notag\\
  \leq & C |\tilz_i ( \alpha, t) - \hatgam_p ( \alpha, t) | \frac{ \Im [ m_t ( E_{-} (0, t) + \i \hatgam_p ( \alpha, t) , 0 ) ]}{ \hatgam_p (\alpha, t) } \leq C \sqrt{ \hatgam_p ( \alpha, t) } \leq C \frac{ N^{\om_\ell + \delta/3}}{N^{1/3}}.
\end{align}
Collecting the above, we see that, using $\om_1 \leq \om_\ell/2$, and $\om_\ell < \om_0/10$,
\beq
\sup_{0 \leq s \leq \tau} F(s) \leq  \left( \frac{ \nu^2 N^{ 4\om_\ell +2 \delta/3}}{N^{4/3}} + \frac{ \nu N^{\om_\ell+\om_1 + \delta/3}}{N^{2/3}} \right) + F(0).
\eeq
under the condition \eqref{eqn:nucond}.  Hence if we choose
\beq
\nu = \frac{ N^{2/3} }{ N^{2 \om_\ell + 2 \delta/3}} N^{\eps_1}
\eeq
for $\eps_1 < \delta /10$, then we see by continuity that $\tau = 10 t_1$  with overwhelming probability, and so 
\beq
\sup_{ 0 \leq s \leq 10 t_1} F(s) \leq 5
\eeq 
with overwhelming probability.  Now if $i \leq N^{ 3 \om_\ell+\delta}/2$, we see that
\beq
\nu | \hatz_i (t, \alpha ) - \gamhat_p (t, \alpha ) | \geq c N^{\eps_1}
\eeq
with overwhelming probability.

This yields,
\beq
\UL_{i p_*} (0, t) \leq N^{-D}
\eeq
for such $i, p_*$ with overwhelming probabilty.  The proof for $\UL_{ p_* i}$ is the same; instead, use $\psi  \to - \psi$ and set the initial condition to be a $f_j = \delta_{ij}$.
\qed

It will be convenient to establish the above estimate holding simultaneously for all $0 \leq s \leq t \leq 10 t_1$.  We will require the following.
\bel
Let $u_i$ be a solution of 
\beq
\del_t u = \L u,
\eeq
with $u_i (0) \geq 0$.  Then for $0 \leq t \leq 10 t_1$ we have
\beq
\frac{1}{2} \sum_i u_i (0) \leq \sum_i u_i (t) \leq \sum_i u_i (0)
\eeq
with overwhelming probability.
\eel
\proof We see that
\beq
\del_t \sum_i u_i =  \sum_i \V_i u_i.
\eeq
With overwhelming probability, 
\beq
- \V_i \leq C \frac{N^{1/3}}{N^{\om_\ell}}.
\eeq
The claim follows from applying Gronwall to
\beq
\del_t \sum u_i \geq - \left( C \frac{N^{1/3}}{N^{\om_\ell}} \right) \sum_i u_i.
\eeq
\qed

\bel \label{lem:fstech} Let $\delta, \eps >0$.  Let $a \leq N^{3 \om_\ell+\delta}/2$ and $b  \geq  N^{3 \om_\ell + \delta+\eps}$.  Then,
\beq
\sup_{ 0 \leq s \leq t \leq 10 t_1 } \UL_{ab} (s, t) + \UL_{ba} (s, t) \leq N^{-D}
\eeq
with overwhelming probability.
\eel
\proof We have,
\beq
\UL_{ai } (0, t) \geq \UL_{ab} (s, t) \UL_{bi} (0, s).
\eeq
By the previous lemma, we have that $\sum_i \UL_{bi} (0, s) \geq 1/2$.  By the first finite speed estimate we know that $\UL_{bi} \leq N^{-100}$ for any $i \leq N^{3 \om_\ell+\delta+\eps}/2$.  Hence there is an $i_* \geq N^{3 \om_\ell+\delta}$ so that $\UL_{b i_*} (0, s) \geq 1/(4 N)$.  But then also by the first finite speed estimate,
\beq
\UL_{ai_*} (0, t) \leq N^{-D}.
\eeq
Hence we get $\UL_{ab} (s, t) \leq N^{-D+2}$.  The estimate for $\UL_{ba} (s, t)$ is similar. \qed

\section{Local law for $t \geq N^{1/3}$} \label{sec:ll1}
In this section we are going to prove a local law for $H_t$.   We will use notation introduced in Section \ref{sec:fcnot}, i.e., $\kappa$, $\xi$, etc. 

Let us denote the matrix elements of $H_t$ by
\beq
(H_t)_{ij} = V_i \delta_{ij} + \sqrt{t} h_{ij}.
\eeq

We will also use the notation
\beq
t = \frac{ N^{\om}}{N^{1/3}}, \qquad \om >0.
\eeq

Fix a $\sigma >0$. Define the domain
\beq
\D_{\sigma} := \left\{ E + \i \eta : 3/4 \geq E \geq E_-,  \sqrt{ \kappa + \eta } \geq \frac{ N^{\sigma}}{N \eta } \right\} \bigcup \left\{ E + \i \eta : -3/4 \leq E \leq E_-, \eta \geq N^{\sigma} / N^{2/3} \right\}.
\eeq
The main theorem of this section is the following.  The derivation of rigidity estimates such as \eqref{eqn:xirig1} from such an estimate is standard - we refer, to, e.g., \cite{burger,general}.
\bet \label{thm:ll}
Let $\sigma >0$.  For any $\eps >0$, we have with overwhelming probability that the following estimates hold for all $z \in \D_\sigma$.  First, for $E \leq E_-$ we have
\beq
| m_N (z) - \mfct (z) | \leq N^{\eps} \left( \frac{1}{ N ( \kappa + \eta ) } + \frac{1}{ ( N \eta)^2 \sqrt{ \kappa + \eta } } \right).
\eeq
For $E \geq E_-$ we have,
\beq
| m_N (z) - \mfct (z) | \leq \frac{ N^{\eps}}{N \eta}.
\eeq
\eet
We introduce now some notation and estimates used in the proof.  
The Schur complement formula gives,
\beq \label{eqn:schur}
G_{ii} (z) := \frac{1}{ V_i - z - t \mfct + Z_i }
\eeq
where
\begin{align}
Z_i =& t (m_N - \mfct ) + \sqrt{t} h_{ii} + t A_i + t B_i + t ( m_N^{(i)} - m_N ) \notag\\
&= t ( m_N - \mfct ) + t ( m_N^{(i)} - m_N ) + Q_i ( G_{ii}^{-1} ).
\end{align}
Here we defined,
\beq
A_i = \sum_{j \neq k } h_{ij} G_{jk}^{(i)} h_{ki}, \qquad B_i = \sum_j ( h^2_{ij} - 1/N) G_{jj}^{(i)}.
\eeq
Let
\beq
\Lambda := | m_N - \mfct |, \qquad \Phi := \sqrt{ \frac{t}{N} } + t \sqrt{ \frac{ \Im [ \mfct ] + \Lambda}{N \eta }}.
\eeq 
We have the following esimates for the above parameters.  They are standard, see, e.g., \cite{landonyau}.
\bel \label{lem:ld}
We have,
\beq \label{eqn:mNimN}
t |m_N^{(i)} - m_N | \leq \frac{t}{ N \eta } \frac{ \Im [ G_{ii} ] }{ |G_{ii} | } \leq \frac{t}{N \eta } 
\eeq
For any $\eps >0$ we have with overwhelming probability that
\beq
t |A_i| + t |B_i| \leq N^{\eps} t \left( \frac{ \Im [ m_N^{(i)} ] }{ N \eta } \right)^{1/2} \leq C N^{\eps} \left( \frac{t}{N \eta } + t \sqrt{ \frac{ \Lambda + \Im [ \mfct ] }{ N \eta } } \right).
\eeq
We have,
\beq \label{eqn:Qi}
| Q_i ( G_{ii}^{-1} ) | \leq N^{\eps} \left( \sqrt{ \frac{t}{N} } +t \sqrt{ \frac{ \Lambda + \Im [ \mfct] }{ N \eta } } + \frac{t}{ N \eta } \right)
\eeq
\eel

Due to how some quantities appearing in the self-consistent equation for $m_N$ behave, it will be notationally convenient to split the proof of the above theorem into two parts.  We will first consider $z \in \D_\sigma$ s.t. $E \geq E_- - C_1 t^2$ for any (fixed) $C_1 >0$.  After this we sketch how the proof is modifed to deal with the remaining part of $\D_\sigma$.

\subsection{Proof for $E \geq E_- - t^2$}
In this section we prove the following.
\bep Fix $C_1>0$.   \label{prop:ll1}
Theorem \ref{thm:ll} holds in the domain
\beq
\D_1 := \D_\sigma \cap \left\{ z = E+ \i \eta : E \geq E_- - C_1 t^2 \right\}.
\eeq
\eep
In preparation we note the behaviour of a few parameters appearing in the proof.  The proofs are provided in Lemma \ref{lem:llcoef1}.  First, we note that for $z \in \D_1$, that there is a $c>0$ so that
\beq \label{eqn:Vistab}
|V_i - \xi | \geq c (t^2 + \eta + t\Im [ \mfct ] ) \asymp t^2 + \eta + t |E - E_- |^{1/2} \1_{ \{ E \geq E_- \} }.
\eeq  
For $z \in \D_1$ we have,
\beq \label{eqn:gipbd}
\frac{1}{N} \sum_i |g_i|^p \leq C \frac{ t+ \eta^{1/2} + \Im [ \mfct ] }{ ( t^2 + \Im [ \xi ] )^{p-1} } \leq C \frac{ t + \Im [ \mfct ] }{ (t^2 + \Im [ \xi ] )^{p-1}}
\eeq
The following is an immediate consequence of Lemma \ref{lem:ld} and the definition of the spectral domain $\D_1$.  Its role is to provide a sufficient condition under which we can Taylor expand the Schur complement formula for $G_{ii}$.
\bel \label{lem:expand}
On the event
\beq
\Lambda \leq \frac{ \Im [ \mfct ] + t}{ ( \log(N) )^2 },
\eeq
we have for any $\eps >0$,
\beq \label{eqn:Ziexpand}
|Z_i| \leq N^{\eps}\left( t^2 ( N^{-3 \om/2} + N^{-\sigma/2} ) + t \Im [\mfct] N^{-\sigma/2} \right) + \frac{ \Im [ \mfct ] t+ t^2}{ \log(N)^2}.
\eeq
and
\beq \label{eqn:Gii}
\frac{1}{2} |g_i | \leq | G_{ii} | \leq 2 |g_i|,
\eeq
with overwhelming probability.  
\eel
\subsubsection{Self-consistent equation}
In this subsection we derive the following self-consistent equation.  Many arguments are similar to those appearing in \cite{landonyau}.
\bep \label{prop:self}
On the event
\beq
\Lambda \leq \frac{ \Im [ \mfct ] + t }{ ( \log (N) )^2 }
\eeq
we have, for any $\eps>0$ with overwhelming probability,
\begin{align}
 &\left| ( 1- t R_2 ) ( m_N - \mfct ) + t^2 R_3 ( m_N - \mfct )^2 \right| \notag\\
\leq & \left| \frac{1}{N} \sum_i g_i^2 Q_i ( G_{ii}^{-1} ) \right| + N^{\eps} \left( \frac{\Lambda + \Im [ \mfct] }{N \eta } + \frac{1}{ ( N\eta)^2 } \right) \frac{t}{t^2 + \Im [ \xi ] } + \frac{\Lambda^2}{ \log(N)} \frac{ t}{ t^2+ \Im [ \xi ] }  \label{eqn:strongself} \\ 
\leq &  N^{\eps} \left( \frac{1}{N \eta } +  \sqrt{ \frac{ \Lambda + \Im [ \mfct ] }{N \eta }} \right) \frac{ t^2 + t \Im [ \mfct ]  }{t^2 + \Im [ \xi ] } +  \frac{\Lambda^2}{ \log(N)} \frac{ t}{ t^2+ \Im [ \xi ] }  \label{eqn:weakself}
\end{align}
\eep
\proof Due to Lemma \ref{lem:expand}, we may, with overwhelming probability, Taylor expand the Schur complement formula \eqref{eqn:schur} in powers of $Z_i$.  We arrive at
\beq
m_N = \mfct + \frac{1}{N} \sum_i g_i^2 Z_i + \frac{1}{N} \sum_i g_i^3 (Z_i)^2 + \O \left( \frac{1}{N} \sum_i |g_i|^4 |Z_i|^3 \right).
\eeq
We split the first order term as
\beq
\frac{1}{N} \sum_i g_i^2 Z_i = \frac{1}{N} \sum_i t g_i^2 ( m_N - \mfct ) + \frac{1}{N} \sum_i g_i^2 Q_i ( G_{ii}^{-1} ) + \frac{1}{N} \sum_i g_i^2 t ( m_N^{(i)} - m_N ).
\eeq
Using \eqref{eqn:mNimN},  \eqref{eqn:Gii}  and \eqref{eqn:Vistab} we obtain,
\begin{align}
\left|  \frac{1}{N} \sum_i g_i^2 t ( m_N^{(i)} - m_N ) \right| &\leq C \frac{1}{N^2 \eta} \sum_i t |g_i| \Im [ G_{ii} ]  \notag\\
&\leq C \frac{\Lambda + \Im [ \mfct ] }{N \eta} \frac{t}{t^2 + \Im [ \xi ] }.
\end{align}
For the other term we write
\beq \label{eqn:fll1}
\frac{1}{N} \sum_i g_i^2 Q_i [ G_{ii}^{-1} ] = \frac{1}{N} \sum_i g_i^2 \sqrt{t} h_{ii} + \frac{1}{N} \sum_i g_i^2 t (A_i + B_i ).
\eeq
For the first term we just calculate the variance and find that with overwhelming probability,
\begin{align}
\left|  \frac{1}{N} \sum_i g_i^2 \sqrt{t} h_{ii}  \right| &\leq N^{\eps} \left( \frac{1}{N^3} \sum_i t |g_i|^4 \right)^{1/2} \leq C N^{\eps} \frac{(t^2 + t \Im [ \mfct ] )^{1/2}}{N ( t^2 + \Im [ \xi ] )^{3/2} } \notag\\
&\leq C N^{\eps} \left( \frac{t}{t^2  + \Im [ \xi ] } \frac{ \Im [ \mfct] }{ N \eta } + \frac{(t \Im [ \mfct ] )^{1/2} }{N ( t^2 + \Im [ \xi ] )^{3/2} } \right),
\end{align}
where we used \eqref{eqn:gipbd}.  
For the second sum of \eqref{eqn:fll1}, we estimate, again using \eqref{eqn:gipbd},
\beq
\left| \frac{1}{N} \sum_i g_i^2 t (A_i + B_i ) \right| \leq N^{\eps} \left( \frac{1}{N \eta } +  \sqrt{ \frac{ \Lambda + \Im [ \mfct ] }{N \eta }} \right) \frac{ t^2 + t \Im [ \mfct ]  }{t^2 + \Im [ \xi ] }.
\eeq
For the second order term we have by Cauchy-Schwarz
\begin{align}
\frac{1}{N} \sum_i g_i^3 Z_i^2 &= \frac{t^2}{N} \sum_i g_i^3 ( m_N - \mfct ) \notag\\
&+ \O \left( \frac{1}{N} \sum_i t^2 |g_i|^3 \frac{ \Lambda^2}{ \log(N)^2} + \log(N)^2 \frac{1}{N} \sum_i |g_i|^3 ( | Q_i ( G_{ii}^{-1} ) |^2 + t^2 |m_N - m_N^{(i)} |^2 )  \right)
\end{align}
We estimate,
\beq
\frac{1}{N} \sum_i t^2 |g_i|^3 \frac{ \Lambda^2}{ \log(N)^2}  \leq C \frac{ \Lambda^2}{ \log(N)^2} \frac{t}{ t^2 + \Im [ \xi ] }.
\eeq
On $\D_1$ we have,
\beq
t |g_i| |m_N - m_N^{(i)} | \leq \frac{t}{(N \eta ) ( t^2 + t \Im [ \xi ] )} \leq 1
\eeq
and so,
\beq
\frac{1}{N} \sum_i |g_i|^3 t^2 |m_N - m_N^{(i)} |^2\leq \frac{1}{N} \sum_i |g_i|^2 t |m_N - m_N^{(i)} | \leq C \frac{ \Lambda + \Im [ \mfct ] }{N \eta} \frac{t}{t^2 + \Im [ \xi ] }.
\eeq
Next, using \eqref{eqn:Qi} we find with overwhelming probability,
\begin{align}
\frac{1}{N} \sum_i | g_i|^3 |Q_i ( G_{ii} ) |^2 &\leq N^{\eps} \left(  \frac{ t}{N} +\frac{t^2}{ (N \eta )^2} + \frac{t^2 ( \Lambda + \Im [ \mfct ] )}{ N \eta } \right) \frac{ t + \Im [ \mfct ] }{ (t^2 + \Im [ \xi ] )^2} \notag\\
&\leq N^{\eps} \left( \frac{t^2  }{ N (t^2 + \Im [ \xi ] )^2 } + \left[ \frac{ \Lambda + \Im [ \mfct ] }{ N \eta } + \frac{1}{ (N \eta)^2} \right] \frac{t}{ t^2 + \Im [ \xi ] } \right) \notag\\
&\leq CN^{\eps} \left[ \frac{ \Lambda + \Im [ \mfct ] }{ N \eta } + \frac{1}{ (N \eta)^2} \right] \frac{t}{ t^2 + \Im [ \xi ] } .
\end{align}
In the last line we used the fact that $\Im [ \mfct ] / \eta \geq t / (t^2 + \Im [ \xi ] )$ on $\D_1$.   Since we have $|Z_i | |g_i | \leq C / (\log(N))^2$ with overwhelming probability, we easily obtain,
\begin{align}
\frac{1}{N} \sum_i |g_i|^4 |Z_i|^3 &\leq \frac{C}{N} \sum_i t^2 |g_i|^3 \frac{ \Lambda^2}{ ( \log(N))^2 } + \frac{1}{N} \sum_i |g_i|^3 ( |Q_i [ G_{ii}^{-1} ] |^2 + t^2 |m_N - m_N^{(i)} |^2 )  \notag\\
&\leq  C \frac{\Lambda^2}{ \log(N)^2} \frac{t}{ t^2 + \Im [ \xi ] }+ N^{\eps} \left[ \frac{ \Lambda + \Im [ \mfct ] }{N \eta} + \frac{1}{ (N \eta)^2} \right] \frac{t}{t^2 + \Im [ \xi ] }.
\end{align}
The claim follows. \qed
\subsubsection{Weak local law}

Before establishing the optimal estimate of Proposition \ref{prop:ll1}, we establish the following weaker estimate.
\bep \label{prop:wll}
Let $\eps >0$.  With overwhelming probability we have the following estimates for every $z \in \D_1$.  For $\kappa + \eta \geq t^2$, 
\beq \label{eqn:wllest1}
| m_N - \mfct | \leq N^{\eps} \sqrt{ \frac{t + \Im [ \mfct ] }{ N \eta } }.
\eeq
For $\kappa + \eta \leq t^2$,
\beq \label{eqn:wllest2}
| m_N - \mfct | \leq N^{\eps} \frac{ t^{2/3}}{( N \eta )^{1/3} }.
\eeq
\eep
We start with the following.
\bep \label{prop:weakest} Fix $\eps >0$. 
The following holds with overwhelming probability on the event
\beq \label{eqn:lamest2}
\Lambda \leq \frac{ t + \Im [ \mfct ] }{ ( \log (N) )^2}.
\eeq
If $\kappa + \eta \geq t^2$, then
\beq
\Lambda \leq 2 N^{\eps} \sqrt{ \frac{ t + \Im [ \mfct ] }{ N \eta } } \leq C  N^{\eps} N^{-\sigma/2} ( t + \Im [ \mfct ] ). \label{eqn:lamest1}
\eeq
If $\kappa + \eta \leq t^2$ then we have the following dichotomy.  Either,
\beq
\Lambda \geq c \sqrt{ \kappa + \eta },
\eeq
or
\beq
\Lambda \leq N^{\eps} \left( \frac{t}{ \sqrt{ \kappa + \eta } } \sqrt{ \frac{ \Im [ \mfct]}{ N \eta } }  + \frac{t^2}{ \kappa + \eta} \frac{1}{ N \eta } \right).
\eeq
If $\kappa + \eta \leq t^2$ then,
\beq \label{eqn:wll5}
\Lambda \leq N^{\eps} \frac{t^{2/3}}{ ( N \eta )^{1/3} } + ( \kappa + \eta )^{1/2}  + N^{\eps} t^{1/2} \left( \frac{ \Im [ \mfct ]^{1/4} }{ ( N \eta)^{1/4} } + \frac{1}{ ( N \eta)^{1/2} } \right).
\eeq
\eep
\proof First suppose that $\kappa+ \eta \geq t^2$.  From \eqref{eqn:weakself} we get the estimate
\beq
\Lambda \leq \frac{ \Lambda}{2} + N^{\eps} \sqrt{ \frac{ t + \Im [ \mfct ] }{ N \eta } },
\eeq
and so
\beq
\Lambda \leq 2 N^{\eps} \sqrt{ \frac{ t + \Im [ \mfct ] }{ N \eta } } \leq C  N^{\eps} N^{-\sigma/2} ( t + \Im [ \mfct ] ).
\eeq
Now we assume that $|E - E_- | + \eta \leq t^2$.  In this regime $|1-t R_2| \asymp \sqrt{ \kappa + \eta }/t$ and so we obtain from \eqref{eqn:weakself},
\begin{align}
\Lambda \leq  C \frac{ \Lambda^2}{ \sqrt{ \kappa + \eta }} + \frac{ t N^{\eps}}{ \sqrt{ \kappa + \eta } } \left( \frac{1}{ N \eta } + \sqrt{ \frac{ \Lambda + \Im [ \mfct ] }{ N \eta } } \right).
\end{align}
Hence, either
\beq
\Lambda \geq c \sqrt{ \kappa + \eta }
\eeq
or
\beq
\Lambda \leq N^{\eps} \left( \frac{t}{ \sqrt{ \kappa + \eta } }  \sqrt{ \frac{ \Im [ \mfct] }{N \eta } }  + \frac{t^2}{\kappa + \eta} \frac{1}{ ( N \eta ) } \right).
\eeq
The quadratic estimate from \eqref{eqn:weakself} gives,
\begin{align}
\Lambda^2 &\leq  t N^{\eps} \frac{ \Lambda^{1/2}}{ ( N \eta)^{1/2}} + \Lambda ( \kappa+ \eta )^{1/2} +C \frac{ \Lambda^2}{ \log(N)^2} \notag\\
&+ t N^{\eps}\left(  \sqrt{ \frac{ \Im [ \mfct ]}{N \eta } } + \frac{1}{N \eta} \right)
\end{align}
which implies
\beq
\Lambda \leq N^{\eps} \frac{t^{2/3}}{ ( N \eta )^{1/3} } +\sqrt{ \kappa + \eta } + N^{\eps} t^{1/2} \left( \frac{ \Im [ \mfct ]^{1/4} }{ ( N \eta)^{1/4} } + \frac{1}{ ( N \eta)^{1/2} } \right).
\eeq
The claim is proven.
\qed

\noindent{\bf Proof of Proposition \ref{prop:wll}.}   We follow the usual proof of such weak local laws, which is a continuity argument in $\eta$.  What we have to check is that the estimate we obtain at each scale on $\Lambda$ is much smaller than $t + \Im [ \mfct ]$.  

Fix an energy $E$.  The estimate for $\eta \geq 1$ is standard.  Fix then a sequence $\eta_k = 1 - k / N^2$, of cardinality less than $N^2$.  First we estimate for $\eta_k$ with $\kappa + \eta_k \geq t^2$.  By the continuity of $\Lambda$ in $z$ we see that the estimate \eqref{eqn:lamest1} at $\eta_{k}$ implies that \eqref{eqn:lamest2} holds at $\eta_{k+1}$.  Then by Proposition \ref{prop:weakest} that \eqref{eqn:lamest1} also holds at $\eta_{k+1}$ with overwhelming probability.  Hence we obtain \eqref{eqn:lamest1} for $\kappa + \eta \geq t^2$.


Now let us consider $\kappa + \eta_k \leq t^2$.   
Let $\eta^* =  \eta^* (E)$ be the first time that
\beq
\kappa + \eta = N^{a} \frac{ t^{4/3}}{ ( N \eta )^{2/3}}.
\eeq
for a small $a>0$ that we choose later.  

We may suppose that the estimate \eqref{eqn:wllest2} holds at $\eta_{k-1}$ (if $\eta_{k-1}$ was s.t. $\kappa + \eta_{k-1} \geq t^2$ then instead \eqref{eqn:wllest1} holds, but this is estimate is the same order as \eqref{eqn:wllest1} at this scale).  On the domain $\D_1$ we have
\beq
\frac{1}{ N \eta} \leq N^{-\sigma/2} ( t + \sqrt{ \kappa + \eta } ) 
\eeq
and so we see that \eqref{eqn:lamest2} holds at $\eta_k$, as well as
\beq \label{eqn:wll4}
\Lambda \leq 2 \frac{N^{\eps} t^{2/3} }{ ( N \eta_k )^{1/3} } \leq 2 N^{\eps} N^{-a/2} \sqrt{ \kappa + \eta },
\eeq
by the definition of $\eta_*$.  
Proposition \ref{prop:weakest} implies that with overwhelming probability, either
\beq
\Lambda \geq  c \alpha
\eeq
or
\beq
\Lambda \leq N^{\delta/2} \left( \frac{t ( \Im [ \mfct ] )^{1/2} }{ \alpha^{1/2}  (N \eta )^{1/2} } + \frac{t^2}{ \alpha N \eta } \right) \leq N^{\delta} \frac{t^{2/3} }{ ( N \eta)^{1/3} }.
\eeq
for any $\delta >0$.   However, due to \eqref{eqn:wll4} the latter estimate holds.

By iteration, we obtain that \eqref{eqn:wllest2} holds for $\eta \geq \eta^*$ with overwhelming probability.
Now consider $\eta \leq \eta^*$.  Suppose that the estimate
\beq
\Lambda \leq   N^{a+\eps} \frac{t^{2/3} }{ ( N \eta )^{1/3} }
\eeq
holds at $\eta_{k-1}$, for any $\eps >0$ with overwhelming probability.  This implies that
\beq
\Lambda \leq N^{-\sigma/4} (t + \Im [ \mfct] ).
\eeq
provided we choose $a < \sigma / 10$, which we can do.  Then, this esimtae implies that \eqref{eqn:lamest2} holds at $\eta_k$.  Then, we see that \eqref{eqn:wll5} holds, which implies that,
\beq
\Lambda \leq   N^{a+\delta} \frac{t^{2/3} }{ ( N \eta )^{1/3} },
\eeq
with overwhelming probability, for any $\delta >0$, at $\eta_k$.    Hence, we obtain the claim on our spectral domain $\D_1$. \qed

\subsubsection{Fluctuation averaging lemma}
We have the following improved bound for one of the error terms appearing in the local law for $m_N$.  For a deterministic control parameter $\gamma$ we define
\beq
\Phi := \sqrt{ \frac{t}{N} } + t \sqrt{ \frac{ \Im [ \mfct ] + \gamma }{ N \eta } }.
\eeq
\bel \label{lem:fa1}  Fix $z \in \D_1$.  
Suppose that $\Lambda \leq \gamma$ with overwhelming probability, where $\gamma$ is a deterministic control parameter satisfying
\beq
\frac{1}{N \eta } \leq \gamma \leq \frac{ t+ \Im [ \mfct ] }{ \log(N)^2}.
\eeq
For any even $p\geq 2$ we have,
\beq \label{eqn:fa1}
\ee \left| \frac{1}{N} \sum_i g_i^2 Q_i [ G_{ii}^{-1} ] \right|^p \leq N^{\eps} \left( \frac{ 1}{ N \eta } \frac{ \Im [ \mfct  ] + \gamma}{ t + \Im [ \mfct ] + \gamma } \right)^p
\eeq
for any $\eps >0$.
\eel
\proof By \eqref{eqn:momentbd} and the estimates
\beq
|g_i| \leq \frac{C}{ t^2 + \Im [ \xi ] }, \quad \frac{1}{N} \sum_i |g_i|^2 \leq \frac{ t + \Im [ \mfct ] }{ t^2 + \Im [ \xi ] },
\eeq
we obtain
\begin{align}
&\ee \left| \frac{1}{N} \sum_i g_i^2 Q_i [ G_{ii}^{-1} ] \right|^p  \notag\\
\leq & \max_{0 \leq s \leq p } \max_{0 \leq l \leq (p+s)/2} \frac{1}{N^p} \left( \sqrt{ \frac{t}{N} } + t \sqrt{ \frac{ \Im [ \mfct ] + \gamma }{N \eta }} \right)^{p+s} \frac{1}{ ( \Im [ \xi ] + t^2 )^{2p + s } } ( N ( t^2 + \Im [ \xi ] ) ( t + \Im [\mfct ] ) )^l \notag\\
&\leq \max_{0 \leq s \leq p } \frac{1}{N^{p/2}} \left( \sqrt{ \frac{t}{N} } + t \sqrt{ \frac{ \Im [ \mfct ] + \gamma }{N \eta }} \right)^{p+s} \frac{N^{s/2} (t + \Im [ \mfct ])^{(p+s)/2}}{ ( \Im [ \xi] + t^2 )^{3p/2+s/2} }. \notag\\
&\leq C \frac{1}{ ( N ( \Im [ \xi ] + t^2 ) )^p } \notag\\
&+ \max_{0 \leq s \leq p } \left[ \frac{ t ( \Im [ \mfct ] + \gamma)^{1/2} (t + \Im [ \mfct ] )^{1/2} }{ N \eta^{1/2} ( \Im [ \xi ] + t^2 )^{3/2} } \right]^p \left[ \frac{ t ( \Im [ \mfct ] + \gamma)^{1/2} ( t + \Im [ \mfct ] )^{1/2} }{ \eta^{1/2} ( \Im [ \xi] + t^2 )^{1/2} } \right]^s.
\end{align}
The first term is bounded by
\beq
\frac{1}{ ( N ( \Im [ \xi ] + t^2 ) )^p }  \leq C \frac{1}{ ( N \eta )^p } \left(\frac{ \Im [ \mfct ] }{ t + \Im [ \mfct ] } \right)^p \leq C \frac{1}{ ( N \eta )^p } \left(\frac{ \Im [ \mfct ] + \gamma }{ t + \Im [ \mfct ] + \gamma } \right)^p.
\eeq
In the first inequality we used $\Im [ \mfct ] / \eta \geq c/t$ which holds on $\D_1$. 
For the second term, 
\begin{align}
&\max_{0 \leq s \leq p } \left[ \frac{ t ( \Im [ \mfct ] + \gamma)^{1/2} (t + \Im [ \mfct ] )^{1/2} }{ N \eta^{1/2} ( \Im [ \xi ] + t^2 )^{3/2} } \right]^p \left[ \frac{ t ( \Im [ \mfct ] + \gamma)^{1/2} ( t + \Im [ \mfct ] )^{1/2} }{ \eta^{1/2} ( \Im [ \xi] + t^2 )^{1/2} } \right]^s \notag\\
\leq &\left[ \frac{ t ( \Im [ \mfct ] + \gamma)^{1/2} (t + \Im [ \mfct ] )^{1/2} }{ N \eta^{1/2} ( \Im [ \xi ] + t^2 )^{3/2} } \right]^p + \left[ \frac{ t^2 ( \Im [ \mfct ] + \gamma) (t + \Im [ \mfct ] )}{ N \eta  ( \Im [ \xi ] + t^2 )^2 } \right]^p \notag\\
&\leq C \left( \frac{1}{ N \eta } \frac{ \Im [ \mfct ] + \gamma }{ t + \Im [ \mfct ] + \gamma } \right)^p.
\end{align}
This yields the claim. \qed

\subsubsection{Proof of Proposition \ref{prop:ll1}}

The proof is by an iteration procedure on $\gamma$.  We may assume that,
\beq
\Lambda \leq \gamma = \frac{ N^{\eps} t^{2/3}}{ (N \eta )^{1/3}},
\eeq
with overwhelming probability.
Let us first consider the case $\kappa + \eta \geq t^2$.   By \eqref{eqn:strongself} and \eqref{eqn:fa1} we have,
\begin{align}
\Lambda \leq \frac{1}{N \eta } + N^{\eps} \left( \frac{ \Lambda + \Im [ \mfct ] }{N \eta } + \frac{1}{ ( N \eta)^2} \right)\frac{1}{ t + \Im [ \mfct ] } + C \Lambda^2 \frac{1}{t+ \Im [ \mfct ] }
\end{align}
Using that $\Lambda \ll t + \Im [ \mfct ]/ \log(N)^2$ and $(N \eta)^{-1} \leq C(t + \Im [ \mfct ])$ on $\D_1$ we see that
\beq
\Lambda \leq \frac{ N^{\eps }}{N \eta}
\eeq
for any $\eps >0$ with overwhelming probability.  Now we consider $\kappa + \eta \leq t^2$.  We see that
 In this case we see that 
\begin{align} \label{eqn:sself}
| (1 - t R_2 ) ( m_N - \mfct ) + t^2 R_3 ( m_N - \mfct )^2 | \leq   \frac{N^{\eps}}{ N \eta } \frac{ \Im [ \mfct ] + \gamma}{t} + \frac{N^{\eps}}{ ( N \eta )^2} \frac{1}{t} + \frac{ \Lambda^2}{\log(N)^2 t}.
\end{align}
for any $\eps >0$ with overwhelming probability. 
Hence we get the estimates
\beq
\Lambda \leq C \frac{ \Lambda^2}{ \alpha } + N^{\eps} \frac{ \gamma + \Im[ \mfct ] }{ N \eta  \alpha } + \frac{N^{\eps}}{ ( N \eta)^2 \alpha }
\eeq
where we denoted
\beq
\alpha := ( \kappa + \eta )^{1/2}.
\eeq
Hence, we see that either
\beq
\Lambda \geq c \alpha
\eeq
or
\beq
\Lambda \leq  N^{\eps} \frac{ \gamma + \Im [ \mfct ] }{ N \eta \alpha } + \frac{N^{\eps}}{ (N \eta)^2 \alpha }.
\eeq
On the other hand if we estimate the quadratic term in \eqref{eqn:sself}, we see that
\beq
\Lambda^2 \leq C \alpha \Lambda + N^{\eps}  \frac{\Im [ \mfct ] + \gamma}{ N \eta } +N^{\eps} \frac{1}{ ( N \eta)^2}
\eeq
Let $\eta_*$ be the first time that 
\beq \label{eqn:etast}
 \frac{ \gamma + \Im [ \mfct ] }{ N \eta \alpha } + \frac{1}{ (N \eta)^2 \alpha } = N^{- 4 \eps} \frac{\alpha}{2}.
\eeq
By a similar continuity argument as in the proof of Proposition \ref{prop:wll} we get
\beq
\Lambda \leq  CN^{\eps}  \frac{ \gamma + \Im [ \mfct ] }{ N \eta \alpha } + C \frac{N^{\eps}}{ (N \eta)^2 \alpha } \leq  C N^{\eps} \sqrt{ \frac{ \gamma}{ N \eta } } + C N^{\eps} \frac{1}{ ( N \eta )}.
\eeq
for $\eta \geq \eta_*$, for any $\eps >0$ with overwhelming probability.  For $\eta < \eta_*$ we get
\beq
\Lambda \leq \frac{N^{4 \eps} }{N \eta }  N^{4 \eps}  \sqrt{ \frac{ \gamma}{ N \eta } } + N^{4 \eps} \sqrt{ \frac{ \Im [ \mfct ] }{ N \eta } } \leq \frac{N^{10 \eps}}{N \eta } + N^{10 \eps} \sqrt{ \frac{ \gamma}{ N \eta } } 
\eeq
where we used $\Im [ \mfct ] \leq C \alpha$ and then $\alpha^2 \leq N^{5 \eps}  / (N \eta )^2 + N^{5 \eps} \gamma / (N \eta)$ which follows from \eqref{eqn:etast} and using $\Im[\mfct] \leq C \alpha$ and the Schwarz inequality.  By iterating this argument finitely many times we see that we derive
\beq
| m_N - \mfct | \leq \frac{N^{\eps}}{ N \eta}
\eeq
with overwhelming probability for any $\eps >0$ and $z \in \D_1$.

For $E \leq E_-$, note that $\alpha \ll ( N \eta )^{-1}$ in $\D_1$.  Hence, we immediately derive from \eqref{eqn:sself} and the fact that $\gamma \leq N^{\eps} / (N \eta )$ that
\beq
\Lambda \leq \frac{N^{2 \eps}}{ ( N \eta )^2 \sqrt{ \kappa + \eta} } + \frac{ N^{2\eps}}{ N ( \kappa + \eta) }.
\eeq
This completes the proof.
\qed

\subsection{Proof for $E \leq E_- - t^2$}
We now complete the proof of Theorem \ref{thm:ll} by proving the following.  Let
\beq
\D_2 := \left\{ E + \i \eta: E \leq E_- - C_1 t^2 \right\} \cap \D_\sigma.
\eeq
\bep  \label{prop:ll2} The estimates of Theorem \ref{thm:ll} hold in $\D_2$.
\eep
As the proof is very similar to Proposition \ref{prop:ll1} we only give the important changes to the proof in the following subsections.

We will use the following a-priori estimates, which are used instead of \eqref{eqn:Vistab} and \eqref{eqn:gipbd}.  We have for $ z\in \D_1$, 
\beq \label{eqn:Vistab2}
|V_i - \xi | \geq c (t^2 + \kappa + \eta ),
\eeq
and
\beq \label{eqn:gipbd2}
\frac{1}{N} \sum_i |g_i |^p \leq \frac{C}{ ( \kappa + \eta + t^2)^{p-3/2} }.
\eeq
These follow from Lemma \ref{lem:llcoef2}.

\subsubsection{Self-consistent equation}
Instead of Lemma \ref{lem:expand} we have the following which is an easy consequence of Lemma \ref{lem:ld}
\bel \label{lem:expand2}
On the event
\beq
\Lambda \leq \frac{ t}{ \log(N)^2}
\eeq
we have with overwhelming probability,
\beq
|Z_i| \leq C \frac{t}{ \log(N)^2}
\eeq
and
\beq
\frac{1}{2} |g_i | \leq |G_{ii} | \leq 2 |g_i |.
\eeq
\eel
In the place of Proposition \ref{prop:self} we have,
\bep \label{prop:self2}
On the event 
\beq
\Lambda \leq \frac{ t}{ \log(N)^2},
\eeq
we have with overwhelming probability, for any $\eps >0$,
\begin{align}
&\left| (1 - t R_2 ) (m_N - \mfct ) \right| \leq C \frac{ \Lambda^2}{t} + N^{\eps} \left( \frac{1}{ ( N \eta )^2 ( t + \sqrt{ \kappa + \eta } )} + \frac{1}{ N ( \kappa + \eta + t^2 ) } \right) + \left| \frac{1}{N} \sum_i  g_i^2 Q_i [ G_{ii}^{-1} ] \right| \notag\\
\leq & C \frac{ \Lambda^2}{t} + N^{\eps} \left( \frac{1}{ N^{1/2} ( t + \sqrt{ \eta + \kappa } )^{1/2}} + \sqrt{ \frac{ \Im [ \mfct ] + \Lambda }{N \eta } }  + \frac{1}{ N \eta }\right).
\end{align}
\eep
Its proof is identical to that of Proposition \ref{prop:self}, except that we do not need to keep the term which is second order in $\Lambda^2$ as we always have $|1 - t R_2 | \asymp 1$ on $\D_2$.

\subsubsection{Weak local law}
From Proposition \ref{prop:self2} we immediately see that we have with overwhelming probability on the event $\Lambda \leq t / \log(N)^2$ that
\beq
\Lambda \leq \frac{N^{\eps}}{N^{1/2} ( \kappa + \eta + t^2 )^{1/4} }  + \frac{N^{\eps}}{ N \eta} 
\eeq
for any $\eps >0$.  The estimate on the RHS is $\ll t$, and so by a similar proof to Proposition \ref{prop:wll} we get
\bep \label{prop:wll2} Let $\eps >0$.  With overwhelming probability,
\beq
\Lambda \leq \frac{N^{\eps}}{N^{1/2} ( \kappa + \eta + t^2 )^{1/4} }  + \frac{N^{\eps}}{ N \eta} ,
\eeq
on $\D_2$.
\eep

\subsubsection{Fluctuation averaging lemma}
Analogously to Lemma \ref{lem:fa1} we have,
\bel
\label{lem:fa2} Fix $z \in \D_2$.   Suppose that $\Lambda \leq \gamma$ with overwhelming probability, where $\gamma$ is a deterministic control parameter satisfying
\beq
\frac{1}{N \eta } \leq \gamma \leq \frac{ t+ \Im [ \mfct ] }{ \log(N)^2}.
\eeq
For any even $p\geq 2$ we have,
\beq \label{eqn:fa2}
\ee \left| \frac{1}{N} \sum_i g_i^2 Q_i [ G_{ii}^{-1} ] \right|^p \leq N^{\eps}   \left( \frac{1}{ N ( t^2 + \kappa + \eta ) }  + \frac{ \gamma}{ N \eta (t + \sqrt{ \kappa + \eta } )}   \right)^p,
\eeq
for any $\eps >0$.
\eel
\proof We proceed as in the proof of Lemma \ref{lem:fa1}.  Applying \eqref{eqn:momentbd} and \eqref{eqn:Vistab2}, \eqref{eqn:gipbd2} we obtain,
\begin{align}
&\ee \left| \frac{1}{N} \sum_i g_i^2 Q_i [ G_{ii}^{-1} ] \right|^p \notag\\
\leq &\max_{0 \leq s \leq p } \max_{ 0 \leq l \leq (p+s)/2} N^{\eps}\left( \sqrt{ \frac{t}{N} } + t \sqrt{ \frac{ \Im [ \mfct ] + \gamma }{ N \eta } } \right)^{p+s} \frac{1}{ (t^2 + \eta + \kappa)^{s+2p-3l/2}} \frac{N^l}{N^p} \notag\\
\leq & \max_{0 \leq s \leq p }  C N^{\eps} \frac{N^{s/2}}{N^{p/2}} \left( \sqrt{ \frac{t}{N} } + t \sqrt{ \frac{ \Im [ \mfct ] + \gamma }{ N \eta } } \right)^{p+s} \frac{1}{ ( t + \sqrt{ \kappa + \eta })^{s/2+5p/2} } \notag\\
\leq &  CN^{\eps} \frac{1}{ ( N ( \kappa + \eta ) )^p } + \max_{0 \leq s \leq p } C N^{\eps} \left[ \frac{ t ( \Im [ \mfct ] + \gamma )^{1/2}}{N \eta^{1/2} (t +\sqrt{\kappa + \eta } )^{5/2} } \right]^p \left[ \frac{ t ( \Im [ \mfct ] + \gamma)^{1/2}}{ \eta^{1/2} (t + \sqrt{ \kappa + \eta } )^{1/2} } \right]^s.
\end{align}
For the second term, we estimate
\begin{align}
 & \max_{0 \leq s \leq p }  \left[ \frac{ t ( \Im [ \mfct ] + \gamma )^{1/2}}{N \eta^{1/2} (t +\sqrt{\kappa + \eta } )^{5/2} } \right]^p \left[ \frac{ t ( \Im [ \mfct ] + \gamma)^{1/2}}{ \eta^{1/2} (t + \sqrt{ \kappa + \eta } )^{1/2} } \right]^s \notag\\
 \leq & \left[ \frac{ ( \Im [ \mfct ] + \gamma)^{1/2} }{ N \eta^{1/2} ( t + \sqrt{ \kappa + \eta } )^{3/2} } \right]^p + \left[ \frac{ \Im [ \mfct ] + \gamma }{ N \eta \sqrt{ \kappa + \eta }  } \right]^p \leq C \left[ \frac{ 1}{ N ( \kappa + \eta ) } + \frac{ \gamma}{ N \eta  \sqrt{ \kappa + \eta }  } \right]^p
\end{align}
which yields the claim. \qed

\subsubsection{Proof of Proposition \ref{prop:ll2}}
Let $\Lambda \leq \gamma$ with $\gamma$ as in Lemma \ref{lem:fa2}.  From Lemma \ref{lem:fa2} and Proposition \ref{prop:self2}, we see that
\beq
\Lambda \leq N^{\eps} \left( \frac{1}{ N ( \kappa + \eta ) } + \frac{ 1}{ ( N \eta )^2 \sqrt{ \kappa + \eta } } + \frac{ \gamma}{ ( N \eta )  \sqrt{ \kappa + \eta } } \right),
\eeq
with overwhelming probability, for any $\eps >0$.  Iterating this, we see that $\Lambda \leq N^{\eps} / (N \eta )$ with overwhelming probability on $\D_2$.  Then, applying the above estimate again with $\gamma = N^{\eps} / (N \eta )$ we conclude the Proposition. \qed

\section{Local law for $0 \leq t \leq N^{1/3}$ and regular initial data} \label{sec:ll2}
In this section we want to prove a local law for 
\beq
H_t := V + \sqrt{t} G
\eeq
in the regime $0 \leq t \leq N^{-\eps}$, provided that $V$ already obeys a local law.  That is, we consider $V$ such that for any $\eps>0$ and $\sigma >0$, $V$ obeys the estimates
\beq
|m_V  - \hatm | \leq \frac{N^{\eps}}{N \eta }, \qquad 0 \leq E \leq 1, \quad 10\geq \sqrt{ |E| + \eta } \geq \frac{ N^{\sigma}}{N \eta },
\eeq
and
\beq
|m_V - \hatm | \leq N^{\eps} \left( \frac{1}{ N( |E| + \eta ) } + \frac{1}{ ( N \eta )^2 \sqrt{ |E| + \eta } } \right) , \qquad -1 \leq E \leq 0, \quad 10 \geq \eta \geq N^{\sigma-2/3},
\eeq
where $\hatm$ is the Stieltjes transform of a law $\hatrho(x)$ so that for $|x| \leq 1$, we have
\beq
\hatrho(x) \asymp \1_{ \{ x \geq 0 \} } \sqrt{x}.
\eeq
We will denote the free convolution of the semicircle at time $t$ with $\hatrho$ by $\hatrhot$ and its Stieltjes transform by $\hatmt$.  By Section \ref{sec:fcanal}, $\hatrhot$ behaves like a square root, and we denote the edge by $\hatE_-$.  We will abuse notation slightly and denote
\beq
\kappa := |E-\hatE_- |.
\eeq

  We want to prove the following theorem.  Let $\D_\sigma$ be as in Section \ref{sec:ll1}.  
\bet \label{thm:ll2}
For any $\sigma >0$, $V$ as above, $\eps >0$ and $\eps_1 >0$, the following estimates hold with overwhelming probability, for any $0 \leq t \leq N^{- \eps_1}$.  First, 
\beq
|m_N - \hatmt | \leq \frac{N^{\eps}}{N \eta}
\eeq
for $ z \in \D_\sigma$ and $E \geq \hatE_-$.    For $E \leq \hatE_-$ we have
\beq
| m_N - \hatmt | \leq N^{\eps} \left( \frac{1}{ N ( \kappa + \eta ) } + \frac{1}{ ( N \eta)^2 \sqrt{ \kappa + \eta } } \right).
\eeq
\eet
Let $\sigma >0$.  Suppose that we want to prove the above result on $\D_{\sigma}$.  Then in the regime
\beq
N^{-\eps } \geq t \geq \frac{ N^{\sigma/100}}{N^{1/3}},
\eeq
we know that Theorem \ref{thm:ll} holds, except that $m_N$ is close to $\mfct$ which is the free convolution of $V$ and the semicircle distribution, and not $\hatmt$, which is the free convolution of $\hatrho$ and the semicircle.  From Appendix \ref{a:fccont} we see that the difference $\hatmt - \mfct$ obeys the stated estimates of Theorem \ref{thm:ll2}.  Hence, we only need to prove Theorem \ref{thm:ll2} in the case that
\beq
0 \leq t \leq \frac{ N^{\sigma/100}}{N^{1/3}}
\eeq
on the domain $\D_\sigma$.  This is the content of the remainder of Section \ref{sec:ll2}.

\subsection{Proof for $E \geq \hatE_- - N^{-2/3+\sigma}$}

Similarly to Section \ref{sec:ll1} it is useful to split the proof of Theorem \ref{thm:ll2} into two cases.  The first is the following.
\bep \label{prop:ll3}
Let $\D_{\sigma}$ as above.  Let $0 \leq t \leq N^{\sigma/100}/N$.  The estimates of Theorem \ref{thm:ll2} hold on
\beq
\hatD_1 := \D_{\sigma} \cap \left\{ E + \i \eta : E \geq \hatE_- - N^{-2/3+\sigma} \right\}.
\eeq
\eep
On $\hatD_1$ we use the following estimates.
\beq \label{eqn:Vistab3}
|V_i - \xi | \geq c ( \eta + t \sqrt{ \kappa + \eta } ),
\eeq
and
\beq \label{eqn:gipbd3}
\frac{1}{N} \sum_i |g_i |^p \leq C \frac{ \sqrt{ \kappa + \eta }}{ (\eta + t \sqrt{ \kappa + \eta } )^{p-1} }
\eeq
which follows from Section \ref{sec:shorttimeself}

The following plays the role of Lemma \ref{lem:expand}
\bel
On the event
\beq
\Lambda \leq \frac{ \Im [ \hatmt ] }{ \log(N)^2}
\eeq
we have
\beq
|Z_i | \leq C \frac{  \eta + t \sqrt{\kappa + \eta } }{ \log(N)^2}
\eeq
and
\beq
\frac{1}{2} |g_i | \leq |G_{ii}| \leq 2 |g_i|
\eeq
with overwhelming probability.
\eel
\proof We just need to check the statement on $Z_i$ (for the $G_{ii}$ statement recall \eqref{eqn:Vistab3}).  Recall $Z_i= t (m_N - \hatmt ) + t ( m_N^{(i)} - m_N ) + t (A_i + B_i)  + \sqrt{t} h_{ii}$.  We only need to estimate $\sqrt{t} h_{ii}$, as the bounds for the other terms are immediate from Lemma \ref{lem:ld} and $(N \eta ) \leq N^{ - \sigma } \sqrt{ \kappa + \eta }$.  With overwhelming probability,
\beq
\sqrt{t} h_{ii} \leq \frac{ N^{\eps} t^{1/2}}{N^{1/2} } \leq N^{3 \eps } \frac{t}{ N \eta } + \frac{ \eta}{ \log(N)^2}
\eeq
and the claim follows. \qed
\subsubsection{Self consistent equation}
We have the following.
\bep \label{prop:self3}
On the event
\beq
\Lambda \leq \frac{ \Im [ \hatmt ]}{ \log(N)^2}
\eeq
we have with overwhelming probability, for any $\eps >0$,
\begin{align}
| (1-t R_2 ) ( m_N - \hatmt ) | \leq & C \frac{ \Lambda^2}{ \sqrt{ \kappa+ \eta }} + \left| \frac{1}{N} \sum_i g_i^2 Q_i [ G_{ii}^{-1} ] \right| + \frac{ N^{\eps}}{N \eta} \notag\\
\leq& C \frac{ \Lambda^2}{ \sqrt{ \kappa+ \eta }}  + C \frac{ N^{\eps}}{N \eta } +N^{\eps} \sqrt{ \frac{ \Im [ \hatmt] }{N \eta } }.
\end{align}
\eep
\proof We can write,
\beq
m_N - \hatmt = \frac{1}{V_i - \xi } - \hatm ( \xi ) + \frac{1}{N} \sum_i g_i^2 Z_i + \O\left( \frac{N^{\eps}}{N \eta } \right).
\eeq
Note that by assumption,
\beq
\left|  \frac{1}{V_i - \xi } - \hatm ( \xi )  \right| \leq \frac{N^{\eps}}{N \eta }.
\eeq
For the remaining term we have,
\beq
\frac{1}{N} \sum_i g_i^2 Z_i = \frac{t}{N} \sum_i g_i^2 ( m_N - \hatmt ) + \frac{1}{N} \sum_i g_i^2 Q_i [ G_{ii}^{-1} ] + \O \left( \frac{1}{N \eta } \right).
\eeq
For the last term we write,
\beq
\frac{1}{N} \sum_i g_i^2 Q_i [G_{ii}^{-1} ] = \frac{1}{N} \sum_i g_i^2 \sqrt{t} h_{ii} + \frac{1}{N} \sum_i g_i^2 t (A_i + B_i).
\eeq
By a variance calculation the first term is $\O ( N^{\eps} / (N \eta ))$ for any $\eps >0$.  The second term is bounded by
\beq
\left| \frac{1}{N} \sum_i g_i^2 t (A_i + B_i) \right| \leq N^{\eps} \sqrt{ \frac{ \Im [ \hatmt ] + \Lambda }{N \eta } } \frac{1}{N} \sum_i |g_i|^2 t \leq C N^{\eps}\sqrt{ \frac{ \Im [ \hatmt ]}{N \eta } }
\eeq
using Lemma \ref{lem:ld} and \eqref{eqn:gipbd3}.  \qed

\subsubsection{Weak local law}
From Proposition \ref{prop:self3} and the fact that $|1 - t R_2 | \asymp 1$ on $\D_\sigma$ (due to $t|R_2| \leq Ct / \sqrt{ \kappa + \eta } \ll 1$), we see that with overwhelming probability on the event $\Lambda \leq \sqrt{ \kappa + \eta } / \log(N)^2$ we have
\beq
\Lambda \leq N^{\eps} \sqrt{ \frac{ \Im [ \hatm] }{N \eta } }.
\eeq
On $\hatD_2$ the RHS is $\ll \sqrt{ \kappa + \eta }$ and so we conclude the following.
\bep
For any $\eps >0$ we have with overwhelming probability on $\hatD_1$ that
\beq
\Lambda \leq N^{\eps} \sqrt{ \frac{ \Im [ \hatmt] }{ N \eta } }.
\eeq
\eep

\subsubsection{Fluctuation averaging lemma}
We have,
\bel \label{lem:fa3}
Suppose that $\Lambda \leq \gamma$ with overwhelming probability where
\beq
\frac{1}{N \eta } \leq \gamma \leq \frac{ \sqrt{ \kappa  + \eta }}{\log(N)^2}.
\eeq
Then,
\beq
\ee\left| \frac{1}{N} \sum_i g_i^2 Q_{i} [G_{ii}^{-1} ] \right|^p \leq N^{\eps} \frac{1}{ ( N\eta)^p}.
\eeq
for any $\eps >0$.
\eel
\proof Similarly to before, the moment in question is bounded by
\begin{align}
&\max_{0 \leq s \leq p } \max_{0 \leq l \leq (p+s)/2} \left( \sqrt{ \frac{t}{N} } +t \sqrt{ \frac{ \Im [ \mfct ] + \gamma }{ N \eta } } \right)^{p+s} \frac{1}{ ( \eta + t \sqrt{ \kappa + \eta } )^{s + 2p - 2l} } \frac{1}{N^{p-l} 
} \left( \frac{ \sqrt{ \kappa + \eta } }{ \eta + t \sqrt{ \kappa + \eta } } \right)^l \notag\\
\leq& \max_{0 \leq s \leq p } \frac{ t^{(p+s)/2} }{N^{(p+s)/2} ( \eta + t \sqrt{ \kappa + \eta })^{s+2p} N^{p} } + \max_{0 \leq s \leq p } \frac{ t^{p+s} ( \Im [ \hatmt ] + \gamma)^{(p+s)/2} }{(N\eta)^{(p+s)/2} ( \eta + t \sqrt{ \kappa + \eta })^{s+2p} N^{p} }  \notag\\
+ & \max_{0 \leq s \leq p } \frac{ t^{(p+s)/2} ( \sqrt{ \kappa + \eta } )^{p/2+s/2} }{ N^p ( \eta + t \sqrt{ \kappa + \eta } )^{s/2+3p/2} } \notag\\
+ & \max_{0 \leq s \leq p } \left( \frac{ ( \Im [ \mfct ] + \gamma )^{1/2} t ( \sqrt{ \kappa + \eta } )^{1/2} }{ ( N \eta )^{1/2} N^{1/2} ( \eta + t \sqrt{ \kappa + \eta } )^{3/2} } \right)^p \left( \frac{ t ( \Im [ \mfct ] + \gamma )^{1/2} ( \sqrt{ \kappa + \eta } )^{1/2} }{ \eta^{1/2} ( \eta + t \sqrt{ \kappa + \eta } )^{1/2} } \right)^s \label{eqn:fa4}
\end{align}
where we just bounded the max over $l$ by the sum of the term with $l=0$ and the term with $l=(p+s)/2$ and then subsequently split the $(t/N)^{1/2}$ and $\sqrt{ \Im [ \hatmt ] + \gamma / N \eta}$ in two.  It is immediate that the first term is bounded by $(N \eta)^{-p} (N \eta \sqrt{ \kappa + \eta } )^{-(p+s)/2} \leq (N\eta)^{-p}$.
The contribution from the second term is
\beq
\frac{ ( \Im [ \mfct ] + \gamma )^{(p+s)/2} t^{p+s} }{ ( \eta + t \sqrt{ \kappa +\eta } )^{s+2p} (N\eta)^{(p+s)/2} N^p } \leq \frac{ ( \Im [ \mfct ] + \gamma)^{p/2} }{ (N \eta)^p (N \eta \sqrt{ \kappa+\eta} )^{p/2} (\sqrt{ \kappa + \eta } )^{p/2} ( N \eta \sqrt{ \kappa + \eta } )^{s/2} } \leq \frac{1}{ (  N \eta )^p}
\eeq
using $\gamma \leq \sqrt{ \kappa + \eta }$.  It is immediate that the third term on the second last line of \eqref{eqn:fa4} is bounded by $(N \eta)^{-p}$.  It is easy to see that the term on the last line of \eqref{eqn:fa4} is less than $C (N \eta)^{-p}$ by considering the cases $s=0$ and $s=p$ separately and using $\gamma \leq \sqrt{ \kappa + \eta}$. \qed

\subsubsection{Proof of Proposition \ref{prop:ll3}}
From Proposition \ref{prop:self3} and Lemma \ref{lem:fa3} we immediately see that
\beq
\Lambda \leq \frac{ N^{\eps}}{N \eta }
\eeq
with overwhelming probability. \qed

\subsection{Proof for $E \leq \hatE_- - N^{-2/3+\sigma}$}.

We now define
\beq
\D_2 := \D_\sigma \cap \left\{ E + \i \eta : E \leq \hatE_- - N^{-2/3+\sigma } \right\}.
\eeq
The goal of this section is to prove the following.
\bep \label{prop:ll4}
The estimates of Theorem \ref{thm:ll2} holds on $\hatD_2$.
\eep
The estimates of Lemma \ref{lem:dxi} hold in the set-up of the present section.  Hence,  we see that there is a $C>0$ so that for $z \in \D_2$ that if $\kappa \geq C \eta$, then $|\Re[ \xi ] | \geq c \kappa$ for some $c>0$.  Therefore, we see that for any $\eps >0$ we have
\beq \label{eqn:lla1}
\left| \frac{1}{ N} \sum_i \frac{1}{V_i - \xi } - \hatm (\xi ) \right| \leq N^{\eps} \left( \frac{1}{ N ( \kappa + \eta ) } + \frac{1}{ ( N \eta)^2 \sqrt{ \kappa + \eta } } \right).
\eeq
From Section \ref{sec:shorttimeself} we conclude
\beq
\frac{1}{N} \sum_i |g_i|^p \leq \frac{C}{ ( \kappa + \eta )^{p-3/2}}.
\eeq
Note that
\beq
|V_i - \xi | \geq c ( \kappa + \eta) \geq c \eta.
\eeq
The following is an easy consequence of the fact that  $t^2 \ll \eta$ on $\D_2$.
\bel
We have with overwhelming probability on the event 
\beq
\Lambda \leq N^{-1/3}
\eeq
that
\beq
|Z_i | \leq \frac{ \eta}{ \log(N)^2}.
\eeq
\eel

\subsubsection{Self-consistent equation}
Similar to above we can derive the following.
\bep\label{prop:self4}
With overwhelming probability on the event 
\beq
\Lambda \leq N^{-1/3},
\eeq
we have
\begin{align}
\left| (1 - t R_2 ) ( m_N - \hatmt ) \right| &\leq  \left| \frac{1}{N} \sum_i g_i^2 Q_i [ G_{ii}^{-1} ] \right| + C \frac{ \Lambda^2}{ \sqrt{\eta}} +N^{\eps} \left( \frac{1}{ N ( \kappa + \eta ) } + \frac{1}{ ( N \eta)^2 \sqrt{ \kappa + \eta } } \right) \notag\\
&\leq N^{\eps} \left( \frac{1}{ N^{1/2} ( \kappa + \eta)^{1/4} } + \frac{1}{ N \eta } \right)
\end{align}
for any $\eps >0$ with overwhelming probability.
\eep
For the second inequality we are just using that with overwhelming probability,
\beq
\left| \frac{1}{N} \sum_i g_i^2 Q_i [ G_{ii} ] \right| \leq N^{\eps} \left( \frac{t^{1/2}}{ ( N ( \kappa + \eta ))^{1/2} } + \frac{t}{ \sqrt{ \kappa + \eta } } \sqrt{ \frac{ \Im [ \hatmt ]+\Lambda }{N \eta } } \right) \leq C N^{\eps} \frac{1}{N^{1/2} ( \kappa + \eta)^{1/4}}.
\eeq

\subsubsection{Weak local law}
By Proposition \ref{prop:self4}, we see that with overwhelming probability on the event
\beq
\Lambda \leq N^{-1/3}
\eeq
we have
\beq
\Lambda \leq N^{\eps} \left( \frac{1}{ N^{1/2} ( \kappa + \eta )^{1/4} } + \frac{1}{ N \eta } \right) \leq N^{-\sigma/10} N^{-1/3}.
\eeq
Hence, we conclude the following.
\bep \label{prop:wll4}
With overwhelming probability we have on $\hatD_2$ that
\beq
\Lambda \leq N^{\eps} \left( \frac{1}{ N^{1/2} ( \kappa + \eta )^{1/4} } + \frac{1}{ N \eta } \right).
\eeq
\eep

\subsubsection{Fluctuation averaging lemma}
We now have the following fluctuation averaging lemma.
\bel \label{lem:fa4}
Suppose that $\Lambda \leq \gamma$ with overwheming probability, where $\gamma$ obeys
\beq
\frac{1}{N \eta } \leq \gamma \leq N^{-1/3}.
\eeq
Then we have
\beq
\ee \left| \frac{1}{N} \sum_i g_i^2 Q_i [ G_{ii}^{-1} ] \right|^p \leq N^{\eps} \left( \frac{1}{ N ( \kappa + \eta ) } + \frac{\gamma}{ N \eta \sqrt{ \kappa + \eta } } \right)^p.
\eeq
for any $\eps >0$.
\eel
\proof By \eqref{eqn:momentbd} we have,
\begin{align}
&\ee \left| \frac{1}{N} \sum_i g_i^2 Q_i [ G_{ii}^{-1} ] \right|^p \leq N^{\eps} \max_{0 \leq s \leq p } \max_{ 0 \leq l \leq (p+s)/2} \left( \sqrt{\frac{t}{N} } + t \sqrt{ \frac{\Im [ \hatmt ] + \gamma}{N \eta } } \right)^{p+s} \frac{1}{ ( \kappa+\eta)^{s+2p-3l/2} } \frac{1}{N^{p-l}} \notag\\
\leq & CN^{\eps} \max_{0 \leq s \leq p } \frac{t^{(p+s)/2}}{N^p ( \kappa+\eta)^{s/4+5p/4} } + CN^{\eps} \max_{0 \leq s \leq p }  \frac{ t^{p+s} ( \Im  [ \hatmt] + \gamma )^{(p+s)/2} }{ N^p \eta^{p/2} ( \kappa + \eta)^{s/4+5p/4} \eta^{s/2} }.
\end{align}
In the second inequality we used that the maximum occurs at $l=(p+s)/2$.  Clearly the first term is bounded by $CN^{\eps} (  N ( \kappa + \eta ) )^{-p}$.   For the second, we see that since $t \leq \sqrt{ \eta}$ and $\gamma \leq N^{-1/3}$ that the maximum occurs at $s=0$.  Then,
\beq
\max_{0 \leq s \leq p }  \frac{ t^{p+s} ( \Im  [ \hatmt] + \gamma )^{(p+s)/2} }{ N^p \eta^{p/2} ( \kappa + \eta)^{s/4+5p/4} \eta^{s/2} } \leq C \left( \frac{  (\Im [  \hatmt ] + \gamma)^{1/2} }{N \eta^{1/2} ( \kappa + \eta )^{3/4} } \right)^p \leq C \left( \frac{1}{ N ( \kappa + \eta ) } + \frac{\gamma}{ N \eta \sqrt{ \kappa + \eta } } \right)^p.
\eeq
\qed

\subsubsection{Proof of Proposition \ref{prop:ll4}}.

Starting with $\gamma$ as in Proposition \ref{prop:wll4} we have by Proposition \ref{prop:self4} and Lemma \ref{lem:fa4} that
\beq
\Lambda \leq N^{\eps} \left( \frac{1}{ N( \kappa + \eta ) } + \frac{1}{ ( N \eta )^2 \sqrt{ \kappa + \eta } } + \frac{ \gamma}{ N \eta \sqrt{ \kappa + \eta } } \right).
\eeq
for any $\eps >0$ with overwhelming probability.  Hence by iteration we obtain that $\Lambda \leq N^{\eps} / (N \eta )$ with overwhelming probability.  Taking this choice of $\gamma$ we then get the claim. \qed

\section{Analysis of free convolution law} \label{sec:fcanal}

Let $V$ be $\eta_*-regular$ as in Definition \ref{def:main}.  We will consider
\beq
t := \frac{ N^{\om}}{N^{1/3}}
\eeq
with $\phi_*/2 > 1/3 - \om > 0$.

Recall our definition of $\mfct$ which satisfies
\beq
\mfct (z) = m_V ( \xi)
\eeq
where
\beq
\xi (z) = z + t \mfct (z).
\eeq
Define the map
\beq
F( \xi ) = \xi - t \int \frac{ \d \mu_V (x) }{ x - \xi },
\eeq
so that
\beq
z = F ( \xi ).
\eeq

An important role is played by the contour $\gamma$ which we define as the image of $\rr$ under the map $E \to \xi (E)$.  

A useful observation is that since
\beq
\Im [ \mfct ] = \left( 1 - t \int \frac{ \d \mu_V (x) }{ |x-\xi |^2} \right)^{-1} \Im [ m_V ( \xi ) ]
\eeq
we have,
\beq \label{eqn:absr2}
\int \frac{ \d \mu_V (x)}{ |x- \xi |^2 } \leq \frac{1}{t}
\eeq
from which the inequality
\beq
|\mfct | \leq t^{-1/2}
\eeq
follows by Cauchy-Schwartz.   Hence, we see that $\Im [ \xi (E) ] = 0$ for $E \in [-3/4, -1/2]$.

By Lemma \ref{lem:imm} it is easy to check that there is a unique solution $\xi_- \in [-3/4, 3/4]$ to the equation
\beq
1 = t \int \frac{ \d \mu_V}{ (x - \xi_- )^2}.
\eeq
Moreover, we see that
\beq \label{eqn:ximin}
- \xi_- \asymp t^2.
\eeq
Let $E_-$ be such that $\xi (E_- ) = \xi_-$.  For $E \geq E_-$, $\xi (E)$ has non-trivial imaginary part, equalling $t \rhofct (E)$.  

We will write
\beq
\xi (E) := a + b \i
\eeq
for $a, b \in \rr$.  In general, $a$ is a strictly increasing function of $E$, and $a$ and $b$ solve
\beq \label{eqn:abdef}
1 = t \int \frac{ \d \mu_V }{ (x-a)^2 + b^2}.
\eeq
We denote $a_- = \Re[ \xi_- ] = \xi_-$.  Our first goal is to get qualitative behaviour of the contour $\gamma$.   We have,
\bel \label{lem:ab}
For $3/4 \geq a \geq a_-$, we have
\beq
b \asymp t |a- a_- |^{1/2}.
\eeq
\eel
\proof We first consider $a$ near $a_-$.  First, by \eqref{eqn:ximin} and Lemma \ref{lem:imm} we see that for a small $c>0$ we have
\beq
|F^{(k)} ( \xi ) | \leq  \frac{C}{t^{2k-2}}, \qquad |\xi - \xi_- | \leq c t^2
\eeq
for any $k \geq 2$.  Moreover, we see that
\beq
|F^{(k)} ( \xi_- ) |\asymp \frac{C}{t^{2k-2}}
\eeq
for $k \geq 2$.  Note that they are real numbers.   Hence, for $|\xi - \xi_- | \leq c t^2$, we can expand
\beq
F( \xi ) - F (\xi_- ) = \frac{F'' ( \xi_- )}{2} ( \xi - \xi_- )^2 + \frac{ F''' ( \xi_- )}{ 6} (\xi - \xi_- )^3 + \O \left( t^{-6} | \xi - \xi_- |^4\right)
\eeq
The LHS equals $z-E_-$.  We can set $z = E$ and invert this, obtaining first that
\beq \label{eqn:xiexpand1}
\xi - \xi_- = \sqrt{ \frac{ 2 ( E - E_- ) }{ F'' ( \xi_- ) } } \left( 1 - \frac{F''' ( \xi_- ) }{ 3 F'' ( \xi_- ) } (\xi - \xi_- ) + \O \left( t^{-4} |\xi - \xi_- |^2 \right) \right),
\eeq
and then back-substituting once more we obtain
\beq
\xi - \xi_- = \sqrt{ \frac{  2 (E- E_- ) }{ F'' ( \xi_- ) } } \left( 1  - \frac{ F''' ( \xi_- )}{ 3 F'' ( \xi_- ) } \sqrt{ \frac{  2 (E- E_- ) }{ F'' ( \xi_- ) } }   + \O \left( t^{-4} |\xi - \xi_- |^2 \right) \right).
\eeq
Taking real and imaginary parts (note that $F'' (\xi_-)$ is negative) we see that
\beq \label{eqn:ader2}
| a - a_- | \asymp |E - E_- | , \qquad b \asymp t |E- E_- |^{1/2} \asymp t |a- a_- |^{1/2},
\eeq
for $| \xi - \xi_- | \leq c t^2$.

Now, a straightforward calculation using Lemma \ref{lem:imm} and \eqref{eqn:abdef} shows that there is a small $ c>0$ so that $b \asymp t |a-a_- |$ for $- c t^2 \leq a \leq 3/4$.

It remains to consider the region $a_- + c_1 t^2 \leq a \leq - c_1 t^2$ for a small $c_1>0$.  That is, we need to prove that $b \asymp t^2$ here.  Note that the upper bound $b \leq C t^2 \leq C t |a-a_- |^{1/2}$ follows immediately from \eqref{eqn:imm2}.  For the lower bound we compute
\begin{align}
0 &= \frac{1}{t} - \int \frac{ \d \mu_V (x) }{ (x-a )^2 + b^2 } \notag \\
&= \int \frac{ \d \mu_V (x) }{ (x-a_- )^2 } - \int \frac{ \d \mu_V (x) }{ (x-a )^2 + b^2 } \notag \\
&= \int \frac{ ( x-a)^2 + b^2 - (x-a_- )^2 } { (x-a_-)^2 ( (x-a)^2 +b^2 ) } \d \mu_V (x) \notag \\
&= \int \frac{ (a-a_- ) (-2x + a + a_- ) + b^2 }{ (x-a_- )^2 ( (x-a)^2 + b^2 ) } \d \mu_V (x).
\end{align}
We rearrange this to get
\beq \label{eqn:ab3}
\frac{b^2}{a-a_-} = \frac{ \int \frac{ 2 x - a - a_- } { (x -a_-)^2((x-a)^2 + b^2 ) } \d \mu_V (x) }{ \int \frac{ 1 } { (x -a_-)^2((x-a)^2 + b^2 ) } \d \mu_V (x)}.
\eeq
Note that for $x \in \supp  ( \mu_V )$ we have $|x-a| \geq c t^2$ and since $b \leq c t^2$ we get
\beq \label{eqn:ab4}
(x-a)^2 + b^2 \asymp (x-a)^2 \asymp (x-a_-)^2.
\eeq
We need to get a lower bound on the numerator of \eqref{eqn:ab3} and an upper bound on the denominator.  For the numerator, we have
\begin{align}
\int \frac{ 2 x - a - a_- } { (x -a_-)^2((x-a)^2 + b^2 ) } \d \mu_V (x)  &= \int_{x \geq -1/2}  \frac{ 2 x - a - a_- } { (x -a_-)^2((x-a)^2 + b^2 ) } \d \mu_V (x)\notag \\
&+ \int_{x \leq -1/2}  \frac{ 2 x - a - a_- } { (x -a_-)^2((x-a)^2 + b^2 ) } \d \mu_V (x) \notag \\
&\geq  c t^2 \int_{x \geq -1/2}  \frac{ 1} { (x -a_-)^2((x-a)^2 + b^2 ) } \d \mu_V (x) - C \notag \\
&\geq c t^2 \int_{\rr}   \frac{ 1} { (x -a_-)^2((x-a)^2 + b^2 ) } \d \mu_V (x) - 2 C \notag \\
& \geq c t^2 \int_\rr \frac{1}{ (x-a_-)^4} \d \mu_V (x) - 2 C \notag\\
&\geq c t^2 \int_\rr \frac{1}{ (x-a_-)^4} \d \mu_V (x)
\end{align}
In the first inequality we used that the integrand is bounded in the region $x \leq -1/2$ and that for $x \geq -1/2$ that $x-a \geq x - a_- \geq c t^2$ for $x \in \supp ( \mu_V)$.  In the third inequality we again used the fact that the integrand is bounded for $x \leq -1/2$.  in the fourth inequality we use \eqref{eqn:ab4} and in the final inequality we use that the integral is order $t^{-5}$.

For the integral in the denominator of \eqref{eqn:ab3} we use \eqref{eqn:ab4} to prove that
\beq
 \int \frac{ 1 } { (x -a_-)^2((x-a)^2 + b^2 ) } \d \mu_V (x) \leq C \int \frac{1}{ (x-a_- )^4} \d \mu_V (x),
\eeq
which yields $b^2 \geq c t^2 (a-a_- )$ and completes the proof. \qed

In order to complete our calculation of the contour $\gamma$ we need the following.
\bel
We have $|a - a_-| \asymp |E - E_- |$.
\eel
\proof 
We already know from \eqref{eqn:ader2} that this holds for $|a - a_- | \leq c t^2$ for a small $c>0$.  The claim will therefore be proved by showing that
\beq
\frac{ \d a }{ \d E} \asymp 1
\eeq
for $|a-a_- | \geq c t^2$.

We calculate
\beq \label{eqn:ader}
\frac{ \d a }{ \d E} = \frac{ \Re \left[ 1 - t \int \frac{1}{ (x-\xi )^2 } \d \mu_V (x) \right]}{ \left| 1 - t \int \frac{1}{ (x- \xi )^2 } \d \mu_V (x) \right|^2 }.
\eeq
Since $1 = t \int |x-\xi|^{-2} \d \mu_V (x)$ the numerator of \eqref{eqn:ader} is positive and the denominator is less than $2$, so
\beq
\frac{ \d a }{ \d E} \geq \frac{1}{2} \Re \left[ 1 - t \int \frac{1}{ (x- \xi )^2 } \d \mu_V (x) \right].
\eeq 
Clearly,
\beq
\frac{ \d a }{ \d E } \leq \left( \Re \left[ 1 - t \int \frac{1}{ (x- \xi )^2 } \d \mu_V (x) \right]\right)^{-1} .
\eeq
Hence, the claim will follow by proving
\beq
\Re \left[ 1 - t \int \frac{1}{ (x- \xi )^2 } \d \mu_V (x) \right] \geq c >0
\eeq
 We can write
\begin{align}
\Re \left[ 1 - t \int \frac{1}{ (x- \xi )^2 } \d \mu_V (x) \right] &= \left( 1 - t \int \frac{1}{ |x - \xi |^2 } \d \mu_V (x) \right) + 2 b^2 t \int \frac{1}{ |x - \xi|^4} \d \mu_V (x) \notag \\
&= 2 b^2 t \int \frac{1}{ |x - \xi|^4} \d \mu_V (x) .
\end{align}
By Lemma \ref{lem:imm} we have for $a \geq 0$,
\beq
b^2 t \int \frac{1}{ |x - \xi|^4} \d \mu_V (x)  \asymp \frac{ t b^2(a + b)^{1/2} }{ b^3} \asymp \frac{ (a + b)^{1/2}}{ |a- a_- |^{1/2} } \asymp 1
\eeq
where we used Lemma \ref{lem:ab} in the second last step.  For $a_- + c t^2 \leq a \leq 0$ we have
\beq
b^2 t \int \frac{1}{ |x - \xi|^4 } \d \mu_V (x) \asymp  \frac{ t b^2}{ (|a| + b)^{5/2}} \asymp \frac{ t^3 |a - a_- | }{ (|a - a_-|^{1/2} t + |a| )^{5/2} } \asymp 1.
\eeq
The claim follows. \qed

\subsection{Estimates on the map $\xi$ and control of self-consistent equation coefficients}

The following result will be useful.
\bel \label{lem:r2}
The following holds for $|E| \leq 1/2$ and $0 \leq \eta \leq 10$.  We have,
\beq
\left|1 - t \int \frac{1}{ (x- \xi)^2 } \d \mu_V (x) \right| \asymp \min \left\{ \frac{ \sqrt{ |E-E_- | + \eta }}{ t}, 1 \right\}
\eeq
\eel
\proof We denote $\xi = a + b \i$.  Unlike in the previous subsection we do not restrict $\xi$ to lie on the contour $\gamma$.  First assume that $|a-a_-| + b \geq c t^2$.    We can write
\begin{align}
\Re \left[ 1 - t \int \frac{1}{ (x- \xi )^2 } \d \mu_V \right] = \left( 1 - t \int \frac{1}{ |x - \xi |^2 } \d \mu_V (x) \right) + 2 b^2 t \int \frac{1}{ |x- \xi|^4 } \d \mu_V (x).
\end{align}
The term in the brackets on the RHS is always positive.  For $a \geq 0$ we have by Lemma \ref{lem:imm},
\beq
t \int \frac{1}{ | x - \xi |^2 } \asymp \frac{ t ( a^{1/2} + b^{1/2})}{b} =: Q
\eeq
and
\beq
t b^2 \int \frac{1}{ |x- \xi |^2} \d \mu_V (x) \asymp Q.
\eeq
Hence,
\beq
\Re \left[ 1 - t \int \frac{1}{ (x- \xi )^2 } \d \mu_V \right]  \geq \max \left\{ 1 - C Q, c Q \right\} \geq c.
\eeq
Now we assume $a \leq 0$.  We have,
\beq
t \int \frac{1}{ |x - \xi|^2} \d \mu_V (x) \asymp \frac{t}{ (|a| + b)^{1/2} },
\eeq
and so there is a constant $C$ s.t. if $|a|^{1/2} + b^{1/2} \geq C t$ then 
\beq
\left( 1- t \int \frac{1}{ |x - \xi|^2} \d \mu_V (x) \right) \geq c.
\eeq
Now we assume that $c_1 t^2 \leq |a - a_- | + |b| \leq C_1 t^2$, for given $c_1$ and $C_1$, as well as $a \leq 0$.  Assume first that $ a \geq a_-$.     In this regime we must have $b \geq c t |a-a_- |^{1/2}$, as $b$ must lie above the contour $\gamma$ in $\cc$.  Hence the assumption $c_1 t^2 \leq |a - a_- | + |b|$ implies $b \geq c t^2$ for another $c>0$ depending on $c_1$.  Hence,
\beq \label{eqn:r22}
b^2 t \int \frac{1}{ |x - \xi |^4 } \d \mu_V (x) \asymp \frac{ b^2 t}{ ( |a| + b )^{5/2} } \geq \frac{b^2 t}{ ( |a_-| + |a-a_- | + b )^{5/2} }  \geq c \frac{ t^5}{ (t^2 +t^2 )^{5/2} } \geq c.
\eeq
We postpone the case $c_1 t^2 \leq |a- a_- | + |b| \leq C_1 t^2$ and $a \leq a_-$.  First, we consider the case $|a-a_-| + b \leq c t^2$ for a small $c>0$.   Then we can expand
\beq
1 -t \int \frac{1}{(x- \xi )^2 } \d \mu_V (x) = F'' ( \xi_- ) ( \xi - \xi_-) \left( 1 + \O ( t^{-2} | \xi - \xi_- | \right).
\eeq
Hence for a small enough $c$, 
\beq
\left| 1 - \int \frac{1}{(x- \xi )^2 } \d \mu_V (x)  \right| \asymp t^{-2} ( |a - a_- | + b ) \asymp t^{-1} ( |E-E_- | + \eta )^{1/2}, \qquad |\xi  - \xi_- | \leq c t^2.
\eeq
We now return to the postponed case above.  Since $|a-a_- | + b \geq c_1 t^2$ we have that $b \geq c_1/2 t^2$ if $|a-a_- | \leq c_1/2 t^2$.  Then calculation \eqref{eqn:r22} applies, and so we can instead work in the regime $c_1 t^2 \leq |a-a_-| + b  \leq C_1 t^2$ and $a \leq a_- - c_1 t^2$.  Then,
\begin{align}
1 -  t\int \frac{ \d \mu_V } { (x-a)^2 + b^2 } &\geq 1 - \int \frac{ \d \mu_V }{ (x-a)^2 }  \notag\\
&\geq 1-  t \int \frac{ \d \mu_V }{ (x - (a_- - c_2 t^2 ) )^2} - C t
\end{align}
for any $c_2 >0$ such that $c_2 < c_1$.   The $Ct$ appearing above is a bound for the contribution of the integral from $x \in \supp ( \mu_V )$ such that $x \leq -1$.  If we take $c_2$ small enough so that the estimate proved in the expansion above holds, then 
\beq
1-  t \int \frac{ \d \mu_V }{ (x - (a_- - c_2 t^2 ) )^2}  \asymp 1.
\eeq
This completes the proof. \qed

Again we write $\xi = a + b \i$ for general $E + \i \eta$ with $|E| \leq 1/2$ and $0 \leq \eta \leq 10$.
First we remark that the proof of Lemma \ref{lem:r2} and \eqref{eqn:ader} immediately yield that
\beq
\frac{ \d a }{ \d E} \asymp 1, \qquad \kappa + \eta \geq c t^2
\eeq
for any $c>0$.  By the Cauchy-Riemann equations we also reach the same conclusion for $\d b / \d \eta$.  Since
\beq
\xi' (z) = \frac{1}{1 - t \int \frac{\d \mu_V (x) }{ (x- \xi )^2}}
\eeq
we then see that
\beq
| \frac{ \d a }{ \d \eta } | + | \frac{ \d b }{ \d E } | \leq C
\eeq
for $\kappa + \eta \geq c t^2$.

Hence, we have proved the following.
\bel \label{lem:dxi}
We have in the region $\kappa + \eta \geq c t^2$, for any $c>0$, that
\beq
\frac{ \d a }{ \d E } \asymp 1, \qquad \frac{ \d b }{ \d \eta } \asymp 1.
\eeq
In the same region we have
\beq
\left| \frac{ \d a }{ \d \eta } \right| + \left| \frac{ \d b } { \d E } \right| \leq C.
\eeq
The above implies
\beq
|a| + |a-a_- | + b \leq C (t^2 + \eta + \kappa ).
\eeq
We also have,
\beq \label{eqn:distsq}
t \sqrt{ \kappa + \eta } \asymp | \xi - \xi_- |, 
\eeq
for $\kappa + \eta \leq c t^2$ for some small $c>0$.
\eel
The estimate \eqref{eqn:distsq} follows from \eqref{eqn:xiexpand1}.  
The following lemma is needed for the proof of the local law.
\bel \label{lem:llcoef1}
Consider $z \in \D_{\sigma}$.  
In the region $E \geq E_- - t^2$, we have
\beq
|V_i - \xi | \geq c (t^2 + \eta  + t \Im  [ \mfct ] ),
\eeq
and
\beq \label{eqn:coef2}
 \int \frac{ \d \mu_V(x) }{ |x - \xi |^p } \leq C \frac{ t + \sqrt{ \kappa + \eta } }{ (t^2 + \Im [ \xi ])^{p-1} }.
\eeq
On $\D_1$ we have,
\beq \label{eqn:coef1}
t + \sqrt{ \kappa + \eta } \leq C ( t + \Im[ \mfct ]).
\eeq
\eel
\proof First, note that $|V_i - \xi | \geq \Im [ \xi] = \eta + t \Im [ \mfct ]$.  If $\eta  + t \Im [ \mfct ] \leq  c_1 t^2$ then $\kappa + \eta \leq C c_1 t^2$ for another $C>0$.  Choosing  $c_1$ sufficiently small we see that $| \xi - \xi_- (E_- ) | \leq c t^2$ for any small $c$, from Lemma \ref{lem:dxi}.  Since $\xi_- (E_- ) \leq - c t^2$ we see that $| V_i - \xi | \geq c t^2$.

For the next estimate, first suppose that $\Im [ \xi ] \leq c t^2$ for a small $c>0$.  Choosing $c$ small enough we see that, as above, $a \leq - ct^2$ and so
\beq
 \int \frac{ \d \mu_V(x) }{ |x - \xi |^p }  \leq C \frac{1}{ ( |a| + b )^{p-3/2} } \leq C \frac{1}{ ( t^2 + \Im [ \xi ] )^{p-3/2} } \leq C \frac{t}{ (t^2 + \Im [ \xi ] )^{p-1}}.
\eeq
We may assume that $\Im [ \xi ] \geq c t^2$ for a small $c>0$.  Then, noting that we always have \eqref{eqn:imm1} as an upper bound, no matter the sign of $a$, we get
\beq
 \int \frac{ \d \mu_V(x) }{ |x - \xi |^p }  \leq C \frac{ \sqrt{ |a| + b } }{ ( \Im [ \xi ] )^{p-1} } \leq C \frac{ \sqrt{ |a| + b } }{ ( t^2 +  \Im [ \xi ] )^{p-1} }.
\eeq
Using Lemma \ref{lem:dxi} to bound the numerator, we complete the proof of \eqref{eqn:coef2}.
 The estimate \eqref{eqn:coef1} is easy. \qed

We also have the following.
\bel \label{lem:llcoef2}
For $z \in \D_\sigma$ and $E \leq E_- - t^2$ we have,
\beq \label{eqn:coef3}
|V_i - \xi | \geq c (t^2 + \kappa + \eta )
\eeq
and
\beq
\int \frac{ \d \mu_V (x) }{ |x - \xi |^p } \leq \frac{ C}{ ( \kappa + \eta + t^2 )^{p-3/2} }.
\eeq
\eel
\proof For the first estimate, as argued in the proof of Lemma \ref{lem:llcoef1}, we see that if $\kappa + \eta \leq c t^2$ for a small $c>0$, then we already arrive at $|V_i - \xi | \geq c t^2$.  Since $\Im[ \xi] \geq \eta$ we then get that $|V_i - \xi | \geq c (t^2 + \eta)$.  Then by Lemma \ref{lem:dxi}, we note that there is a $C>0$ so that if $E \leq E_- - C \eta - t^2$, then $a \leq - c_1 \kappa$ for anther $c_1 >0$.  Hence, $|V_i - \xi | \geq C \kappa$ for such $E$, which completes \eqref{eqn:coef3}.

Building on this observation, we see that there is a $C>0$ so that if $E \leq E_- - C \eta - t^2$, then
\beq
\int \frac{ \d \mu_V (x) }{ |x - \xi |^p } \leq C \frac{1}{ \kappa^{p-3/2}} \leq C \frac{1}{ ( \kappa + \eta + t^2 )^{p-3/2}}
\eeq
where we applied \eqref{eqn:imm2}.  So we can assume that $\kappa \leq C \eta + t^2$.  If $\eta \geq t^2$, then the desired estimate immediately follows from \eqref{eqn:imm1} and $|a|+|b| \leq C (t^2 +\eta + \kappa)$.  In the case $\eta \leq t^2$ then the result follows from \eqref{eqn:imm1} and just counting powers of $t$.\qed

We also have the following.
\bel  There is a $c_1>0$ so that if $\kappa + \eta \leq c_1 t^2$, then
\beq
\left| \int \frac{ \d \mu_V (x) }{ (x-\xi )^3} \right| \asymp \frac{1}{t^3}
\eeq
\eel
\proof First, note that the claim follows for 
\beq
\int_{x \geq -1/2} \frac{ \d \mu_V (x) } {|x- \xi |^3}
\eeq
by Lemma \ref{lem:imm}.  Then note that for $x \geq -1/2$ and $x \in \supp ( \mu_V )$, that
\beq
\left| \Re[ \frac{1}{ (x - \xi )^3} ] \right| \asymp \frac{1}{ |x-\xi|^3},
\eeq
by direct calculation, if we take $c_1>0$ sufficiently small.  This yields the claim. \qed

\subsection{Qualitative properties of $\rhofct$ and $\mfct$}

We first prove the following.
\bel  \label{lem:rhofctexpand} We have the following for $|E| \leq 3/4$.
The density $\rhofct$ satisfies
\beq \label{eqn:rhosq}
\rhofct (E) \asymp |E - E_- |^{1/2} \1_{ \{ E \geq E_- \} }.
\eeq
Moreover, for $|E - E_- | \leq c t^2$ we have,
\beq \label{eqn:rhoexp}
\rhofct (E) = \sqrt{ \frac{ ( E - E_-) }{ t^2 |F'' ( \xi_- )|}}\left(1 + \O \left(  \frac{ |E - E_- |}{t^2} \right) \right)
\eeq.
We have also
\beq
\left|  t^2 F'' ( \xi_- ) \right|  \asymp 1.
\eeq
\eel
\proof We have already proved \eqref{eqn:rhosq}, because $\Im [ \xi] = t \rhofct (E)$. Equation \eqref{eqn:rhoexp} follows from continued back-substitution in \eqref{eqn:xiexpand1} (note that the correction to the $|E-E_-|^{1/2}$ term above is $|E-E_- |^{3/2}$ instead of $|E-E_-|$ - this is due to all the coefficients in \eqref{eqn:xiexpand1} being real).   The final estimate is a consequence of the fact that $- \xi_- \asymp t^2$ and Lemma \ref{lem:imm}. \qed

Now since $\mfct$ has a square root behaviour, we get the following.
\bel \label{lem:immfct1}
We have for $E \geq E_-$, 
\beq
\Im [ \mfct ] \asymp \sqrt{ \kappa + \eta}
\eeq
and for $E \leq E_-$,
\beq
\Im [ \mfct ] \asymp \frac{\eta}{ \sqrt{ \kappa + \eta}}.
\eeq
\eel
We have the equality
\beq
\del_z \mfct (z) = \left( 1 - t \int \frac{ \d \mu_V (x) }{ (x - \xi )^2 } \right)^{-1} \int \frac{ \d \mu_V (x)}{ (x- \xi )^2},
\eeq
from which, using Lemma \ref{lem:r2} and \eqref{eqn:absr2} we conclude that
\beq
| \del_z \mfct | \leq C \max\left\{ \frac{1}{ \sqrt{ \kappa  + \eta } } , \frac{1}{t} \right\}.
\eeq
Combining this with the trivial estimate
\beq
|\del_z \mfct | \leq \frac{ \Im [ \mfct ] }{ \eta }
\eeq
we obtain
\bel
For $\kappa + \eta \leq  t^2$ we have
\beq
| \del_z \mfct | \leq \frac{C}{ \sqrt{ \kappa + \eta }}.
\eeq
For $\kappa + \eta \geq t^2$ we have for $E \geq E_-$,
\beq
| \del_z \mfct | \leq C \frac{ \sqrt{ \kappa + \eta} }{ t\sqrt{ \kappa + \eta } + \eta },
\eeq
and for $E \leq E_-$,
\beq
| \del_z \mfct | \leq \frac{C}{ \sqrt{ \kappa + \eta }}.
\eeq
\eel

\subsection{Comparison of free convolutions of matching measures} \label{sec:matching}
In this section the set-up is the following.  We have two measures $\mu_1$ and $\mu_2$ that have densities on the interval $[-1, 1]$ such that,
\beq \label{eqn:rhoicl}
\rho_1 (x)= \rho_2 (x ) \left( 1 + \O\left( \frac{|x|}{t_0^2} \right) \right), \qquad 0 \leq x \leq c t_0^2,
\eeq
and $\rho_1 (x) = \rho_2 (x) = 0$ for $-1 \leq x \leq 0$ and $\rho_2 (x) \asymp \sqrt{x}$.

We assume that for $|x| + \eta \leq c t_0^2$ that we have
\beq \label{eqn:mider}
| \del_z m_i (z) | \leq \frac{C}{\sqrt{ |x| + \eta }}.
\eeq

We consider the free convolutions $m_{1, t}$ and $m_{2, t}$ for $0 \leq t \leq t_0 N^{-\eps_0}$, for some $\eps_0 >0$.  Denote the maps $\xi_i = z + t m_{i, t} (z)$, as well as the points $\xi_{i, -}$ as above.  Our goal is to compare the densities $\rho_{i, t}$ to each other.

Since the $\rho_i (x)$ are continuous densities behaving like a square root, the qualitative behaviour of Definition \ref{def:main} holds down to $\eta_*=0$.  Hence, the analysis of the previous subsection goes through, and the contours $\gamma_i = \xi_i ( \rr )$ have the same qualitative behaviour, i.e., there is a $E_{-, i}$ at which they leave the real line, and for $E \geq E_{-, i}$ we have
\beq
|a_i-a_{-, i} | \asymp |E-E_{-, i} |, \qquad b_i \asymp t |E-E_{-, i} |^{1/2},
\eeq
etc.  Moreover, it is not hard to check that the estimate
\beq \label{eqn:mitder}
| \del_z m_{i, t} (z) | \leq \frac{C}{ \sqrt{ |E-E_{-, i} | + \eta }}
\eeq
using the methods of the previous section, as well as \eqref{eqn:mider}.

We first check that the $a_{-, i}$ are close.   Using the equation that defines them we get
\begin{align}
0 &= \int \frac{ \d \muo (x) }{ ( x- \aom )^2 } - \int \frac{ \d \mut (x) }{ (x - \atm )^2 } \notag \\
&= (\aom - \atm ) \int \frac{  ( 2 x - \aom - \atm )}{ (x - \aom)^2 ( x - \atm )^2 } \d \muo (x)  + \left( \int \frac{ \d \muo (x) }{ ( x - \atm )^2 } - \int \frac{ \d \mut (x)}{ (x- \atm )^2} \right) 
\end{align}
It is easy to check using the assumption \eqref{eqn:rhoicl} that
\beq
\left| \int \frac{ \d \muo (x) }{ ( x - \atm )^2 } - \int \frac{ \d \mut (x)}{ (x- \atm )^2} \right|
 \leq \frac{C}{t_0}
\eeq
We bound below the factor multiplying $(\aom - \atm )$ by
\begin{align}
\left| \int \frac{  ( 2 x - \aom - \atm )}{ (x - \aom)^2 ( x - \atm )^2 } \d \mu_1 (x) \right|  = \left| \int \frac{ \d \muo (x) }{ ( x - \aom )^2 ( x - \atm ) } + \int \frac{ \d \muo (x) }{ ( x - \aom ) ( x - \atm )^2 }  \right|
\end{align}
Each term on the RHS is positive and is order $t^{-3}$.  Hence,
\beq
| \aom - \atm | \leq C t^2 \frac{t}{t_0}.
\eeq
It this then easy to see that,
\beq
\left| \Fo^{(k)} ( \xiom ) - \Ft^{(k)} ( \xitm ) \right| \leq C \frac{t}{t_0} \frac{1}{ t^{2(k-1)}}
\eeq
Consider now the expansion 
\beq
E - E_{-, i} = \sum_{k=2}^m \frac{ F_i^{(l)} ( \xi_{-, i} ) }{j!} ( \xi_i - \xi_{-, i} )^j + \O \left( |\xi_i - \xi_{-, i} |^{m+1} t^{-2m} \right).
\eeq
By repeated back-substitution, as in the proof of Lemma \ref{lem:ab} we see that for any $\delta_0 >0$ we have, for $0 \leq x \leq N^{-\delta_0} t^2$, that
\begin{align}
\ao(x + \Eom ) - \aom &= (\at (x + \Etm ) - \atm ) (1 + \O (t / t_0 +N^{-D}) ),\notag\\
   \bo (x + \Eom ) &= \bt (x + \Etm )(1  + \O (t/t_0  +N^{-D})).
\end{align}
for any $D>0$, if we take $m$ large enough (but finite), depending on $\delta_0$.

We make the choice
\beq
\delta_0 < \eps_0/100.
\eeq

We now need to deal with 
\beq
 t^2 N^{- \delta_0} \leq x \leq t_0^2.
 \eeq 
   Since we know that $\d a / \d E$ is increasing we can parameterize $b = b(a)$.  We then want to determine the difference $\bo (a) - \bt (a)$ beyond their natural scale.  We have the equation
\begin{align}
0 &= \int \frac{ \d \muo (x) }{ (x - a)^2 + \bo^2 } - \int \frac{ \d \mut (x) }{ (x- a)^2 + \bt^2 } \notag \\
&= (b_2 - b_1)(b_1 + b_2 ) \int \frac{ \d \muo (x) }{ ( (x-a )^2 + \bo^2 )(  (x-a)^2 + \bt^2 )} + \int \frac{ \d \muo (x) - \d \mut (x) }{ (x-a)^2 + \bt^2 }.
\end{align}
By \eqref{eqn:rhoicl} we have the estimate,
\beq
\left| \int \frac{ \d \muo (x) - \d \mut (x) }{ (x-a)^2 + \bt^2 } \right| \leq 
 C \frac{1}{ t_0}.
\eeq
To lower bound the integral multiplying the $b$'s, we note that since $|a - \aom | \geq c t^2 N^{-\delta_0}$ and $|a - \atm | \geq c t^2 N^{- \delta_0 }$, and that  by our choice of $\delta_0$ that $c t^2 N^{-\delta_0} \gg |\aom - \atm |$, we have  $|a - \aom | \asymp | a - \atm |$.  This also implies $b_1 \asymp b_2$.  
We have the lower bound 
\beq
\left| \int \frac{ \d \muo (x) }{ ( (x-a )^2 + \bo^2 )(  (x-a)^2 + \bt^2 )}  \right| \geq  \frac{c \1_{ \{ a \leq 0 \} }}{ ( |a| +|a - \aom|^{1/2}  t )^{5/2} } + \frac{ c \1_{\{ a \geq 0 \} } ( |a| + |a - \aom |^{1/2} t )^{1/2} }{ ( t |a-\aom |^{1/2} )^3 }
\eeq
Hence, using $|a| \1_{ \{ a \leq 0 \} } \leq C t^2$, and that $|a - \aom | \asymp |a| + t |a-\aom |^{1/2}$ when $a \geq 0$,
\begin{align}
| \bo- \bt | &\leq \frac{  C \1_{ \{ a \leq 0 \} }  ( |a| + | a - \aom |^{1/2} t )^{5/2} }{ t_0 |a - \aom |^{1/2} t } + \frac{ C \1_{ \{ a \geq 0 \} } ( t |a - \aom |^{1/2} )^2 }{ t_0 ( |a - \aom |^{1/2} t + |a| )^{1/2} }  \notag \\
& \leq \frac{ C \1_{ \{ a \leq 0 \} } t^5}{ t_0 N^{-\delta} t^2} + \frac{ C \1_{ \{ a \leq 0 \} }( | a- \aom |^{1/2} t ) t}{ t_0 }  + \frac{ C \1_{ \{ a \geq 0 \} } t^2 | a - \aom |^{1/2}}{ t_0 } \notag \\
& \leq C \frac{ N^{2 \delta} t }{ t_0 } \left( t |a- \aom |^{1/2} \right). \label{eqn:bbds}
\end{align}

We now want to use this to show that $\d \ao / \d E$ and $\d \at / \d E$ are close for $a$ at least $c t^2 N^{- \delta_1}$ away from $\aom$ and $\atm$.  This means that we need to study the function
\beq
1 - t \int \frac{1}{ (x - \xio)^2} \d \muo (x).
\eeq
Fix $a$.  We write
\begin{align}
\left( 1 - t \int \frac{ \d \muo (x) }{ ( x - \xio )^2 } \right) - \left( 1 - t \int \frac{ \d \mut (x) }{ (x - \xit )^2 } \right) &= \left( t \int \frac{1}{ (x - \xit )^2} - \frac{1}{ (x - \xio )^2 }  \d \mut (x) \right) \notag \\
&+ \left( t \int \frac{1}{ ( x - \xio )^2} \d \mut (x) - \d \muo (x) \right).
\end{align}
As above, we have
\beq
\left| t \int \frac{1}{ ( x - \xio )^2} \d \mut (x) - \d \muo (x) \right| \leq C \frac{ t}{t_0} .
\eeq
For the other term we write it as
\beq
\left( t \int \frac{1}{ (x - \xit )^2} - \frac{1}{ (x - \xio )^2 }  \d \mut (x) \right) =t (\bo - \bt ) \int \frac{ (x - \xio ) + (x - \xit )}{ (x - \xio )^2 ( x - \xit )^2 }  \d \mut (x).
\eeq
For $a \leq 0$ we have, using our bounds on $\bo - \bt$ and Lemma \ref{lem:imm},
\beq
\left|t  (\bo - \bt ) \int \frac{ (x - \xio ) + (x - \xit )}{ (x - \xio )^2 ( x - \xit )^2 }  \d \mut (x) \right| \leq C \frac{ N^{2 \delta} t}{t_0}  \frac{ t^2 |a - \aom |^{1/2}}{ (|a| + t |a- \aom |^{1/2})^{3/2} } \leq  C \frac{ N^{2 \delta} t}{t_0} 
\eeq
where we used $|a| + t |a- \aom |^{1/2} \geq c t^2$.  For $a \geq 0$ we have
\beq
\left|t  (\bo - \bt ) \int \frac{ (x - \xio ) + (x - \xit )}{ (x - \xio )^2 ( x - \xit )^2 }  \d \mut (x) \right|  \leq C \frac{ N^{2 \delta} t}{t_0} \frac{ t |a - \aom |^{1/2} t ( |a|^{1/2} + t^{1/2} |a - \aom |^{1/4} )} { t^2 |a| } \leq C \frac{ N^{2 \delta} t}{t_0} .
\eeq
From this we get that for $|a - \aom | \geq c t^2 N^{- \delta}$,
\beq \label{eqn:Eder}
\left| \frac{ \d \Eo }{ \d a} - \frac{ \d \Et }{ \d a} \right| \leq C \frac{N^{5 \delta} t }{ t_0}.
\eeq
From our earlier expansion we see that for $x \leq c t^2 N^{-\delta}$,
\beq
(\Eo (x + \aom ) - \Eom ) = ( \Et (x + \atm ) - \Etm ) \left(1 + \O \left( t/t_0 \right) \right).
\eeq
We then use \eqref{eqn:Eder} and get 
\beq \label{eqn:Eicompare}
(\Eo (x + \aom ) - \Eom ) = ( \Et (x + \atm ) - \Etm ) \left(1 + \O \left( N^{5 \delta} t/t_0 \right) \right),
\eeq
for $0 \leq x \leq t_0^2 N^{-10 \delta}$.  This easily implies an estimate on $\ao $ and $\at$.  Indeed, define the functions

\beq
f (x) = \ao (x + \Eom ) - \aom, \qquad g(x) = \at (x + \Etm ) - \atm.
\eeq
We know that for $x \leq c t^2$ that
\beq
f(x) = g(x) \left(1 + \O ( t / t_0 ) \right).
\eeq
Assume that $x \geq c t^2$.  Both $f$ and $g$ are bijections of some intervals $\I_{f/g} \to \J_{f/g}$ and for their inverses we have
\beq
f^{-1} (y) = g^{-1} (y) \left( 1 + \O ( N^{5 \delta } t / t_0 \right).
\eeq
We write
\beq
f(x) - g(x) = g ( g^{-1} ( f (x) ) ) - g ( f^{-1} ( f (x) )) = \frac{ g ( g^{-1} (  f(x) )) - g ( f^{-1} f(x) ) )}{ g^{-1} ( f(x) ) - f^{-1} ( f(x) ) } \left( g^{-1}( f(x) )- f^{-1} ( f (x) ) \right).
\eeq
For the quotient we have the bound
\beq
\left| \frac{ g ( g^{-1} (  f(x) )) - g ( f^{-1} f(x) ) )}{ g^{-1} ( f(x) ) - f^{-1} ( f(x) ) }  \right| \leq C N^{2 \delta},
\eeq
which is a result of the fact that
\beq
\left| \frac{ \d a_1}{ \d E } \right| \leq C N^{2 \delta_0}
\eeq
for $E \geq E_- + N^{-\delta_0} t^2$.
Therefore, we obtain
\beq \label{eqn:abds}
f(x) = g(x) \left( 1 + \O ( N^{7 \delta} t / t_0 ) \right).
\eeq
Now, we write
\begin{align}
b_1 (x + E_{-, 1} ) - b_2 (x + E_{-, 2} ) = \left( b_1 (x + E_{-, 1} ) - b_2 (E ) \right) + \left( b_2 (E) - b_2 ( x + E_{-, 2} ) \right)
\end{align}
where $E$ is chosen so that $a_2 (E) = a_1 ( x + E_{-, 1} )$.  Then, by the above bounds on $b_1 (a) -b_2(a)$ we get
\beq
\left| b_1 (x + E_{-, 1} ) - b_2 (E ) \right| \leq C t \frac{ \sqrt{x} N^{5 \delta_0} t}{t_0}.
\eeq
By \eqref{eqn:Eicompare}, we see that
\beq
\left| E - (x + E_{-, 2} ) \right| \leq \frac{ C x N^{5 \delta_0} t}{t_0}.
\eeq
Hence, since $| \del_E \rho_i (E+ E_{-, i} ) | \leq C E^{-1/2}$ (by \eqref{eqn:mitder}) we see that
\beq
b_2 (E) - b_2 ( x + E_{-, 2} ) \leq C t \sqrt{x} \frac{N^{ 5 \delta_0} t }{ t_0}.
\eeq
We have therefore proved the following.
\bel \label{lem:match1}
Let $\eps >0$ and $t$, $t_0$ as above.  For $0 \leq x \leq c N^{ - 2 \eps } t_0^2$ we have
\beq \label{eqn:rho12match}
\rho_{t, 1} (x + E_{-, 1} ) = \rho_{t, 2} (x + E_{-, 2} )\left( 1 + \O ( N^{\eps} t/t_0 ) \right)
\eeq
and
\beq \label{eqn:mmatch}
| \Re[ m_{t, 1} (x + E_{-, 1} ) - m_{t, 1} ( E_{-, 1} ) ] - \Re [ m_{t, 2} (x + E_{-, 2} )  - m_{t, 2} (E_{-, 2} ) ] | \leq C \frac{ x N^{\eps}}{t_0}.
\eeq
\eel
We also need an estimate similar to \eqref{eqn:mmatch} for $x \leq 0$.  That is, we have
\bel \label{lem:match2} Let $\eps >0$ and $t, t_0$ as above.  For $-c N^{-2 \eps} t_0^2 (t_1/t_0) \leq x \leq 0$ we have,
\beq
| \Re[ m_{t, 1} (x + E_{-, 1} ) - m_{t, 1} ( E_{-, 1} ) ] - \Re [ m_{t, 2} (x + E_{-, 2} )  - m_{t, 2} (E_{-, 2} ) ] | \leq C\frac{|x|^{1/2} (t_1)^{1/2} N^{\eps}}{t_0^{1/2} }.
\eeq
\eel
\proof Fixing a scale $\eta \leq N^{- 2 \eps} t_0^2$ we can estimate the quantity by
\begin{align}
&| \Re[ m_{t, 1} (x + E_{-, 1} ) - m_{t, 1} ( E_{-, 1} ) ] - \Re [ m_{t, 2} (x + E_{-, 2} )  - m_{t, 2} (E_{-, 2} ) ] | \notag\\
\leq & \left| \int_{E \geq \eta, E <-1/2} \left( \frac{1}{ E- x } - \frac{1}{E  } \right) \rho_{t, 1} (E+E_{-1} ) \right| \notag\\
+& \left| \int_{E \geq \eta, E<-1/2} \left( \frac{1}{ E- x} - \frac{1}{E } \right) \rho_{t, 2} (E+E_{-2} ) \right| \notag\\
+ &\left| \int_{0 \leq E \leq  \eta} \left( \rho_{t, 2} (E +E_{-, 2} ) - \rho_{t, 1} (E+ E_{-, 2} ) \right) \left( \frac{1}{ E} + \frac{1}{E+x} \right)\right|   =: A_1 + A_2 +A_3.
\end{align}
By the square root behaviour of the densities we have,
\beq
|A_1| + |A_2| \leq C \frac{|x|}{\eta^{1/2}}.
\eeq
For $A_3$ we use \eqref{eqn:rho12match} and find
\beq
|A_3| \leq C N^{\eps}  \frac{t}{t_0} \eta^{1/2}.
\eeq
We obtain the estimate by choosing $\eta = |x|^{1/2} (t_0/t_1)^{1/2}$.
\qed

\subsubsection{Self-consistent equation coefficients} \label{sec:shorttimeself}
Recall the definition of $\hatD_1$.  It is clear that on this domain, that $\Im [ \mfct ] \geq c \sqrt{ \kappa + \eta}$, and so we get
\beq
|V_i - \xi | \geq \Im [ \xi] \geq \eta + c t \sqrt{ \kappa + \eta}.
\eeq
In the set-up of Section \ref{sec:ll2} the analogous estimates to Lemma \ref{lem:dxi} hold, and so 
in general we have
\beq
|\xi| \leq C(t^2 + \kappa + \eta ) \leq C ( \kappa + \eta)
\eeq
where the second inequality holds due to the definition of $\hatD_1$.  Hence,
\beq
\frac{1}{N} \sum_i |g_i|^p \leq C \frac{ |\xi|^{1/2}}{ \Im [ \xi]^{p-1}} \leq C \frac{ \sqrt{ \kappa + \eta }}{ ( \eta + t \sqrt{ \kappa + \eta })^{p-1} } .
\eeq
We now consider $\hatD_2$. 
Since the estimates of Lemma \ref{lem:dxi} hold, we have that there is a $C >0$ so that if $\kappa \geq C \eta$, then $\Re[ \xi ] \leq - c \kappa$.    Hence, we have \eqref{eqn:lla1}.  This also proves that
\beq
|V_i - \xi | \geq c ( \kappa+ \eta),
\eeq
as well as
\beq
\frac{1}{N} \sum_i |g_i|^p \leq \frac{C}{ ( \kappa + \eta )^{p-3/2}},
\eeq
for $\kappa \geq C' \eta$.  On the other hand if $ \kappa \leq C' \eta$, then,
\beq
\frac{1}{N} \sum_i |g_i|^p \leq C \frac{ |\xi|^{1/2}}{ \eta^{p-1} } \leq C \frac{ 1}{ ( \kappa + \eta )^{p-3/2}}.
\eeq

\appendix

\section{Large deviations estimates}

Let $X_i$ be a family of independent random variables obeying
\beq \label{eqn:momentassump}
\ee[X_i] = 0, \qquad \ee[ |X_i|^2] = 1, \qquad \ee[ |X_i|^p ] \leq C_p, \quad p \geq 2.
\eeq

We have the following estimates, see, e.g., \cite{general}.
\bel
Let $X_i$ and $Y_i$ be random variables obeying \eqref{eqn:momentassump}.  Let $b_i$ and $a_{ij}$ be deterministic.  We have for any $\eps >0$ and $D>0$ and $N$ large enough,
\begin{align}
\pp \left[ \left| \sum_i b_i X_i \right| \geq N^{\eps} \left( \sum_i |b_i|^2 \right)^{1/2} \right] & \leq N^{-D}, \\
\pp \left[ \left| \sum_{ij} a_{ij} X_i Y_i \right| \geq N^{\eps} \left( \sum_{ij} |a_{ij} |^2 \right)^{1/2} \right] & \leq N^{-D} , \\
\pp \left[ \left| \sum_{i \neq j } X_i a_{ij} X_j \right| \geq N^{\eps} \left( \sum_{i \neq j } |a_{ij}|^2 \right)^{1/2} \right] & \leq N^{-D}.
\end{align}
\eel

\section{Fluctuation averaging lemma}

We record here the following fluctuation averaging lemma.  As it is very similar to estimates appearing in \cite{landonyau}, we do not give a proof.  The proof of the main estimate, \eqref{eqn:momentbd}, is very similar to the proof given there.
\bep
Suppose that $\gamma$ is a deterministic parameter such that
\beq
| m_N - \mfct | \leq \gamma
\eeq
with overwhelming probability, and
\beq
\frac{1}{ N \eta } \leq \gamma \leq \frac{ t + \sqrt{ \kappa + \eta } + N^{-1/3} }{ \log(N)^2}.
\eeq
Then, with overwhelming probability, for any $\eps >0$,
\beq
|Q_i [ (G^{(\mathbb{T})}_{ii})^{-1} ] | \leq N^{\eps}\left( \sqrt{ \frac{t}{N} } + t \sqrt{ \frac{ \gamma + \Im [ \mfct ] }{ N \eta } } \right),
\eeq
and
\beq
\frac{1}{2} |g_i | \leq |G^{(\mathbb{T})}_{ii} | \leq 2 |g_i|,
\eeq
and,
\beq
|G_{ij} | \leq N^{\eps} |g_i | |g_j | \left( \sqrt{ \frac{t}{N} } + t \sqrt{ \frac{ \Im [ \mfct ] + \gamma }{ N \eta } } \right)
\eeq
uniformly over $|\mathbb{T} | \leq C$ for any fixed $C>0$.

Moreover, for any even $p>0$ we have, for any $\eps >0$, 
\begin{align} \label{eqn:momentbd}
&\ee \left| \frac{1}{N} \sum_i g_i^2 Q_i [ G_{ii}^{-1} ] \right|^p  \notag\\
\leq &N^{\eps} \max_{0 \leq s \leq p } \max_{ 0 \leq l \leq (p+s)/2} \left( \sqrt{ \frac{ t}{N} } + t \sqrt{ \frac{ \Im [ \mfct ] + \gamma }{ N \eta } } \right)^{p+s} \left( \sup_i |g_i | \right)^{s+2p-2l} \frac{1}{N^{p-l}} \left( \frac{1}{N} \sum_i |g_i|^2 \right)^l.
\end{align}

\eep

\section{$\Im [ m ]$ analysis}
Let $\mu_V$ be a measure whose Stieltjes transform obeys the assumptions of Section \ref{sec:main}.  Define the domain $\D_*$ by 
\begin{align}
\D_* := &\left\{ E + \i \eta : -3/4 \leq E \leq 0, 2 \eta_* \leq \eta \leq 10 \right\}\notag\\
& \cup \left\{ E + \i \eta : 0 \leq E \leq 3/4 , \eta_*^{1/2} \sqrt{ |E| } + \eta_* \leq \eta \leq 10 \right\} \cup \left\{ E + \i \eta : -3/4 \leq E \leq -2 \eta_*,  0 \leq \eta \leq 10 \right\}.
\end{align}
First it is clear that the estimates of Definition \ref{def:main} hold in the domain $\D_*$.  We want to prove the following.
\bel \label{lem:imm}
Let $\mu_V$ be as above.  For any $p \geq 2$ we have the following for $ a + b \i \in \D_*$.  If $a \geq 0$,
\beq \label{eqn:imm1}
\int \frac{ \d \mu_V (x)}{ |x- a - b \i |^p } \asymp \frac{ \sqrt{ a+b}}{ b^{p-1}}.
\eeq
If $a \leq 0$ then,
\beq \label{eqn:imm2}
\int \frac{ \d \mu_V (x) }{ |x - a -b \i |^p } \asymp \frac{1}{ ( |a| + b )^{p-3/2}}.
\eeq
\eel
\proof The upper bounds are immediate.  We first prove the lower bound of \eqref{eqn:imm1}.  Fix a $C_* >0$.  We have,
\begin{align}
\int \frac{ \d \mu_V (x) }{ |x - a - b \i |^p } &\geq \int_{ |x-a| \leq C_* b } \frac{ \d \mu_V (x) } { |x- a - b \i |^p} \geq \frac{c}{b^{p-2}} \int_{ |x-a |\leq C_* b} \frac{ \d \mu_V (x) }{ |x- a - b \i |^2}
\end{align}
for a $c>0$ depending on $C_*$.  We then have,
\begin{align}
\int_{ |x-a |\leq C_* b} \frac{ \d \mu_V (x) }{ |x- a - b \i |^2} &= \int \frac{ \d \mu_V (x) }{ |x- a - b \i |^2} -\int_{ |x-a |> C_* b} \frac{ \d \mu_V (x) }{ |x- a - b \i |^2} \notag\\
&\geq \frac{1}{b} \left( \Im [ m_V ( a + b \i ) ] - \frac{4}{C_*} \Im [ m_V (a + C_*b/2\i ) ] \right) \notag\\
&\geq \frac{1}{b} \left( c_1 \sqrt{ a+ b } - \frac{C_1 \sqrt{ a + C_* b }}{C_*} \right)
\end{align}
where $c_1$ and $C_1$ only depend on the assumptions on $\mu_V$.  Hence by taking $C_*$ large enough depending only on $c_1$ and $C_1$, we see that
\beq
\int_{ |x-a |\leq C_* b} \frac{ \d \mu_V (x) }{ |x- a - b \i |^2}  \geq c \frac{ \sqrt{ a + b }}{b}
\eeq
which in view of the above yields the lower bound of \eqref{eqn:imm1}.

The above argument also gives the lower bound of \eqref{eqn:imm2} in the regime $- b \leq a \leq 0$.  Note that in this case, the RHS of \eqref{eqn:imm2} is the same order as \eqref{eqn:imm1}.  In particular, we obtain \eqref{eqn:imm2} for all $a + b \i \in \D_*$ such that $- 2 \eta_* \leq a \leq 0$ (for such $a$ we have $b \geq 2 \eta_*$ by the definition of $\D_*$).

It remains to prove the lower bound of \eqref{eqn:imm2} in the case that $a \leq - 2 \eta_*$ and $|a| \geq b$.  Fix again a $C_*>0$.  We have
\beq
\int \frac{ \d \mu_V}{ |x-a-b\i |^p } \geq \int_{ x \leq C_* |a| } \frac{ \d \mu_V}{ |x-a - b \i |^p} \geq \frac{c}{ |a|^{p-2}} \int_{ x \leq C_* |a| } \frac{ \d \mu_V}{ (x-a)^2}
\eeq
for a $c>0$ depending on $C_*$.  We also used the fact that $|a|\geq b$ to observe that $|x-a-b\i|\geq |x-a|$ on the support of $\mu_V$ since $|a| \geq 2 \eta_*$.  We now have,
\begin{align}
\int_{ x \leq C_* |a| } \frac{ \d \mu_V (x)}{ (x-a)^2} &=\int \frac{ \d \mu_V (x)}{ (x-a)^2}  - \int_{ x \geq C_* |a| } \frac{ \d \mu_V (x)}{ (x-a)^2}  \notag\\
&\geq \int \frac{ \d \mu_V (x)}{ (x-a)^2}  - 2 \int \frac{ \d \mu_V (x)}{ (x- C_*a)^2} \notag\\
&\geq \frac{c}{ |a|^{1/2} } - \frac{C}{ |C_* a |^{1/2}}.
\end{align}
Choosing $C_*$ large enough yields the claim. \qed

\section{Free convolution continuity} \label{a:fccont}
In this section we consider two measures $\mu_1$ and $\mu_2$ and denote the free convolution of each with the semicircle by $m_{1, t}$ and $m_{2, t}$.  We estimate the difference $m_{1, t} - m_{2, t}$ under the assumption that $m_1 -m_2$ is small.

We assume that the restriction of $\mu_2$ to $[-1, 1]$ has a density $\rho_2 (x)$ such that
\beq
\rho_2 (x) \asymp \1_{ \{ x \geq 0 \} } \sqrt{x}.
\eeq

We assume the following estimates hold for  any $\delta, \eps >0$.  For any $1 \geq E \geq 0$ and $\sqrt{ |E| + \eta } \geq N^{\delta}/(N \eta ) + N^{\delta}/N^{1/3}$ we have,
\beq
|m_1 (z)  - m_2 (z) | \leq \frac{ N^{\eps}}{N \eta}.
\eeq
For $-1 \leq E \leq 0$ and $\eta \geq N^{\delta} / N^{2/3}$ we have
\beq
|m_1 (z) - m_2 (z) | \leq N^{\eps} \left( \frac{1}{ N (|E| + \eta ) } + \frac{1}{ ( N \eta )^2 \sqrt{ |E| + \eta } } \right).
\eeq

We denote $\xi_i = z + t m_{i, t}$ and we let $E_-$ be the edge of $\rho_{2, t}$, and define $\kappa = |E-E_-|$.  

\bel
Let $\mu_1$ and $\mu_2$ as above.  Let $\D_\sigma$ be as in Section \ref{sec:ll1} and $t$ satisfy
\beq
\frac{N^{\om}}{N^{1/3}} \leq t \leq N^{-\om}
\eeq
for $\om>0$.  For any $\eps >0$ the following estimates hold on $\D_\sigma$.  First, for $E \geq E_-$ we have
\beq
| m_{1, t} - m_{2, t} | \leq \frac{N^\eps}{N \eta}.
\eeq
For $E \leq E_-$ we have,
\beq
|m_{1, t} - m_{2, t} | \leq N^{\eps} \left( \frac{1}{ N ( \kappa + \eta ) } + \frac{1}{ ( N \eta )^2 \sqrt{ \kappa + \eta } }\right)
\eeq
\eel
\proof  Define $\Lambda = |m_{1, t} - m_{2, t} |$.  For $\eta$ order $1$, the estimates on $\Lambda$ follow easily.  First let us consider $E \geq E_-$.  Suppose that the estimate
\beq
\Lambda \leq \frac{ t + \Im [ m_{2, t} ] }{ \log(N)^2} 
\eeq
holds.  This assumption assures that $|x- \xi_2 | \gg |\xi_1 - \xi_2 |$ for $x$ in the support of $\mu_2$.  We then write
\beq
m_{1, t} - m_{2, t} = \left( \int \frac{ \d \mu_1}{ x- \xi_1 } - \int \frac{ \d \mu_2}{ x - \xi_1} \right) + \left( \int \frac{ \d \mu_2 }{ x - \xi_1 } - \int \frac{ \d \mu_2}{ x - \xi_2} \right).
\eeq
Expanding the second term and estimating the first by $N^{\eps} / N ( \Im [ \xi_1])$ leads to
\beq \label{eqn:detself}
\left| (1 - t R_2 ) ( m_{1, t} - m_{2, t} ) + t^2 R_3 (m_{1, t} - m_{2, t} )^2 \right| \leq \frac{N^{\eps}}{N \Im [ \xi_1] } + C t^3 \Lambda^3 \frac{ t + \Im [ m_{2, t} ] }{ (t^2 + \Im [ \xi_2 ] )^3}.
\eeq
Since $|1 - t R_2 | \asymp 1$ for $\kappa + \eta \geq t^2$ we can conclude that 
\beq
|m_{1, t} - m_{2, t} | \leq \frac{ N^{\eps}}{N \eta}, \qquad E \geq E_-, \quad \kappa + \eta \geq t^2.
\eeq
We can now suppose that $E_- \leq E \leq E_- + t^2$ and $\eta \leq t^2$.  Suppose that
\beq
\Lambda \leq \frac{ \sqrt{ \kappa + \eta}}{\log(N)^2}.
\eeq
Note that when $\eta = t^2$ we know that this is the case.  Then in this case, \eqref{eqn:detself} leads to 
\beq
\Lambda \leq C \frac{ \Lambda^2}{ \sqrt{ \kappa + \eta } } + \frac{N^{\eps}}{N \eta}.
\eeq
Since $(N \eta ) \ll \sqrt{ \kappa + \eta}$ for $E \geq E_-$, we conclude that $\Lambda \leq N^{\eps} / (N \eta)$ in the regime $ E \geq E_-$.

Now we consider the regime $E \leq E_-$.  
First, we observe that the estimate $\Lambda \leq N^{\eps} / (N \eta)$ in the regime $E_- - t^2 \leq E \leq E_-$ and $\eta \geq  ct^2$ follows from the above argument.  We then check the regime $E_- - t^2 \leq E \leq E_-$ and $\eta \leq c t^2$.  If we take $c >0$ small enough, we can ensure that
\beq
\Re[ \xi_2 ] \leq - c_1 t^2
\eeq
for another $c_1>0$.  Hence, if $\Lambda \leq \sqrt{ \kappa + \eta } / \log(N)^2$ we see that
\beq
| m_1 ( \xi_1 ) - m_2 ( \xi_1 ) |  \leq N^{\eps} \left( \frac{1}{ N t^2} + \frac{1}{ (N \eta)^2 t } \right).
\eeq
and so we obtain by a similar argument to above that the estimate $\Lambda \leq \sqrt{ \kappa  + \eta} / (\log(N))^2$ implies that
\beq
\left| (1 - t R_2 ) ( m_{1, t} - m_{2, t} ) + t^2 R_3 (m_{1, t} - m_{2, t} )^2 \right| \leq N^{\eps} \left( \frac{1}{ N t^2} + \frac{1}{ (N \eta)^2 t } \right)+ C \frac{\Lambda^3}{t^2}.
\eeq
We have that $|1- t R_2 | \asymp \sqrt{\kappa+\eta}/t$, and so we see that
\beq
\Lambda \leq C \frac{ \Lambda^2}{\sqrt{ \kappa + \eta}} + C N^{\eps} \left( \frac{1}{ N ( \kappa + \eta ) } + \frac{1}{ ( N \eta)^2 \sqrt{ \kappa + \eta } } \right).
\eeq
The second term is $\ll \sqrt{ \kappa + \eta}$ and so we see that so far we have proven the desired estimates as long as $E \geq E_- - t^2$.

Finally, to do the regime $E \leq E_- - t^2$ we first observe that there is a $C>0$ so that if $\kappa \geq C \eta$ then $- \Re[ \xi_2 ] \geq c \kappa$.  So, we see that the estimate $\Lambda \leq t / \log(N)^2$ implies
\beq
| m_{1, t} (\xi_1 )- m_{2, t} (\xi_1) | \leq N^{\eps}  \left( \frac{1}{ N ( \kappa + \eta ) } + \frac{1}{ (N \eta)^2 \sqrt{ \kappa + \eta } } \right).
\eeq
Hence, arguing as above we see that the estimate $\Lambda \leq t/\log(N)^2$ implies that
\beq
\Lambda \leq C \frac{ \Lambda^2}{t} + N^{\eps}  \left( \frac{1}{ N ( \kappa + \eta ) } + \frac{1}{ (N \eta)^2 \sqrt{ \kappa + \eta } } \right).
\eeq
This is enough to complete the proof. \qed

\section{Interpolating convolution measure properties} \label{a:inter}
We prove Lemma \ref{lem:edges}.   This follows from the following general estimate.  If we have two meaures $\mu_1, \mu_2$ that have, when restricted to $[-1, 1]$ a density $\rho (x)$ that behaves like
\beq
\rho (x) \asymp \sqrt{x} \1_{ \left\{ x \geq \right\} },
\eeq
and moreover
\beq
|m_1 (z) -m_2 (z) | \leq \frac{ N^{\eps}}{N \eta}
\eeq
for any $\eps >0$ and $\eta \geq N^{-2/3+\sigma}$.   Denote the $\xi$ maps at time $t$ by $\xi_1$ and $\xi_2$, the edges $E_1$, $E_2$, and $\xi_{i, -} = \xi_i (E_i)$.

Subtracting the defining equations for the $\xi_{i, -}$ we easily see
\beq
| \xi_{1, -} - \xi_{2, - } | \leq C t^3.
\eeq
Next, we estimate
\beq
|E_1 - E_2 | \leq | \xi_{1, -} - \xi_{2, -} | + t |m_1 ( \xi_{1, - } ) - m_2 ( \xi_{1, - } ) | + t |m_2 ( \xi_{2, - } - m_2 ( \xi_{1, -} ) |.
\eeq
The first term is bounded by $C t^3$ and since $| m_2' (E)| \leq t^{-1}$ for all $E$ such that $- E \asymp t^2$, we see that
\beq
t |m_2 ( \xi_{2, - } - m_2 ( \xi_{1, -} ) | \leq C t^3.
\eeq
For the last term, since the measures $\mu_1 = \mu_2$ on $[-1, 1]$ we see that,
\beq
| [m_1 ( \xi_{1, - } ) - m_2 ( \xi_{1, - } )] -[ m_1 ( \xi_{1, - } + \i N^{-1/2} ) - m_2 ( \xi_{1, - }  + \i N^{-1/2} )  ]| \leq N^{-1/2}
\eeq
and by assumption,
\beq
| m_1 ( \xi_{1, - } + \i N^{-1/2} ) - m_2 ( \xi_{1, - }  + \i N^{-1/2} )  | \leq \frac{ N^{\eps}}{N^{1/2}}.
\eeq
Hence, 
\beq
|E_1 - E_2 | \leq N^{\eps} (t^3 + \frac{t}{N^{1/2}}).
\eeq


\bibliography{mybib}{}
\bibliographystyle{abbrv}

\end{document}